\documentclass[11pt]{article}

\usepackage[a4paper,margin=2.2cm]{geometry}

\usepackage{inputenc}
\usepackage{amsmath}%
\usepackage{amsfonts}%
\usepackage{amssymb}%
\usepackage{graphicx}
\usepackage{mathtools}
\usepackage[T2A]{fontenc}
\usepackage{mathrsfs}
\usepackage{tikz}
\usetikzlibrary{arrows.meta}

\usepackage{subcaption}

\usepackage{xcolor}
\usepackage{amsthm} 
\usepackage[colorlinks]{hyperref}
\hypersetup{
    citecolor=magenta,
    colorlinks=true,
    linkcolor=blue,
    filecolor=magenta,      
    urlcolor=cyan,    
    }
\usepackage{enumerate}

\usepackage[all,cmtip]{xy}

\usepackage{doi}

\newcommand{\bbeta}{\boldsymbol{\beta}}
\newcommand{\blambda}{\boldsymbol{\lambda}}
\newcommand{\btheta}{{\boldsymbol{\theta}}}
\newcommand{\bw}{\mathbf{w}}

\newcommand{\bbg}{\bar{\mathbf{\Gamma}}}
\newcommand{\bg}{\mathbf{\Gamma}}
\newcommand{\std}{\mathrm{std}}
\newcommand{\op}{\mathrm{op}}

\DeclareMathOperator{\Ad}{Ad}
\DeclareMathOperator{\gl}{\mathfrak{gl}}
\DeclareMathOperator{\Mat}{Mat}
\DeclareMathOperator{\SL}{SL}
\DeclareMathOperator{\GL}{GL}
\DeclareMathOperator{\dom}{dom}
\DeclareMathOperator{\pr}{pr}
\DeclareMathOperator{\tr}{tr}
\DeclareMathOperator{\Span}{Span}
\DeclareMathOperator{\CYB}{CYB}

\DeclareMathOperator{\Imm}{Im}

\DeclareMathOperator{\Pf}{Pf}

\DeclareMathOperator{\ses}{ss}

\DeclareMathOperator{\id}{id}
\DeclareMathOperator{\sm}{sm}

\DeclareMathOperator{\sgn}{sgn}

\newcommand{\bdvsize}{1.2pt} 
\newtheorem{theorem}{Theorem}[section]
\newtheorem{lemma}[theorem]{Lemma}
\newtheorem{proposition}[theorem]{Proposition}
\newtheorem*{proposition*}{Proposition}
\newtheorem{corollary}[theorem]{Corollary}
\theoremstyle{definition}\newtheorem{definition}[theorem]{Definition}
\theoremstyle{definition}\newtheorem{remark}[theorem]{Remark}
\newtheorem{example}[theorem]{Example}

\theoremstyle{remark}
\numberwithin{equation}{section}

\newcommand{\edge}{\,
	\begin{tikzpicture}
		\draw (0,.5ex) -- (2ex,.5ex);
		\draw[fill] (0,.5ex) circle (.2ex);
		\draw[fill] (2ex,.5ex) circle (.2ex);
		\draw[fill] (0,0) circle (0);
	\end{tikzpicture}\,
}

\newcommand{\rightedge}{\,
	\begin{tikzpicture}
		\draw[-{Stealth}, midway] (0,.5ex) -- (2ex,.5ex);
		\draw[fill] (0,.5ex) circle (.2ex);
		\draw[fill] (2ex,.5ex) circle (.2ex);
		\draw[fill] (0,0) circle (0);
	\end{tikzpicture}\,
}

\title{Belavin--Drinfeld bivectors and Poisson deformations}
\author{Zihang Liu and Dmitriy Voloshyn}
\date{}

\newcounter{myfootnote}
\begin{document}
\maketitle

\newcommand{\blfootnote}[1]{%
    \begingroup
    \begin{NoHyper}
        \renewcommand{\thefootnote}{}%
        \footnote{#1}%
        \addtocounter{footnote}{-1}%
    \end{NoHyper}
    \endgroup
}

\begin{abstract}
    We study Belavin--Drinfeld Poisson bivectors $\pi_{(\bbg^r,\bbg^c)}$ on $\Mat_{n\times m}(\mathbb{C})$ from the perspective of Poisson deformation theory. We prove that the log-canonical part $\pi_0$ of $\pi_{(\bbg^r,\bbg^c)}$ is $\mathbb{T}_{(\bg^r,\bg^c)}$-log-symplectic if and only if the pair $(\bg^r,\bg^c)$ of Belavin--Drinfeld triples is aperiodic. In this case, we show that $\pi_{(\bbg^r,\bbg^c)}$ is a maximal $\mathbb{T}_{(\bg^r,\bg^c)}$-invariant Poisson deformation of $\pi_0$, and we construct a  $\mathbb{T}_{(\bg^r,\bg^c)}$-equivariant Poisson birational map between $(\Mat_{n\times m},\pi_0)$ and $(\Mat_{n\times m},\pi_{(\bbg^r,\bbg^c)})$.
\end{abstract}

\tableofcontents


\section{Introduction and main results}
Recently, Jiang-Hua Lu and Mykola Matviichuk developed a theory of algebraic Poisson deformations of $\mathbb{T}$-log-symplectic Poisson bivectors \cite{lm-deform}. The theory starts with the complex affine space $\mathsf{X} := \mathbb{C}^n$ equipped with coordinates $x_1,\ldots,x_n$, an action by a complex torus $\mathbb{T}$ such that each $x_i$ is a $\mathbb{T}$-weight vector, and a log-canonical bivector
\begin{equation*}
    \pi_0 := \sum_{1 \leq j<k \leq n} \lambda_{jk} x_j x_k \frac{\partial}{\partial x_j} \wedge \frac{\partial}{\partial x_k},\quad 
    \lambda_{jk} \in \mathbb{C}.
\end{equation*}
The bivector $\pi_0$ is called \emph{$\mathbb{T}$-log-symplectic} if there is an open $\mathbb{T}$-leaf (that is, $\mathbb{T}\cdot S$ is open for some symplectic leaf $S$). Given a bivector $\pi_1$ with no log-canonical terms representing a class in the algebraic $\mathbb{T}$-invariant Poisson cohomology $ \mathsf{H}_{\pi_0}^2(\mathsf{X})^{\mathbb{T}}$, the Lu--Matviichuk theory describes $\mathbb{T}$-invariant bivectors $\pi^\prime$ as correction terms such that the bivector $\pi:=\pi_0 + \pi_1 + \pi^\prime$ is a (possibly formal) $\mathbb{T}$-invariant Poisson bivector. The Poisson bivector $\pi$ is called a \emph{Poisson deformation of $\pi_0$ along~$\pi_1$}.

Without any further assumptions, a Poisson deformation of $\pi_0$ along some $\pi_1$ exists at least in the formal sense. One then distinguishes the class of bivectors $\pi_1$ satisfying Property $(\mathrm{W0})$, and the subclass satisfying Property~\eqref{eq:W1}. Roughly, Property $(\mathrm{W0})$ guarantees that the resulting Poisson deformation is algebraic rather than formal. Property~\eqref{eq:W1}, central to this paper, further guarantees that a Poisson deformation along a given $\pi_1$ is unique.

The theory was motivated by earlier developments in logarithmic geometry, where $\pi_0$ was assumed to be a log-symplectic bivector with normal crossings Pfaffian divisor on a variety. M.~Matviichuk, B.~Pym, and T.~Schedler proved in \cite{mps-torelli} that such a bivector $\pi_0$ admits formal Poisson deformations. By considering $\mathbb{T}$-log-symplectic bivectors, Lu and Matviichuk were able to include examples from Poisson--Lie theory and cluster algebras in the Poisson deformation theory. 

One class of such examples is given by Poisson bivectors \(\pi\) arising from symmetric Poisson CGL extensions. As proved in \cite[Theorem 11.1]{gy-cgl} by K. R. Goodearl and M. Yakimov, if $\mathsf{X}$ is equipped with such a $\pi$, the coordinate ring $\mathbb{C}[\mathsf{X}]$ admits the structure of a cluster algebra compatible with $\pi$ in the sense of~\cite{gsv-book}. The induced standard Poisson structures on (generalized) Schubert cells~\cite{elek-lu} and the space of rectangular matrices~\cite{gy-cgl} are examples of symmetric Poisson CGL extensions. One of the remarkable achievements of the Lu--Matviichuk theory is a complete classification of symmetric Poisson CGL extensions in terms of their log-canonical parts \cite[Theorem 4.4.3]{lm-deform}.

The purpose of this paper is to relate the Lu--Matviichuk theory to another line of research known as the GSV conjecture. The conjecture was formulated by M. Gekhtman, M. Shapiro, and A.~Vainshtein in~\cite{gsv-conj}. It asserts that, given a simply connected simple complex algebraic group $G$ and a Belavin--Drinfeld Poisson bivector $\pi$, the coordinate ring $\mathbb{C}[G]$ has the structure of a (generalized) upper cluster algebra compatible with $\pi$. The GSV conjecture extends to varieties related to $G$, such as the dual Poisson--Lie group~\cite{multdual} and the Drinfeld double~\cite{multdouble} (see~\cite{gsv-unified,multdual} for an overview of the state of the conjecture).

The Belavin--Drinfeld bivectors $\pi_{(\bbg^r,\bbg^c)}$ on $G$ are the Poisson bivectors parameterized by pairs $(\bbg^r,\bbg^c)$ of Belavin--Drinfeld quadruples (BD quadruples). A BD quadruple is a pair $\bbg:=(\bg,r_0)$ where $\bg:=(\Gamma_1,\Gamma_2,\gamma)$ is a BD triple, and $r_0$ is a certain $2$-tensor with components in a Cartan subalgebra $\mathfrak{h}$ (see Section~\ref{s:bd_background} for details), with a skew-symmetric part denoted by $s_0$. The sets $\Gamma_1$ and $\Gamma_2$ are subsets of the set of simple roots of~$G$, and $\gamma : \Gamma_1 \rightarrow \Gamma_2$ is a nilpotent isometry. The pair $(\bbg^r,\bbg^c)$ gives rise to the BD torus $\mathbb{T}_{(\bg^r,\bg^c)}$ acting on~$G$. In this setting, the standard Poisson--Lie bivector corresponds to the pair of trivial BD quadruples, where $\Gamma_1$ and $\Gamma_2$ are empty, and $r_0$ is half of the Casimir element of~$\mathfrak{h}$.

For the group $G$, one distinguishes between aperiodic and non-aperiodic pairs $(\bg^r,\bg^c)$ (see Definition~\ref{d:aperiodic}). The notion depends only on the BD triples and not on $(r_0^r,r_0^c)$. As part of the GSV conjecture, it is expected that $(G,\pi_{(\bbg^r,\bbg^c)})$ has an open $\mathbb{T}_{(\bg^r,\bg^c)}$-leaf if and only if $(\bg^r,\bg^c)$ is aperiodic, and in this case, $\mathbb{C}[G]$ is expected to have the structure of a (non-generalized) upper cluster algebra compatible with the BD bivector $\pi_{(\bbg^r,\bbg^c)}$. If $(\bg^r,\bg^c)$ is not aperiodic, $\mathbb{C}[G]$ is expected to have the structure of a generalized upper cluster algebra compatible with the BD bivector. The Poisson deformation theory gives a clue on how to generalize aperiodicity of BD pairs to $\Mat_{n\times m}(\mathbb{C})$. We analyze this in Section~\ref{s:meander}.

We establish in this paper that the Belavin--Drinfeld bivectors on $\Mat_{n\times m}(\mathbb{C})$ associated with aperiodic BD pairs are among the nicest examples in the Lu--Matviichuk Poisson deformation theory. They are completely determined by their log-canonical terms, and they have Property~\eqref{eq:W1}. In technical terms, we prove the following theorem:

\begin{minipage}{\textwidth}
\begin{theorem}\label{theorem}
    For any Belavin--Drinfeld quadruples $\bbg^r := (\bg^r,r_0^r)$ and $\bbg^c := (\bg^c,r_0^c)$, the log-canonical part $(\pi_{(\bbg^r,\bbg^c)})_0$ of the Poisson bivector $\pi_{(\bbg^r,\bbg^c)}$ on $\Mat_{n\times m}(\mathbb{C})$ is given by
    \begin{equation*}
        (\pi_{(\bbg^r,\bbg^c)})_0 = (\pi_{\std})_0 + (s_0^c)^L - (s_0^r)^R.
    \end{equation*}
    The Poisson bivector $(\pi_{(\bbg^r,\bbg^c)})_0$ is  $\mathbb{T}_{(\bg^r,\bg^c)}$-log-symplectic if and only if the pair $(\bg^r,\bg^c)$ is aperiodic. If this is the case, the following statements hold:
    \begin{enumerate}[1)]
        \item The pair $\big((\pi_{(\bbg^r,\bbg^c)})_0,\mathbb{T}_{(\bg^r,\bg^c)}\big)$ satisfies Property~\eqref{eq:W1}, and $\pi_{(\bbg^r,\bbg^c)}$ is the maximal normalized admissible algebraic $\mathbb{T}_{(\bg^r,\bg^c)}$-invariant Poisson deformation of $(\pi_{(\bbg^r,\bbg^c)})_0$;
        \item There exists a $\mathbb{T}_{(\bg^r,\bg^c)}$-equivariant Poisson birational map
        \begin{equation*}
            \mathcal{U}_0 : (\Mat_{n \times m}(\mathbb{C}),(\pi_{(\bbg^r,\bbg^c)})_0) \dashrightarrow (\Mat_{n \times m}(\mathbb{C}),\pi_{(\bbg^r,\bbg^c)});
        \end{equation*}
        \item The Poisson variety $(\Mat_{n \times m}(\mathbb{C}),\pi_{(\bbg^r,\bbg^c)})$ has an open $\mathbb{T}_{(\bg^r,\bg^c)}$-leaf.
    \end{enumerate}
\end{theorem}
\end{minipage}\newline

On the way to proving Theorem~\ref{theorem}, we establish several results of independent interest. The key features of Belavin--Drinfeld Poisson brackets are:
\begin{enumerate}[1)]
    \item The $\mathbb{T}_{(\bg^r,\bg^c)}$-weights $\beta_{kl}$ of the matrix coordinates $x_{kl}$ are pairwise distinct, and the torus $\mathbb{T}_{(\bg^r,\bg^c)}$ has a scaling subtorus;\label{i:cond1}
    \item The kernel of $\bbeta:=(\beta_{kl} \mid k \in [1,n], \ l \in [1,m])$ has a basis $\mathcal{S}^\prime$ given in terms of nonzero $(\mathbb{C}^\times)^{nm}$-weights in $\mathsf{H}_{\pi_0}^2(\Mat_{n \times m}(\mathbb{C}))^{\mathbb{T}_{(\bg^r,\bg^c)}}$.\label{i:cond2}
\end{enumerate}
Reformulating the above properties in the abstract setting of the Lu--Matviichuk theory, in Section~\ref{s:pdeform} we derive several results for log-canonical bivectors $\pi_0$ that may not be $\mathbb{T}$-log-symplectic. In Proposition~\ref{p:w1check}, we show that condition~\ref{i:cond1} is a sufficient condition for Property~\eqref{eq:W1}, and in Proposition~\ref{prop:SisLM}, we derive a description of the nonzero $(\mathbb{C}^{\times})^{n}$-weights in $\mathsf{H}_{\pi_0}^2(\mathsf{X})^{\mathbb{T}}$. Under condition~\ref{i:cond2}, in Proposition~\ref{p:ker_cycles}, we derive an obstruction to $\mathbb{T}$-log-symplecticity. 

We believe that a theory of Poisson deformations of $\pi_0$ that are not $\mathbb{T}$-log-symplectic is tractable under the above two conditions, and that such a theory may give rise to generalized cluster structures.

The second part of Theorem~\ref{theorem} not only sharpens the statement that $\pi_{(\bbg^r,\bbg^c)}$ is a Poisson deformation of its log-canonical part, but is also relevant to the theory of cluster algebras. Indeed, since the matrix entries are log-canonical with respect to $(\pi_{(\bbg^r,\bbg^c)})_0$, their pullbacks $(\mathcal{U}_0^{-1})^*(x_{ij})$ are log-canonical with respect to $\pi_{(\bbg^r,\bbg^c)}$. Extracting their numerators, we thus obtain a regular log-canonical system on $\Mat_{n\times m}(\mathbb{C})$, which can be upgraded to an initial seed. The currently available techniques in the literature are sufficient to prove that $\mathbb{C}[\Mat_{n\times m}]$ has the structure of an upper cluster algebra compatible $\pi_{(\bbg^r,\bbg^c)}$ (for aperiodic BD pairs). However, in this paper, we focus only on the Poisson deformation side of the story.

The paper is organized as follows. In Section~\ref{s:prelims}, we describe  background material on the Lu--Matviichuk Poisson deformation theory and on the theory of Belavin--Drinfeld Poisson brackets on $\Mat_{n\times m}(\mathbb{C})$. In Section~\ref{s:pdeform}, we derive several new results in Poisson deformation theory that are relevant for BD bivectors. In Section~\ref{s:main}, we prove the first part of Theorem~\ref{theorem}, as well as analyze the aperiodicity condition for $(\bg^r,\bg^c)$ in the setting of $\Mat_{n\times m}(\mathbb{C})$. In Section~\ref{s:birat}, we prove the second and the third parts of Theorem~\ref{theorem}.

\paragraph{Tool and computational resource disclosure.} 
Lemma~\ref{l:coisop_slice} was suggested by GPT-5.6 Sol, which significantly simplified our initial approach in the proof of Proposition~\ref{p:u_poiss}. In turn, the lemma generalizes \cite[Proposition 6.7, Part 1]{lu-evens}. All other mathematical arguments were developed by the authors. Freely available AI tools were used for spellchecking the final draft. Numerical experiments were carried out using our own software developed in MATLAB.

\paragraph{Acknowledgments.} We thank Jiang-Hua Lu and Mykola Matviichuk for patiently explaining to us various parts of their paper~\cite{lm-deform}, as well as for numerous discussions. We thank Misha Gekhtman for his suggestion to use planar networks in the proof of Proposition~\ref{p:U0}. The research of Z. Liu has been supported by the RGC of the Hong Kong SAR, China (GRF 17304415). The research of D.~Voloshyn was supported by the Institute for Basic Science (IBS-R003-D1).

\section{Preliminaries}\label{s:prelims}
In this section, we describe background material on the Poisson geometry of $\Mat_{n\times m}(\mathbb{C})$ and the Lu--Matviichuk Poisson deformation theory. 

Our notation for the Poisson deformation theory is mostly the same as that in~\cite{lm-deform}. For a matrix $X$, we denote by $X^T$ its transpose, and for a subset $I$ of its row labels and a subset $J$ of its column labels, we denote by $X_{I}^J$ the submatrix of $X$ formed by $I\times J$ entries. For a $2$-tensor $a\otimes b$, we set $(a\otimes b)^t := b\otimes a$. The identity matrix of size $n\times n$ is denoted by $I_n$.

We use the following notation in Lie-theoretic arguments. We denote by $G$ a reductive complex algebraic group, and by $\mathfrak{g}$ its Lie algebra equipped with an invariant nondegenerate symmetric bilinear form $\langle\, ,\, \rangle$. The corresponding Casimir element is denoted by $\Omega$. We choose a pair of opposite Borel subgroups $\mathcal{B}$ and $\mathcal{B}_-$, with unipotent radicals $\mathcal{N}$ and $\mathcal{N}_-$, respectively, and we denote by $\mathcal{H}:=\mathcal{B}\cap \mathcal{B}_-$ the corresponding Cartan subgroup. We use calligraphic letters for the subgroups in order to reserve capital letters for matrices (for instance, we denote by $B$ an upper triangular matrix). The corresponding Lie algebras are denoted by $\mathfrak{b}$, $\mathfrak{b}_-$, $\mathfrak{n}$, $\mathfrak{n}_-$, and $\mathfrak{h}$. 

The maps $\pr_> : \mathfrak{g} \rightarrow \mathfrak{n}$, $\pr_< : \mathfrak{g} \rightarrow \mathfrak{n}_-$, and $\pr_0 : \mathfrak{g}\rightarrow \mathfrak{h}$ are the projections onto the corresponding subspaces. We denote by $\Omega_{>}:=(1\otimes \pr_>)\Omega \in \mathfrak{n}_-\otimes \mathfrak{n}$, $\Omega_{<} := (1\otimes \pr_<)\Omega \in \mathfrak{n}\otimes \mathfrak{n}_-$, and $\Omega_0:=(1\otimes \pr_0)\Omega \in S^2\mathfrak{h}$ the $2$-tensors corresponding to these projections.

The set of positive (resp. negative) roots is denoted by $\Phi_+$ (resp. $\Phi_-$), and the set of simple roots by $\Delta$. If the Cartan--Killing type of $\mathfrak{g}$ is specified, we denote the simple roots by the same symbol as the corresponding Dynkin diagram (for instance, for type $A_{n-1}$, we write $A_{n-1}$ in place of $\Delta$).

The Lie algebra $\gl_n(\mathbb{C})$ is equipped with the trace form: $\langle X,Y \rangle := \tr(XY)$. The space of matrices $\Mat_{n\times m}(\mathbb{C})$ has a natural action of $\GL_n(\mathbb{C})$ on the left (resp. $\GL_m(\mathbb{C})$ on the right). To distinguish the Lie-theoretic data associated with $\GL_n(\mathbb{C})$ and $\GL_m(\mathbb{C})$, we sometimes use the superscript $(n)$ and $(m)$, respectively. For instance, $\mathcal{B}_-^{(n)}$ denotes the subgroup of $\GL_n(\mathbb{C})$ consisting of nondegenerate lower triangular matrices, and $\Omega^{(m)}$ denotes the Casimir element of $\gl_m(\mathbb{C})$ associated with the trace form.

For a $2$-tensor $r \in \mathfrak{g}\otimes \mathfrak{g}$, we denote by $r^R$ (resp. $r^L$) the right-invariant (resp. left-invariant) contravariant tensor field whose value at the identity element of $G$ is $r$. 

In the case of $\Mat_{n\times m}(\mathbb{C})$, given a tensor $r \in \gl_m(\mathbb{C})^{\otimes k}$, we denote by $r^L$ the contravariant tensor field on $\Mat_{n\times m}(\mathbb{C})$ obtained by left translation of $r$, as follows:
\begin{equation*}
    (r^L)_X:=[(L_X)_* \otimes (L_X)_* \otimes \cdots \otimes (L_X)_*](r), \quad X \in \Mat_{n\times m}(\mathbb{C}).
\end{equation*}
Likewise, for $r \in \gl_n(\mathbb{C})^{\otimes k}$, we denote by $r^R$ the contravariant tensor field on $\Mat_{n\times m}(\mathbb{C})$ obtained by right translation of $r$, as follows:
\begin{equation*}
    (r^R)_X:=[(R_X)_* \otimes (R_X)_* \otimes \cdots \otimes (R_X)_*](r), \quad X \in \Mat_{n\times m}(\mathbb{C}).
\end{equation*}

Given a regular function $f$ on $G$, we define its left and right gradients as the maps
\begin{equation*}
    \nabla^L f,\ \nabla^R f : G \rightarrow \mathfrak{g}
\end{equation*}
given by
\begin{align*}
    \langle \nabla^L_gf,x\rangle  &= d_gf(x^L) = \left.\frac{d}{dt}\right|_{t=0} f(g\exp(tx));\\
    \langle \nabla^R_g f,x\rangle &= d_gf(x^R) = \left.\frac{d}{dt}\right|_{t=0} f(\exp(tx)g);
\end{align*}
where $g \in G$ and $x \in \mathfrak{g}$. Here, the subscript in $\nabla^L_g$ and $\nabla^R_g$ indicates that the gradient is evaluated at the point $g$. The gradients are related by the identity
\begin{equation*}
    \Ad_g \nabla^L_g f = \nabla_g^R f.
\end{equation*}
In particular, for any two regular functions $f_1$ and $f_2$ on $G$,
\begin{equation*}
    \langle \nabla^L_g f_1,\nabla^L_g f_2\rangle = \langle \nabla^R_g f_1,\nabla^R_g f_2\rangle.
\end{equation*}

For a regular function $f$ on $\Mat_{n\times m}(\mathbb{C})$, the left and right gradients
\begin{equation*}
    \nabla^L f : \Mat_{n \times m}(\mathbb{C}) \rightarrow \gl_m(\mathbb{C}), \quad \ \nabla^R f: \Mat_{n\times m}(\mathbb{C})\rightarrow \gl_n(\mathbb{C})
\end{equation*}
are given by similar formulas:
\begin{align*}
    \langle \nabla^L_Xf,Y^{(m)}\rangle  &= d_Xf((Y^{(m)})^L) = \left.\frac{d}{dt}\right|_{t=0} f(X \exp(t Y^{(m)}));\\
    \langle \nabla^R_X f,Y^{(n)}\rangle &= d_Xf((Y^{(n)})^R) = \left.\frac{d}{dt}\right|_{t=0} f(\exp(t Y^{(n)})X),
\end{align*}
where $X\in \Mat_{n\times m}(\mathbb{C})$,  $Y^{(n)} \in \gl_n(\mathbb{C})$, and $Y^{(m)} \in  \gl_m(\mathbb{C})$. Explicitly, the left and right gradients of $f$ are given by
\begin{equation*}
    \nabla^L_X f = \nabla_X f \cdot X \in \gl_m(\mathbb{C}), \quad \nabla^R_X f = X \cdot \nabla_X f \in \gl_n(\mathbb{C}),
\end{equation*}
where $\nabla_X f$ is the matrix gradient of $f$, whose entries are given by
\begin{equation*}
    (\nabla_X f)_{ij} := \frac{\partial f}{\partial x_{ji}}(X), \quad i \in [1,m], \ j \in [1,n].
\end{equation*}

\subsection{Belavin--Drinfeld bivectors}\label{s:bd_background}
In this section, we provide background on Poisson--Lie theory. We use the setup and notation introduced in the preamble to Section~\ref{s:prelims}.

For a Poisson bivector $\pi \in \Gamma(G,\Lambda^2 G)$, the pair $(G,\pi)$ is called a \emph{Poisson--Lie group} if the multiplication map of $G$ is Poisson. If $G$ is semisimple, then for any $\pi$ such that $(G,\pi)$ is a Poisson--Lie group, there exists $r \in \mathfrak{g}\otimes \mathfrak{g}$ such that $\pi$ is of the form
\begin{equation*}
    \pi = r^L - r^R.
\end{equation*}

We will be interested in Poisson bivectors on $G$ of the form
\begin{equation}\label{eq:biv}
    \pi = (r^c)^L - (r^r)^R
\end{equation}
where $r^r$ and $r^c$ are solutions of the \emph{Classical Yang--Baxter equation (CYBE)}, each with symmetric part equal to $(1/2)\Omega$. The superscripts $c$ and $r$ stand for \emph{column} and \emph{row}, and should be understood as labels\footnote{The notation is a historical tribute to~\cite{gsv-plethora}. In that paper, the authors constructed a family of compatible cluster structures on square matrices, and the construction of initial variables involved gluing of certain blocks along rows and columns. The gluing is governed by the Belavin--Drinfeld data associated with $r^r$ or $r^c$ (the first one for gluing along rows, the second one for columns). } distinguishing the two $r$-matrices. The CYBE is an equation for $r \in \mathfrak{g}\otimes \mathfrak{g}$ of the form
\begin{equation*}
    \CYB(r) := [r_{12},r_{13}] + [r_{12},r_{23}] + [r_{13},r_{23}] = 0.
\end{equation*}
Here, given $a,b \in \mathfrak{g}$, we set $(a\otimes b)_{12} := a\otimes b \otimes 1$, $(a\otimes b)_{13}:= a \otimes 1 \otimes b$, and so on. Defining a linear map $R:\mathfrak{g}\rightarrow \mathfrak{g}$ by\footnote{Equivalently, $R$ is such that $\langle r, x \otimes y \rangle = \langle R(x),y \rangle$ for $x, y \in \mathfrak{g}$.} $r=(1\otimes R)\Omega$, and denoting by $R^*$ its adjoint with respect to the form, we can express the CYBE for $r$ in the following equivalent form:
\begin{equation*}
    [R(x),R(y)] = R([R(x),y] - [x,R^*(y)])
\end{equation*}
where $x,y \in \mathfrak{g}$.

In~\cite{bd}, under the assumption that $\mathfrak{g}$ is semisimple, Belavin and Drinfeld classified solutions of the CYBE with a nondegenerate symmetric part (the so-called \emph{factorizable quasi-triangular $r$-matrices}). The classification is given in terms of Belavin--Drinfeld quadruples, defined as follows.

\begin{definition}
    Given the set of simple roots $\Delta$, a \emph{Belavin--Drinfeld triple} (a \emph{BD triple}) is a triple $\bg:=(\Gamma_1,\Gamma_2,\gamma)$ such that $\Gamma_1,\Gamma_2\subset \Delta$ and $\gamma : \Gamma_1 \rightarrow \Gamma_2$ is a nilpotent isometry. A \emph{Belavin--Drinfeld quadruple} (a \emph{BD quadruple}) is the tuple $\bbg:=(\bg,r_0)$ consisting of a BD triple $\bg$ and an element $r_0 \in \mathfrak{h}\otimes \mathfrak{h}$ satisfying the following equations:
    \begin{align}
        &r_0 + r_0^t = \Omega_{0}; \label{eq:r0_1} \\
        &(\gamma(\alpha) \otimes 1)(r_0) + (1 \otimes \alpha) (r_0) = 0 \in \mathfrak{h}, \ \ \alpha \in \Gamma_1.\label{eq:r0_2}
    \end{align}
    A BD triple $\bg$ is called \emph{trivial} (or \emph{standard}) if $\bg = \bg_{\std}:=(\emptyset,\emptyset,\emptyset \rightarrow \emptyset)$. A BD quadruple $\bbg$ is called \emph{trivial} (or \emph{standard}) if $\bbg = \bbg_{\std}:=(\bg_{\std},\frac{1}{2}\Omega_0)$.
\end{definition}
The nilpotency of $\gamma$ means that for each $\alpha \in \Gamma_1$, there exists $j \geq 1$ such that $\gamma^j(\alpha) \notin \Gamma_1$. We denote by $\gamma^* : \Gamma_2 \rightarrow \Gamma_1$ the adjoint of $\gamma$ (equivalently, its inverse). 

An element $r_0$ can be expressed as $r_0 = \frac{1}{2}\Omega_0 + s_0$, where $s_0\in \wedge^2 \mathfrak{h}$. By using the form on $\mathfrak{h}\otimes \mathfrak{h}$ corresponding to $\Omega_0$, we define the corresponding linear maps $R_0,S_0 : \mathfrak{h}\rightarrow \mathfrak{h}$ via 
\begin{equation*}
    r_0 = (1\otimes R_0)\Omega_0, \quad s_0 = (1\otimes S_0)\Omega_0.
\end{equation*}
In terms of $R_0$, equations~\eqref{eq:r0_1}--\eqref{eq:r0_2} translate as
\begin{align*}
    &R_0 + R_0^* = \id;\\
    &R_0(h_{\alpha}-h_{\gamma(\alpha)}) = h_{\alpha}, \quad \alpha \in \Gamma_1.
\end{align*}
Here, $R_0^*$ denotes the adjoint of $R_0$ with respect to the form on $\mathfrak{h}$, $\id$ denotes the identity map, and $h_{\alpha} \in \mathfrak{h}$ denotes the element corresponding to the simple root $\alpha$ under the form. Likewise, for the skew-adjoint map $S_0$, the equations are given by
\begin{equation}\label{eq:S0}
    S_0(h_{\alpha}-h_{\gamma(\alpha)}) = \frac{1}{2}(h_{\alpha} + h_{\gamma(\alpha)}), \quad \alpha \in \Gamma_1.
\end{equation}

For a given BD triple $(\Gamma_1,\Gamma_2,\gamma)$, the map $\gamma$ extends to an isomorphism $\gamma : \mathbb{Z}\Gamma_1 \rightarrow \mathbb{Z}\Gamma_2$  between lattices. Choosing Chevalley generators $e_{\alpha}$ and $e_{-\alpha}$ of $\mathfrak{g}$, $\alpha \in \Delta$, we set
\begin{equation*}
    \gamma(e_{\pm \alpha}) := \begin{cases}
        e_{\pm \gamma(\alpha)} & \alpha \in \Gamma_1;\\
        0 &\text{otherwise.}
    \end{cases}
\end{equation*}
We then extend $\gamma$ to unique Lie algebra homomorphisms $\gamma : \mathfrak{n}_-\rightarrow \mathfrak{n}_-$ and $\gamma : \mathfrak{n}\rightarrow \mathfrak{n}$. We denote by $\tilde{\gamma}: \mathcal{N}_- \rightarrow \mathcal{N}_-$ and $\tilde{\gamma} : \mathcal{N}\rightarrow \mathcal{N}$ the corresponding Lie group homomorphisms. The same procedure applied to the opposite BD triple $(\Gamma_2,\Gamma_1,\gamma^*)$ yields a pair of Lie algebra homomorphisms denoted by $\gamma^*$, and a pair of Lie group homomorphisms denoted by $\tilde{\gamma}^*$.

Note that we do not attempt to extend these
homomorphisms to the Lie algebra $\mathfrak{g}$
or to the corresponding algebraic group $G$. Such integrations are subtle and may not exist. In this paper, we only need the homomorphisms on the unipotent subgroups and their Lie algebras.

The following theorem is due to Belavin and Drinfeld~\cite{bd}:
\begin{theorem}
    Assume that $\mathfrak{g}$ is semisimple.
    For any BD quadruple $\bbg:=(\bg,r_0)$, the $2$-tensor $r \in \mathfrak{g}\otimes \mathfrak{g}$ given by
    \begin{equation}\label{eq:r_bd}
        r_{\bbg} = \Big(1 \otimes \frac{1}{1-\gamma}\Big)\Omega_{>} - \Big(1\otimes \frac{\gamma^*}{1-\gamma^*}\Big)\Omega_{<} + r_0
    \end{equation}
    is a solution of the CYBE. Conversely, given a solution $r$ of the CYBE such that $r + r^t = \Omega$, there exists a BD quadruple $\bbg$, a choice of root decomposition of $\mathfrak{g}$ and a choice of root vectors $e_{\alpha}$ such that $r$ is of the form~\eqref{eq:r_bd}.
\end{theorem}

Defining the linear map $R_{\bbg} : \mathfrak{g}\rightarrow \mathfrak{g}$ by $r_{\bbg} = (1\otimes R_{\bbg})\Omega$, equation~\eqref{eq:r_bd} becomes
\begin{equation*}
    R_{\bbg} = \frac{1}{1-\gamma}\pr_> - \frac{\gamma^*}{1-\gamma^*}\pr_{<} + R_0\pr_0.
\end{equation*}
Our constructions typically involve a pair of BD triples (resp. quadruples), which we denote by $(\bg^r,\bg^c)$ (resp. $(\bbg^r,\bbg^c)$). The corresponding pair of $r$-matrices is denoted by $r^r:=r_{\bbg^r}$ and $r^c := r_{\bbg^c}$. To distinguish between the simple roots in $\bg^r$ and $\bg^c$, we often decorate them with  superscripts $r$ and $c$, respectively.

The Poisson bracket $\{\cdot,\cdot\}$ corresponding to~\eqref{eq:biv} can be expressed as
\begin{equation*}
    \{f_1,f_2\} = \langle R^c(\nabla^L f_1),\nabla^L f_2\rangle - \langle R^r(\nabla^R f_1),\nabla^R f_2\rangle.
\end{equation*}
Setting $R_- := -R^*$ for both $r$ and $c$ superscripts, the Poisson bracket can also be written as
\begin{equation*}
    \{f_1,f_2\} = \langle R^c_-(\nabla^L f_1),\nabla^L f_2\rangle - \langle R^r_-(\nabla^R f_1),\nabla^R f_2\rangle.
\end{equation*}
We will use both expressions depending on convenience. To indicate that the Poisson bivector (resp. the Poisson bracket) depends on the BD quadruple $(\bbg^r,\bbg^c)$, we sometimes use subscripts, such as $\pi_{(\bbg^r,\bbg^c)}$ (resp. $\{f_1,f_2\}_{(\bbg^r,\bbg^c)}$).

In the case of $\Mat_{n\times m}(\mathbb{C})$, we consider a pair of BD quadruples $(\bbg^r,\bbg^c)$ where $\bbg^r$ is viewed in type $A_{n-1}$, while $\bbg^c$ is viewed in $A_{m-1}$. Given a pair of the corresponding $r$-matrices~\eqref{eq:r_bd}, we consider the Poisson bivector on $\Mat_{n\times m}(\mathbb{C})$ given by the same formula~\eqref{eq:biv} as in the case of $G$. 

There are two ways to see that $\pi$ on $\Mat_{n\times m}(\mathbb{C})$ given by~\eqref{eq:biv} is indeed a Poisson bivector. One involves a general argument using homogeneous spaces \cite[Proposition 4.2]{gy-homog}. Alternatively, one can compute explicitly that $(\Omega^{(m)})^L = (\Omega^{(n)})^R$, which allows us to rewrite $\pi$ in the form
\begin{equation*}
    \pi = (s^c)^L - (s^r)^R
\end{equation*}
where $s^r$ (resp. $s^c$) is the skew-symmetric part of $r^r$ (resp. $r^c$). For a given $r$-matrix satisfying the CYBE, its skew-symmetric part $s$ satisfies the \emph{modified CYBE} $[s,s] = -\frac{1}{3}\CYB(\Omega)$, where $[\cdot,\cdot]$ is the Schouten bracket. Using this, we compute that
\begin{equation*}
    [\pi,\pi] = [(s^c)^L - (s^r)^R,(s^c)^L - (s^r)^R]=  [(s^c)^L,(s^c)^L] + [(s^r)^R,(s^r)^R] = [s^c,s^c]^L - [s^r,s^r]^R = 0.
\end{equation*}
Thus~\eqref{eq:biv} defines a Poisson bivector on $\Mat_{n\times m}(\mathbb{C})$.

In the final part of the subsection, we describe the notion of a \emph{Belavin--Drinfeld torus} acting on~$G$. Let us decompose the Cartan subalgebra $\mathfrak{h}$ as $\mathfrak{h} = \mathfrak{h}_{\ses} \oplus \mathfrak{z}(\mathfrak{g})$ where $\mathfrak{h}_{\ses}$ is a Cartan subalgebra of the semisimple part $\mathfrak{g}_{\ses}$ of $\mathfrak{g}$, and $\mathfrak{z}(\mathfrak{g})$ is the center of $\mathfrak{g}$. A BD triple $\bg$ gives rise to a Lie subalgebra $\mathfrak{h}_{\bg}$ of $\mathfrak{h}_{\ses}$ given by
\begin{equation*}
    \mathfrak{h}_{\bg}:=\{h \in \mathfrak{h}_{\ses}\mid \alpha(h) = \gamma(\alpha)(h), \ \alpha \in \Gamma_1\}.
\end{equation*}
We let $\mathcal{H}_{\bg}$ be the corresponding connected subgroup of the Cartan subgroup $\mathcal{H}_{\ses}$ of the semisimple part of $G$.
Note that, since $\mathcal{H}_{\bg}$ is a subtorus of $\mathcal{H}_{\ses}$, there is a projection map $X^*(\mathcal{H}_{\ses}) \rightarrow X^*(\mathcal{H}_{\bg})$.

Denoting by $Z(G)$ the center of $G$ and by $Z(G)^{\circ}$ its identity component, we define the BD torus as follows.
\begin{definition}
    Given a pair $(\bg^r,\bg^c)$ of BD triples, the \emph{Belavin--Drinfeld torus} (the \emph{BD torus}) is the complex algebraic torus given by
    \begin{equation*}
        \mathbb{T}_{(\bg^r,\bg^c)}:= \mathcal{H}_{\bg^r} \times Z(G)^{\circ} \times \mathcal{H}_{\bg^c},
    \end{equation*}
    and acting on $G$ as
    \begin{equation*}
        (h^r,z,h^c)\cdot g = z h^rgh^c, \quad (h^r,z,h^c) \in \mathbb{T}_{(\bg^r,\bg^c)}, \ g \in G.
    \end{equation*}
\end{definition}
For a Poisson bivector $\pi$ defined by a pair $(\bbg^r,\bbg^c)$, the action of the BD torus $\mathbb{T}_{(\bg^r,\bg^c)}$ on $(G,\pi)$ is Poisson.

The following notion allows us to construct
a basis of $\mathfrak{h}_{\bg}$.
\begin{definition}\label{d:bd_chain}
    Given a BD triple $\bg:=(\Gamma_1,\Gamma_2,\gamma)$, let $\sim$ be the equivalence relation on $\Delta$ generated by the relations $\alpha\sim\gamma(\alpha)$ for $\alpha\in\Gamma_1$. An equivalence class is called a \emph{BD chain}.
\end{definition}

The BD chains are in bijection with the set $\Delta \setminus \Gamma_2$. Each $\alpha \in \Delta\setminus \Gamma_2$ gives rise to a BD chain $C$ of the form
\begin{equation*}
    C = \{\alpha,\gamma(\alpha),\gamma^2(\alpha),\ldots, \gamma^j(\alpha)\}
\end{equation*}
where $j$ is such that $\gamma^{j}(\alpha) \notin \Gamma_1$. For a given BD chain $C$, define
\begin{equation*}
    h_C := \sum_{\alpha_i \in C} \varpi_i^{\vee}
\end{equation*}
where $\varpi_i^{\vee}$ is the fundamental coweight corresponding to $\alpha_i$. In other words, $\alpha_j(\varpi_i^{\vee}) = \delta_{ij}$. Then the set of elements $h_C$ for all BD chains $C$ is a basis of $\mathfrak{h}_{\bg}$. In particular, $\dim \mathfrak{h}_{\bg} = |\Delta \setminus \Gamma_2|$.

In the case of $\Mat_{n\times m}(\mathbb{C})$, the corresponding group $\mathcal{H}_{\bg^r}$ (resp.  $\mathcal{H}_{\bg^c}$) is a connected subgroup of the group of diagonal matrices $\mathcal{H}_{\ses}^{(n)}$ of $\SL_n(\mathbb{C})$ (resp. $\mathcal{H}_{\ses}^{(m)}$ of $\SL_m(\mathbb{C})$). The BD torus is given by
\begin{equation*}
    \mathbb{T}_{(\bg^r,\bg^c)} := \mathcal{H}_{\bg^r} \times \mathbb{C}^\times \times \mathcal{H}_{\bg^c},
\end{equation*}
and it acts on $\Mat_{n\times m}(\mathbb{C})$ in the same way as in the case of a group $G$. The character lattice $X^*(\mathbb{T}_{(\bg^r,\bg^c)})$ is given by the direct sum
\begin{equation*}
    X^*(\mathbb{T}_{(\bg^r,\bg^c)}) = X^*(\mathcal{H}_{\bg^r}) \oplus \mathbb{Z}\langle e \rangle \oplus X^*(\mathcal{H}_{\bg^c})
\end{equation*}
where $e$ is a basis element of $X^*(\mathbb{C}^\times)$.

We denote by
\begin{equation*}
    \pr_r : X^*(\mathcal{H}_{\ses}^{(n)}) \rightarrow X^*(\mathcal{H}_{\bg^r}), \quad \pr_c : X^*(\mathcal{H}_{\ses}^{(m)}) \rightarrow X^*(\mathcal{H}_{\bg^c})
\end{equation*}
the projection maps. To distinguish between the weights of $X^*(\mathcal{H}_{\ses}^{(n)})$ and $X^*(\mathcal{H}_{\ses}^{(m)})$, we use superscripts $r$ and $c$, respectively. For instance, $\varpi_i^r$ denotes a fundamental weight belonging to $X^*(\mathcal{H}_{\ses}^{(n)})$. We further introduce the weights
\begin{align*}
    &\beta_i^r := \varpi_i^r - \varpi_{i-1}^r \in X^*(\mathcal{H}_{\ses}^{(n)}), \quad i \in [1,n];\\
    &\beta_j^c := \varpi_j^c - \varpi_{j-1}^c \in X^*(\mathcal{H}_{\ses}^{(m)}), \quad j \in [1,m],
\end{align*}
where $\varpi_{0}^r = \varpi_n^r = 0$ and $\varpi_0^c =\varpi_m^c = 0$. Then a matrix entry $x_{ij}$ is a $\mathbb{T}_{(\bg^r,\bg^c)}$-weight vector, with the weight given by
\begin{equation*}
    \beta_{ij}:=\pr_r(\beta_i^r) + e + \pr_c(\beta_j^c).
\end{equation*}
Note that if the pair $(\bg^r,\bg^c)$ consists of the standard BD triples, then $\mathcal{H}_{\bg^r} = \mathcal{H}_{\ses}^{(n)}$ and $\mathcal{H}_{\bg^c} = \mathcal{H}_{\ses}^{(m)}$. Correspondingly, $\pr_r$ and $\pr_c$ are the identity maps. In this case, we denote by $(\beta_{\std})_{ij}$ the $\mathbb{T}_{(\bg_{\std},\bg_{\std})}$-weights of~$x_{ij}$.

\subsection{Poisson deformations}\label{s:background_deform}

In this subsection, we describe background material on Poisson deformations of $\mathbb{T}$-log-symplectic log-canonical Poisson bivectors on a complex affine space $\mathsf{X}:=\mathbb{C}^n$. The main reference is~\cite{lm-deform}. We describe only the material relevant to our setting. In particular, we do not describe formal deformations and Property (W0). 

The space $\mathsf{X}$ is equipped with the standard coordinates $x_1,\ldots,x_n$ and an action of the \emph{big torus} $(\mathbb{C}^\times)^n$ on $\mathsf{X}$ such that
\begin{equation*}
    (\lambda_1,\ldots,\lambda_n)\cdot (x_1,\ldots,x_n) \mapsto (\lambda_1 x_1,\ldots,\lambda_n x_n).
\end{equation*}
We identify the character lattice of $(\mathbb{C}^\times)^n$ with $\mathbb{Z}^n$ so that the weight of $x_j$ is $e_j$.
The action of $(\mathbb{C}^\times)^n$ induces an action on the space of algebraic polyvector fields $\mathfrak{X}(\mathsf{X})$, so that the weight of $\partial/\partial x_j$ is given by $-e_j$. Recall that an algebraic $p$-vector field has the form
\begin{equation*}
    V = \sum_{1 \leq i_1 < \cdots < i_p \leq n} V_{i_1,\ldots,i_p} \frac{\partial}{\partial x_{i_1}}\wedge \cdots \wedge \frac{\partial}{\partial x_{i_p}}, \quad V_{i_1,\ldots,i_p} \in \mathbb{C}[\mathsf{X}].
\end{equation*}
The subspace of algebraic $p$-vector fields is denoted by $\mathfrak{X}^p(\mathsf{X})$, and the space $\mathfrak{X}(\mathsf{X})$ is a graded algebra
\begin{equation*}
    \mathfrak{X}(\mathsf{X}) = \bigoplus_{p=0}^n \mathfrak{X}^p(\mathsf{X}).
\end{equation*}
In particular, $\mathfrak{X}^0(\mathsf{X}) = \mathbb{C}[\mathsf{X}]$. The graded algebra $\mathfrak{X}(\mathsf{X})$ is equipped with the Schouten bracket $[\cdot,\cdot]$.

Each subspace $\mathfrak{X}^p(\mathsf{X})$ decomposes into its $(\mathbb{C}^\times)^n$-weight subspaces. Given $\bw \in (\mathbb{Z}_{\geq -1})^n$, we denote by $\mathfrak{X}^p(\mathsf{X})^{\bw}$ the corresponding weight subspace. It can be described as follows. Define
\begin{equation*}
    J_{\bw} := \{j \in [1,n] \mid \bw_j = -1\},
\end{equation*}
and let
\begin{equation*}
    \partial_j := x_j \frac{\partial}{\partial x_j}, \quad \partial_J := \partial_{j_1} \wedge \cdots \wedge \partial_{j_k}, \ J = \{j_1,\ldots,j_k\}.
\end{equation*}
Note that the $(\mathbb{C}^\times)^n$-weight of any $\partial_J$ equals zero.
By \cite[Lemma 2.1.1]{lm-deform}, $\mathfrak{X}^p(\mathsf{X})^{\bw} \neq 0$ if and only if $|J_{\bw}|\leq p$. More precisely, 
\begin{equation}\label{eq:xpw}
    \mathfrak{X}^p(\mathsf{X})^{\bw} = \Span_{\mathbb{C}}\{ x^{\bw} \partial_J \mid J_{\bw} \subset J \subset [1,n], \ |J| = p\},
\end{equation}
and in particular, 
\begin{equation*}
    \dim \mathfrak{X}^p(\mathsf{X})^{\bw} = \binom{n-|J_{\bw}|}{p-|J_{\bw}|}.
\end{equation*}

A \emph{log-canonical bivector} $\pi_0$ on $\mathsf{X}$ is a $(\mathbb{C}^\times)^n$-invariant Poisson bivector. That is,
\begin{equation*}
    \pi_0 = \sum_{1 \leq j < k \leq n} \lambda_{jk} x_j x_k \frac{\partial}{\partial x_j}\wedge \frac{\partial}{\partial x_k}.
\end{equation*}
We denote by $\blambda:=(\lambda_{jk})_{j,k=1}^n$ the matrix of coefficients of $\pi_0$, and we let $\mathsf{d}_{\pi_0} : \mathfrak{X}(\mathsf{X}) \rightarrow \mathfrak{X}(\mathsf{X})$ be the differential given by $\mathsf{d}_{\pi_0}(V) := [\pi_0,V]$. The key observation about $\mathsf{d}_{\pi_0}$ is that it preserves the $(\mathbb{C}^{\times})^n$-weight subspaces of $\mathfrak{X}(\mathsf{X})$. More precisely, setting $\underline{\partial}:=(\partial_1,\ldots,\partial_n)$, by \cite[Lemma 2.1.3]{lm-deform},

\begin{equation*}
    \mathsf{d}_{\pi_0}(\mathfrak{X}^p(\mathsf{X})^{\bw}) \subseteq \mathfrak{X}^{p+1}(\mathsf{X})^{\bw},\quad
    \mathsf{d}_{\pi_0}|_{\mathfrak{X}^p(\mathsf{X})^{\bw}}(V) = \underline{\partial}\blambda \bw \wedge V.
\end{equation*}

The Poisson cohomology $\mathsf{H}_{\pi_0}(\mathsf{X})$ is defined as the cohomology of the cochain complex
\begin{equation*}
    0 \rightarrow \mathfrak{X}^0(\mathsf{X}) \xrightarrow{\mathsf{d}_{\pi_0}} \mathfrak{X}^1(\mathsf{X}) \xrightarrow{\mathsf{d}_{\pi_0}} \cdots \xrightarrow{\mathsf{d}_{\pi_0}} \mathfrak{X}^n(\mathsf{X})\rightarrow 0.
\end{equation*}
Since $\mathsf{d}_{\pi_0}$ preserves the $(\mathbb{C}^{\times})^n$-weight decomposition of $\mathfrak{X}(\mathsf{X})$, it follows that
\begin{equation*}
    \mathsf{H}_{\pi_0}(\mathsf{X}) = \bigoplus_{\bw \in (\mathbb{Z}_{\geq -1})^n} \mathsf{H}_{\pi_0}(\mathsf{X})^{\bw}.
\end{equation*}

For a subset $J \subseteq [1,n]$ of indices, define the following spaces:
\begin{equation*}
    \mathbb{C}^J := \sum_{j \in J} \mathbb{C}e_j, \quad (\mathbb{C}^J)_0 := \Big\{ \sum_{j \in J} a_j e_j \in \mathbb{C}^J \mid \sum_{j \in J} a_j = 0\Big\}.
\end{equation*}
The following proposition is a key structural result concerning the $(\mathbb{C}^\times)^n$-weight subspaces of $\mathsf{H}_{\pi_0}(\mathsf{X})$:
\begin{proposition}\label{p:hxw}
    (\cite[Proposition 2.1.5]{lm-deform}) For $\bw \in (\mathbb{Z}_{\geq -1})^n$, one has $\mathsf{H}_{\pi_0}(\mathsf{X})^{\bw} \neq 0$ if and only if
    \begin{equation*}
        \blambda \bw \in \mathbb{C}^{J_{\bw}}, \quad \text{equivalently,} \quad \blambda \bw \in (\mathbb{C}^{J_{\bw}})_0,
    \end{equation*}
    and in this case, $\mathsf{d}_{\pi_0}|_{\mathfrak{X}(\mathsf{X})^{\bw}} = 0$ and $\mathsf{H}_{\pi_0}(\mathsf{X})^{\bw} \cong \mathfrak{X}(\mathsf{X})^{\bw}$.
\end{proposition}
Proposition~\ref{p:hxw} allows us to identify $\mathsf{H}_{\pi_0}(\mathsf{X})^{\bw}$ with $\mathfrak{X}(\mathsf{X})^{\bw}$ when $\mathsf{H}_{\pi_0}(\mathsf{X})^{\bw} \neq 0$.

We now introduce an action of another complex algebraic torus $\mathbb{T}$ on $\mathsf{X}$, given by
\begin{equation*}
    t\cdot (x_1,\ldots,x_n) = (t^{\beta_1}x_1,\ldots,t^{\beta_n}x_n)
\end{equation*}
where each $\beta_j \in X^*(\mathbb{T})$ is a character.

It follows from equation~\eqref{eq:xpw} that, for a given $\bw \in (\mathbb{Z}_{\geq -1})^n$ with $|J_{\bw}|\leq p$, the $\mathbb{T}$-weight of the elements in  $\mathfrak{X}^p(\mathsf{X})^{\bw}$ is given by $\bbeta \bw$. Hence, an element in $\mathfrak{X}^p(\mathsf{X})^{\bw}$ is $\mathbb{T}$-invariant if and only if $\bbeta \bw = 0$. Therefore, we obtain a $(\mathbb{C}^\times)^n$-weight decomposition of $\mathfrak{X}(\mathsf{X})^{\mathbb{T}}$ and $\mathsf{H}_{\pi_0}(\mathsf{X})^{\mathbb{T}}$:
\begin{equation*}
    \mathfrak{X}(\mathsf{X})^{\mathbb{T}} = \bigoplus_{\bw \in (\mathbb{Z}_{\geq -1})^n, \bw \in \ker \bbeta} \mathfrak{X}(\mathsf{X})^{\bw}, \qquad 
    \mathsf{H}_{\pi_0}(\mathsf{X})^{\mathbb{T}} = \bigoplus_{\bw \in (\mathbb{Z}_{\geq -1})^n, \bw \in \ker \bbeta}\mathsf{H}_{\pi_0}(\mathsf{X})^{\bw}.
\end{equation*}

First-order $\mathbb{T}$-invariant Poisson deformations of $\pi_0$ are controlled by the second $\mathbb{T}$-invariant Poisson cohomology $\mathsf{H}^2_{\pi_0}(\mathsf{X})^{\mathbb{T}}$. Setting
\begin{equation}\label{eq:spi0}
\begin{split}
    \mathcal{S}(\pi_0) &:= \mathcal{S}(\pi_0,\mathbb{T}) := \{ \text{the nonzero }(\mathbb{C}^\times)^n\text{-weights in } \mathsf{H}^2_{\pi_0}(\mathsf{X})^{\mathbb{T}}\} = \\ &= \{\btheta \in (\mathbb{Z}_{\geq -1})^n \mid \, |J_{\btheta}| \leq 2, \ \blambda \btheta \in (\mathbb{C}^{J_{\btheta}})_0, \  \bbeta \btheta = 0\},
    \end{split}
\end{equation}
we can therefore write
\begin{equation*}
    \mathsf{H}^2_{\pi_0}(\mathsf{X})^{\mathbb{T}} = \mathsf{H}^2_{\pi_0}(\mathsf{X})^{0} \oplus \bigoplus_{\btheta \in \mathcal{S}(\pi_0)} \mathsf{H}^2_{\pi_0}(\mathsf{X})^{\btheta}.
\end{equation*}
In particular, by equation~\eqref{eq:xpw} and Proposition~\ref{p:hxw}, for any $\btheta \in \mathcal{S}(\pi_0)$ such that $|J_{\btheta}| = 2$,
\begin{equation}\label{eq:coh1}
    \mathsf{H}_{\pi_0}^2(\mathsf{X})^{\btheta} = \mathbb{C}V_{\btheta}, \quad V_{\btheta} := x^{\btheta} \partial_{J_{\btheta}}.
\end{equation}
Note that for any $\btheta \in \mathcal{S}(\pi_0)$ with $|J_{\btheta}| \leq 1$, $\btheta \in \ker \blambda \cap \ker \bbeta$, and in this case, $\mathsf{H}_{\pi_0}^2(\mathsf{X})^{\btheta}$ is not one-dimensional. Thus we naturally arrive at the following definition:
\begin{definition}
    The log-canonical bivector $\pi_0$ is called \emph{$\mathbb{T}$-log-symplectic} if $\ker \bbeta \cap \ker \blambda = \{0\}$.
\end{definition}
Equivalently, $\pi_0$ is $\mathbb{T}$-log-symplectic if there is a symplectic leaf $S$ such that $\mathbb{T}\cdot S$ is open. Since $\pi_0$ is log-canonical, if $\mathbb{T}\cdot S$ is open, then necessarily $\mathbb{T}\cdot S = (\mathbb{C}^\times)^n$.

One of our new results is that the weights $\btheta$ in $\mathcal{S}(\pi_0,\mathbb{T})$ satisfy $|J_{\btheta}| = 2$ for log-canonical bivectors that are not necessarily $\mathbb{T}$-log-symplectic but satisfy a different condition (see Proposition~\ref{prop:SisLM}).

We now investigate combinatorial properties of $\mathcal{S}(\pi_0,\mathbb{T})$ more closely. Define the following set:
\begin{equation*}
    \mathcal{W}_p := \{\bw \in (\mathbb{Z}_{\geq -1})^n \mid \, |J_{\bw}| \leq p\} = \{\mathbf{w}\in \mathbb{Z}^n\mid \mathfrak{X}(\mathsf{X})^\mathbf{w} \neq 0\}.
\end{equation*}
Note that an element $\btheta \in (\mathbb{Z}_{\geq -1})^n$ belongs to $\mathcal{W}_2 \setminus \mathcal{W}_1$ if and only if $|J_{\btheta}| = 2$.

\begin{definition}
    Given distinct indices $j,k \in [1,n]$ and $a \in \mathbb{C}^\times$, the \emph{Lu--Matviichuk equations} (the \emph{LM equations}) are the following system of linear equations for $\btheta\in  \mathcal{W}_2 \setminus \mathcal{W}_1$ with $J_{\btheta} = \{j,k\}$:
    \begin{equation*}
        \begin{cases}
            \blambda \btheta = a(e_j - e_k) \in \mathbb{C}^n&\\
            \bbeta \btheta = 0 \in X^*(\mathbb{T})
        \end{cases}
    \end{equation*}
\end{definition}
The first equation is precisely the condition $\blambda \btheta \in (\mathbb{C}^{J_{\btheta}})_0$ from~\eqref{eq:spi0}. The number $a$ is determined by $\btheta$ up to a sign, and \cite{lm-deform} chooses $j<k$ to make $a=a_{\btheta}$ unambiguous.

Combinatorial properties of solutions to the LM equations can be conveniently expressed in terms of the following notion.
\begin{definition}
    Given $\mathcal{S}\subseteq \mathcal{S}(\pi_0,\mathbb{T})$, the \emph{smoothing graph} $\Gamma(\mathcal{S})$ is an undirected graph on $n$ vertices whose edge set is formed as follows. For each $j,k \in [1,n]$, there is an edge $j\edge k$ if and only if there is $\btheta \in \mathcal{S}$ such that $J_{\btheta} = \{j,k\}$. The \emph{smoothing graph of $(\pi_0,\mathbb{T})$} is the smoothing graph $\Gamma(\pi_0,\mathbb{T}) := \Gamma(\mathcal{S}(\pi_0,\mathbb{T}))$.
\end{definition}
For a given $\btheta$ with $|J_{\btheta}| = 2$, we denote by $\varepsilon^{\btheta}$ the corresponding edge in the smoothing graph. For a given edge $\varepsilon$ in the smoothing graph, we denote by $\btheta^{\varepsilon}$ the corresponding element of $\mathcal{S}(\pi_0,\mathbb{T})$. We thus have a bijection between $\mathcal{S}(\pi_0,\mathbb{T}) \cap (\mathcal{W}_2 \setminus \mathcal{W}_1)$ and the set of edges of $\Gamma(\pi_0,\mathbb{T})$. Note that if $\pi_0$ is $\mathbb{T}$-log-symplectic, then $\mathcal{S}(\pi_0,\mathbb{T}) \cap (\mathcal{W}_2\setminus\mathcal{W}_1) = \mathcal{S}(\pi_0,\mathbb{T})$, and therefore, there is a bijection between $\mathcal{S}(\pi_0,\mathbb{T})$ and the set of edges of $\Gamma(\pi_0,\mathbb{T})$.

The following proposition summarizes combinatorial properties of $\Gamma(\pi_0,\mathbb{T})$:
\begin{proposition}\label{p:lm_bunch}
    The following statements hold:
    \begin{itemize}
        \item (\cite[Lemma 2.2.9]{lm-deform}) If $\btheta, \btheta^\prime \in \mathcal{S}(\pi_0)$ are such that $\blambda \btheta = a(e_i-e_j)$ and $\blambda \btheta^\prime = a^\prime(e_j-e_k)$, then $a/a^\prime \in \mathbb{Q}_{> 0}$; 
        \item (\cite[Lemma 2.2.10]{lm-deform}) Each vertex of $\Gamma(\pi_0,\mathbb{T})$ is of degree $\leq 2$. In particular, $\Gamma(\pi_0,\mathbb{T})$ is a disjoint union of chains and cycles;
        \item (\cite[Lemma-Definition 2.2.12]{lm-deform}) For any chain $c \in \Gamma(\pi_0)$, $\{\btheta^{\varepsilon} \mid \varepsilon \in E(c)\}$ is $\mathbb{C}$-linearly independent. Conversely, if $\pi_0$ is $\mathbb{T}$-log-symplectic and $c$ is a connected component of $\Gamma(\pi_0)$ such that $\{\btheta^{\varepsilon} \mid \varepsilon \in E(c)\}$ is $\mathbb{C}$-linearly independent, then $c$ is a chain;
        \item (\cite[Lemma 2.2.6]{lm-deform}) If $\pi_0$ is $\mathbb{T}$-log-symplectic, then $\Gamma(\pi_0,\mathbb{T})$ has no $2$-cycles.
    \end{itemize}
\end{proposition}
If $\pi_0$ is not $\mathbb{T}$-log-symplectic, the graph $\Gamma(\pi_0,\mathbb{T})$ might have $2$-cycles. That is, there might be $\btheta,\btheta^\prime \in \mathcal{S}(\pi_0)$ such that $J_{\btheta} = J_{\btheta^\prime}$.

For the next results, let us introduce the following notation. Given a subset $\mathsf{W} \subseteq \mathbb{Z}^n$, write
\begin{equation*}
    \mathsf{W}_m := \underbrace{\mathsf{W}+\cdots +\mathsf{W}}_{m \text{ times}} \qquad \mathsf{W}_{\geq m}=\bigcup_{m'\geq m}\mathsf{W}_{m'} \qquad \mathfrak{X}^p(\mathsf{X})^{\mathsf{W}} = \sum_{\mathbf{w}\in\mathsf{W}}
\mathfrak{X}^p(\mathsf{X})^{\mathbf{w}}.
\end{equation*}
\begin{definition}\label{def:S-adm}(\cite[Definition 3.2.1]{lm-deform}) 
Assume that a finite subset $\mathsf{W} \subset (\mathbb{Z}_{\geq -1})^n$ is linearly independent. Given $\pi_1\in \mathfrak{X}^2(\mathsf{X})^{\mathsf{W} }$, if 
$\pi\in \mathfrak{X}^2(\mathsf{X})^{\mathsf{W} _{\geq 0}}$ is Poisson and satisfies
\[
\pi=\pi_0+\pi_1\;\;\;{\rm mod}\;\;\mathfrak{X}^2(\mathsf{X})^{\mathsf{W} _{\geq 2}},
\]
we say that $\pi$ is a {\it $\mathsf{W}$-admissible algebraic Poisson deformation of $\pi_0$ along $\pi_1$}.
\end{definition}

By \cite[Theorem 2.3.4]{lm-deform}, for a $\mathbb{T}$-log-symplectic $\pi_0$ and a given linearly independent $\mathcal{S} \subseteq \mathcal{S}(\pi_0,\mathbb{T})$, the Poisson deformation along $\pi_1 \in \mathfrak{X}^2(\mathsf{X})^{\mathcal{S}}$ exists at least in the formal sense. The proof of the theorem relies essentially on the $\mathbb{T}$-log-symplecticity to show that $\mathsf{H}^3_{\pi_0}(\mathsf{X})^{\mathcal{S}_{\geq 2}} = 0$. To ensure that $\pi$ exists algebraically rather than only formally, Lu and Matviichuk impose an additional condition, the so-called Property (W0) (see~\cite[Definition 3.1.7]{lm-deform}). Under Property (W0), although $\pi$ exists algebraically, it is unique only up to a finite-dimensional unipotent gauge group (see \cite[Theorem 3.2.4]{lm-deform}).

In our setting, we will see that Poisson deformations of the log-canonical parts of aperiodic Belavin--Drinfeld bivectors on $\Mat_{n\times m}(\mathbb{C})$ satisfy an even stronger condition.

\begin{definition}(\cite[Definition 3.3.1]{lm-deform})
    A finite subset $\mathsf{W} \subset (\mathbb{Z}_{\geq -1})^n$ has Property~\eqref{eq:W1} if
    \begin{equation}\label{eq:W1}
    \mathsf{W}_{\geq 1}\cap \mathcal{W}_1=\varnothing,
    \qquad\text{equivalently}\qquad
    \mathfrak{X}^1(\mathsf{X})^{\mathsf{W}_{\geq 1}}=0.
    \tag{W1}
    \end{equation}
    The log-canonical Poisson bivector $\pi_0$ on $\mathsf{X}$ has Property~\eqref{eq:W1} if $\mathcal{S}(\pi_0,\mathbb{T})$ has Property~\eqref{eq:W1}.
\end{definition}

The following theorem is an abridged version of \cite[Theorem 3.3.10]{lm-deform}, and is a key result on Poisson deformations with Property~\eqref{eq:W1}.

\begin{theorem}\label{thm:W1}
Let \(\pi_0\) be any \(\mathbb{T}\)-log-symplectic log-canonical Poisson structure on
\(\mathsf{X}=\mathbb{C}^n\).
Suppose that \(\mathcal{S}\subseteq \mathcal{S}(\pi_0)\) has Property~\eqref{eq:W1}. Then, for every
\(\pi_1\in \mathfrak{X}^2(\mathsf{X})^{\mathcal{S}}\),
there exists a unique \(\mathcal{S}\)-admissible algebraic Poisson deformation
\(\pi\) of \(\pi_0\) along \(\pi_1\). Moreover, \(\pi\) has Property~\eqref{eq:W1}.
\end{theorem}

In the final part of this subsection, we discuss the notion of a normalized Poisson deformation. Given a solution $\btheta$ of the LM equations such that $\blambda \btheta = a_{\btheta}(e_j-e_k)$, define
\begin{equation}\label{eq:Vprime}
    V_\btheta^{\prime}:={a_\btheta}x^\btheta \partial_j \wedge \partial_k.
\end{equation}
Note that $\btheta$ also satisfies the equation $\blambda \btheta = (-a_{\btheta})(e_k-e_j)$. However, the bivector $V_{\btheta}^\prime$ is defined unambiguously because
\begin{equation*}
    a_\btheta x^\btheta \partial_j \wedge \partial_k = (-a_\btheta)x^\btheta \partial_k \wedge \partial_j.
\end{equation*}

\begin{definition}
For a subset $\mathcal{S}\subseteq \mathcal{S}(\pi_0,\mathbb{T})$ with Property~\eqref{eq:W1} and parameters $u = (u_{\btheta})_{\btheta \in \mathcal{S}} \in \mathbb{C}^{\mathcal{S}}$, define
\begin{equation*}\label{eq:pi1s}
    \pi_1^{\mathcal{S}}(u) := \sum_{\btheta \in \mathcal{S}} u_{\btheta} V_{\btheta}^\prime.
\end{equation*}
We denote the unique $\mathcal{S}$-admissible algebraic Poisson deformation of $\pi_0$ along $\pi_1^{\mathcal{S}}(u)$ by
\begin{equation}\label{eq:pi-S-c}
\pi^{\mathcal{S}}(u)
=
\pi_0+\pi_1^{\mathcal{S}}(u)+\pi_2^{\mathcal{S}}(u)+\cdots+\pi_N^{\mathcal{S}}(u),
\end{equation}
where $N\geq 1$ and
$\pi_m^{\mathcal{S}}(u)\in
\mathfrak{X}^2(\mathsf{X})^{\mathcal{S}_m}$
for $2\leq m\leq N$.
We call
$\pi^{\mathcal{S}}(u)$ the
\(\mathcal{S}\)-admissible \emph{family} of algebraic Poisson deformations of \(\pi_0\) and
\(
\pi^{\mathcal{S}}
\bigl(
(1_{\btheta})_{\boldsymbol{\theta}\in\mathcal{S}}
\bigr)
\)
the \emph{normalized}
\(\mathcal{S}\)-admissible algebraic Poisson deformation of \(\pi_0\). If $\mathcal{S} = \mathcal{S}(\pi_0,\mathbb{T})$, they are called 
the \emph{maximal family of admissible $\mathbb{T}$-invariant algebraic Poisson deformations} of \(\pi_0\) and
the \emph{maximal normalized admissible $\mathbb{T}$-invariant algebraic Poisson deformation} of \(\pi_0\).
\end{definition}

Our definition of $\pi_1^{\mathcal{S}}(u)$ differs from~\cite[Notation 3.3.6]{lm-deform} because we use different generators of $\mathsf{H}_{\pi_0}^{2}(\mathsf{X})^{\btheta}$. When $\pi_0$ is of symmetric action type, our definition of a normalized Poisson deformation agrees with \cite[Definition 4.4.5]{lm-deform}. We believe that our notion of normality provides an appropriate generalization of the Goodearl--Yakimov normality condition for symmetric Poisson CGL extensions \cite{lm-deform,gy-cgl}. In particular, we expect normalized deformations in the sense developed here to admit cluster structures analogous to those associated with Goodearl--Yakimov normal symmetric Poisson CGL extensions.

\section{Poisson deformation theory}\label{s:pdeform}
In this section, we describe some new results in Poisson deformation theory that we later apply in the setting of Belavin--Drinfeld bivectors on $\Mat_{n\times m}(\mathbb{C})$. Throughout this section, we use the same setup and notation as in Section~\ref{s:background_deform}. In particular, $\pi_0$ denotes a log-canonical Poisson bivector on $\mathbb{C}^n$, and $\mathbb{T}$ denotes a complex algebraic torus acting on $\mathbb{C}^n$ such that the standard coordinate functions are $\mathbb{T}$-weight vectors. We do not assume that $\pi_0$ is $\mathbb{T}$-log-symplectic unless otherwise stated.

\subsection{\texorpdfstring{Cycles and $\mathbb{T}$-log-symplecticity}{Cycles and T-log-symplecticity}}\label{s:cycles}

In the case of BD bivectors on $\Mat_{n\times m}(\mathbb{C})$, we will show that a certain collection of solutions of the LM equations forms a basis of $\ker \bbeta$ (without assuming aperiodicity). Surprisingly, this allows one to characterize the intersection $\ker \blambda \cap \ker \bbeta$, and thus describe obstructions to $\mathbb{T}$-log-symplecticity. In this subsection, we formulate the relevant results in the abstract setting.

Recall from \cite[Lemma 2.2.10]{lm-deform} that the smoothing graph $\Gamma(\pi_0,\mathbb{T})$ is a disjoint union of chains and cycles. By choosing an orientation for each chain and cycle, we obtain the source and target maps $s$ and $t$. That is, given an oriented edge $\varepsilon : j \rightedge j^\prime$, $s(\varepsilon) = j$ and $t(\varepsilon) = j^\prime$. With these orientations fixed, the map $\varepsilon \mapsto a^{\varepsilon} \in \mathbb{C}^{\times}$ is well-defined, where $a^{\varepsilon}$ is determined by the equation
\begin{equation*}
    \blambda \btheta^{\varepsilon} = a^{\varepsilon} (e_{s(\varepsilon)}-e_{t(\varepsilon)}).
\end{equation*}
Moreover, \cite[Lemma 2.2.9]{lm-deform} implies that if $\varepsilon$ and $\varepsilon^\prime$ are two consecutive edges, then $a^{\varepsilon}/a^{\varepsilon^\prime} \in \mathbb{Q}_{>0}$. In particular, if $\blambda$ is defined over $\mathbb{Q}$, the scalars $a^{\varepsilon}$ have the same sign along each chain and cycle.

Given an oriented cycle $c$ in $\Gamma(\mathcal{S})$, let $E(c)$ denote its edge set. Define
\begin{equation}\label{eq:bthetac}
    \btheta(c) := \sum_{\varepsilon \in E(c)} (a^{\varepsilon})^{-1}\btheta^{\varepsilon}.
\end{equation}
Note that reversing the orientation of $c$ changes $\btheta(c)$ to $-\btheta(c)$. Clearly, $\btheta(c) \in \ker \blambda \cap \ker \bbeta$.

\begin{proposition}\label{p:ker_cycles}
    Suppose there exists a subset $\mathcal{S}^\prime \subseteq \mathcal{S}(\pi_0)$ that is a basis of $\ker \bbeta$. Let $c_1,\ldots,c_r$ be all the cycles in $\Gamma(\mathcal{S}^\prime)$, each with a fixed orientation. Then the set $\btheta(c_1),\ldots,\btheta(c_r)$ is a basis of $\ker \blambda \cap \ker \bbeta$. In particular, the dimension of $\ker \blambda \cap \ker \bbeta$ is the number of cycles in $\Gamma(\mathcal{S}^\prime)$.
\end{proposition}
\begin{proof}
    Let us observe that for any oriented cycle $c$, if $b^{\varepsilon} \in \mathbb{C}$ are any scalars satisfying
    \begin{equation*}
        \sum_{\varepsilon \in E(c)} b^{\varepsilon} (e_{s(\varepsilon)}-e_{t(\varepsilon)}) = 0,
    \end{equation*}
    then $b^{\varepsilon} = b^{\varepsilon^\prime}$ for all $\varepsilon,\varepsilon^\prime \in E(c)$.

    Since the cycles are disjoint and $\mathcal{S}^\prime$ is linearly independent, the elements $\btheta(c_i)$ are linearly independent.
    Since $\mathcal{S}^\prime$ is a basis of $\ker \bbeta$, for any $\bw \in \ker \blambda \cap \ker \bbeta$, we can write
    \begin{equation*}
        \bw = \sum_{\varepsilon \in E(\Gamma(\mathcal{S}^\prime))} b^\varepsilon (a^{\varepsilon})^{-1}\btheta^\varepsilon, \quad b^\varepsilon \in \mathbb{C}.
    \end{equation*}
    Applying $\blambda$ to this expression, we obtain
    \begin{equation*}
        \blambda \bw = \sum_{\varepsilon \in E(\Gamma(\mathcal{S}^\prime))} b^\varepsilon (a^{\varepsilon})^{-1} \blambda \btheta^\varepsilon = \sum_{\varepsilon \in E(\Gamma(\mathcal{S}^\prime))} b^\varepsilon (e_{s(\varepsilon)}-e_{t(\varepsilon)}) = \sum_{i=1}^r \sum_{\varepsilon \in E(c_i)} b^\varepsilon(e_{s(\varepsilon)}-e_{t(\varepsilon)}) = 0.
    \end{equation*}
    By the earlier observation, since the cycles are disjoint, we see that $b^\varepsilon = b^{\varepsilon^\prime}$ whenever $\varepsilon$ and $\varepsilon^\prime$ belong to the same cycle $c$; hence, we can set $b^{c}:=b^{\varepsilon}$ for any $\varepsilon \in E(c)$. Therefore, $\bw$ can be expressed as
    \begin{equation*}
        \bw = \sum_{i=1}^r b^{c_i} \sum_{\varepsilon \in E(c_i)} (a^{\varepsilon})^{-1} \btheta^{\varepsilon} = \sum_{i=1}^r b^{c_i} \btheta(c_i),
    \end{equation*}
    and thus $\{\btheta(c_1),\ldots,\btheta(c_r)\}$ is a basis of $\ker \blambda \cap \ker \bbeta$.
\end{proof}
\begin{corollary}\label{c:cycles_tlogs}
   Suppose there exists a subset $\mathcal{S}^\prime \subseteq \mathcal{S}(\pi_0)$ that is a basis of $\ker \bbeta$. Then $\pi_0$ is $\mathbb{T}$-log-symplectic if and only if the smoothing graph of $\mathcal{S}^\prime$ has no cycles.
\end{corollary}
\begin{proof}
    By the previous proposition, the dimension of $\ker \blambda \cap \ker \bbeta$ is the number of cycles in $\Gamma(\mathcal{S}^\prime)$. Since $\pi_0$ is $\mathbb{T}$-log-symplectic if and only if $\ker \blambda \cap \ker \bbeta = \{0\}$, the statement follows.
\end{proof}

\subsection{\texorpdfstring{One condition for Property \eqref{eq:W1}}{One condition for Property W1}}\label{s:W1_general}
In this subsection, we study the consequences of assuming that the weights $\beta_i$ are pairwise distinct and that the action of $\mathbb T$ has a scaling subtorus. More precisely, we derive a sufficient condition for Property~\eqref{eq:W1} and show that $\mathcal{S}(\pi_0) \subseteq \mathcal{W}_2\setminus \mathcal{W}_1$.

We say that a subtorus $\mathbb{T}_0$ of $\mathbb{T}$ is \emph{scaling} if the restrictions of all weights $\beta_i$ to $\mathbb{T}_0$ coincide and are nontrivial.

\begin{lemma}\label{l:T1}
    Assume that the weights $\beta_i$ are pairwise distinct and that the action of $\mathbb{T}$ on $\mathbb{C}^n$ has a scaling subtorus. Then $\ker\bbeta \cap \mathcal{W}_1 = \{0\}$.
\end{lemma}
\begin{proof}
    Since there is a scaling subtorus $\mathbb{T}_0$, there exists $\rho^\vee \in X_*(\mathbb{T}_0)$ such that, for all $i,j \in [1,n]$, $\langle \rho^\vee,\beta_i\rangle = \langle \rho^\vee,\beta_j\rangle \neq 0$. For any $\bw \in \ker \bbeta \cap \mathcal{W}_1$, pairing the equation $\bbeta \bw = 0$ with $\rho^\vee$, we obtain $\sum_{k=1}^n \bw_k = 0$. If $\bw \neq 0$, then, since $\bw \in \mathcal{W}_1$, there exist distinct $k,l \in [1,n]$ such that $\bw = e_k - e_l$. Then the equation $\bbeta \bw = 0$ implies $\beta_k = \beta_l$, which contradicts the assumption that the weights are pairwise distinct. Thus $\ker \bbeta \cap \mathcal{W}_1 = \{0\}$.
\end{proof}

\begin{proposition}\label{prop:SisLM}
     Assume that the weights $\beta_i$ are pairwise distinct and that the action of $\mathbb{T}$ on $\mathbb{C}^n$ has a scaling subtorus. Then
     \begin{equation}\label{eq:sw2}
        \mathcal{S}(\pi_0)  = \{\btheta \in (\mathbb{Z}_{\geq -1})^n\mid \,|J_{\btheta}|=2, \, \blambda \btheta \in (\mathbb{C}^{J_{\btheta}})_0, \, \bbeta 
        \btheta = 0\}.
     \end{equation}
\end{proposition}
\begin{proof}
    Since $\mathcal{S}(\pi_0) \subseteq \ker \bbeta \cap \mathcal{W}_2$, Lemma~\ref{l:T1} implies that $|J_{\btheta}| = 2$ for any $\btheta \in \mathcal{S}(\pi_0)$.
\end{proof}

\begin{proposition}\label{p:w1check}
   Assume that the weights $\beta_i$ are pairwise distinct and that the action of $\mathbb{T}$ on $\mathbb{C}^n$ has a scaling subtorus. Then for any $\mathcal{S}\subseteq \mathcal{S}(\pi_0)$ such that $0 \notin \mathcal{S}_{\geq 1}$, $\mathcal{S}$ satisfies Property~\eqref{eq:W1}.
\end{proposition}
\begin{proof}
    Since $\mathcal{S}_{\geq 1} \subseteq \ker \bbeta$ and $0 \notin \mathcal{S}_{\geq 1}$,  Lemma~\ref{l:T1} implies that $\mathcal{S}_{\geq 1} \cap \mathcal{W}_1 = \emptyset$. Thus $\mathcal{S}$ satisfies Property~\eqref{eq:W1}. 
\end{proof}

\begin{remark}
    It was shown in \cite[Section 2.2]{lm-deform} that equation~\eqref{eq:sw2} holds under the assumption that $\pi_0$ is $\mathbb{T}$-log-symplectic. When the action of $\mathbb{T}$ satisfies the conditions of Lemma~\ref{l:T1}, equation~\eqref{eq:sw2} holds as well. In particular, we will show in Proposition~\ref{p:difw} that it holds for the log-canonical part of a BD bivector on $\Mat_{n\times m}(\mathbb{C})$.
\end{remark}

\section{Belavin--Drinfeld bivectors and Poisson deformations}\label{s:main}
In this section, given a pair of Belavin--Drinfeld quadruples $(\bbg^r,\bbg^c)$, we study $\mathbb{T}_{(\bg^r,\bg^c)}$-invariant Poisson deformation theory of the log-canonical part $\pi_0:=(\pi_{(\bbg^r,\bbg^c)})_0$ of the BD bivector $\pi_{(\bbg^r,\bbg^c)}$ on $\Mat_{n\times m}(\mathbb{C})$. Our main objective is to prove the first part of Theorem~\ref{theorem}. We start by computing an explicit formula for $\pi_0$ and identifying a distinguished subset of weights $\mathcal{S}^\prime \subseteq \mathcal{S}(\pi_0)$. We will prove that $\mathcal{S}^\prime$ forms a basis of $\ker \bbeta$, where $\bbeta$ denotes a collection of $\mathbb{T}_{(\bg^r,\bg^c)}$-weights of the coordinate functions $x_{ij}$. This will allow us to apply the results of Section~\ref{s:cycles}. In particular, we will obtain that $\pi_0$ is $\mathbb{T}_{(\bg^r,\bg^c)}$-log-symplectic if and only if the smoothing graph $\Gamma(\mathcal{S}^\prime)$ is acyclic, and that, in this case, $\mathcal{S}^\prime = \mathcal{S}(\pi_0)$. 

We will then prove Property~\eqref{eq:W1} for $\pi_0$ using the results of Section~\ref{s:W1_general}. Finally, in Subsection~\ref{s:pi0_piBD}, assuming that $\pi_0$ is $\mathbb{T}_{(\bg^r,\bg^c)}$-log-symplectic, we will show that $\pi_{(\bbg^r,\bbg^c)}$ is the maximal normalized admissible $\mathbb{T}_{(\bg^r,\bg^c)}$-invariant algebraic Poisson deformation of $\pi_0$.

We discuss the notion of aperiodicity of $(\bg^r,\bg^c)$ in Subsection~\ref{s:meander}. There, we give a combinatorial criterion for aperiodicity in terms of the actions of $\gamma_r$ and $\gamma_c$ on simple roots. We will prove that $(\bg^r,\bg^c)$ is aperiodic if and only if $\pi_0$ is $\mathbb{T}_{(\bg^r,\bg^c)}$-log-symplectic. The notion of aperiodicity will be crucial in Section~\ref{s:birat}, where we show that if $(\bg^r,\bg^c)$ is aperiodic, then there exists a $\mathbb{T}_{(\bg^r,\bg^c)}$-equivariant Poisson birational map between $(\Mat_{n\times m}(\mathbb{C}),\pi_0)$ and $(\Mat_{n\times m}(\mathbb{C}),\pi_{(\bbg^r,\bbg^c)})$.

\subsection{\texorpdfstring{The log-canonical part $\pi_0$}{The log-canonical part}}
In this subsection, for a given pair $(\bbg^r,\bbg^c)$ of BD quadruples, we compute the log-canonical part $\pi_0$ of the Poisson bivector $\pi_{(\bbg^r,\bbg^c)}$ on $\Mat_{n\times m}(\mathbb{C})$.

We denote by $\blambda$ the matrix of coefficients of $\pi_0$. Note that the rows and columns of $\blambda$ are labeled by pairs of indices, so the entries of $\blambda$ are $\lambda_{ij,kl}$, $i,k \in [1,n]$, $j,l \in [1,m]$. The Poisson bivector corresponding to the pair of trivial BD quadruples is denoted by $\pi_{\std}$.

\begin{proposition}\label{p:log-part}
    For any pair of BD quadruples $\bbg^r:=(\bg^r,r_0^r)$ and $\bbg^c := (\bg^c,r_0^c)$, the log-canonical part $\pi_0$ of the BD bivector $\pi_{(\bbg^r,\bbg^c)}$ is given by
    \begin{equation*}
        \pi_0 = (\pi_{\std})_0 + (s_0^c)^L - (s_0^r)^R.
    \end{equation*}
    Correspondingly, the matrix of coefficients of $\pi_0$ is given by
    \begin{equation*}
        \lambda_{ij,kl} = (\lambda_{\std})_{ij,kl} + \langle S_0^c e_{jj}, e_{ll}\rangle - \langle S_0^r e_{ii},e_{kk}\rangle,
    \end{equation*}
    where
    \begin{equation*}
        (\lambda_{\std})_{ij,kl} = \frac{1}{2}\big(\delta_{ik}\sgn(l-j) + \delta_{jl}\sgn(k-i)\big).
    \end{equation*}
\end{proposition}
\begin{proof}
    Let us write $\{\cdot,\cdot\}$ for the Poisson bracket $\{\cdot,\cdot\}_{(\bbg^r,\bbg^c)}$. Fix $i,k \in [1,n]$ and $j,l \in [1,m]$, and assume that $(i,j)\neq (k,l)$. We will compute the Poisson bracket $\{x_{ij},x_{kl}\}$ and extract the term proportional to $x_{ij}x_{kl}$. 
    Let $I$ be the ideal of $\mathbb{C}[\Mat_{n\times m}]$ generated by all $x_{ab}$ for $(a,b) \notin \{(i,j),(k,l)\}$. For $f,g \in \mathbb{C}[\Mat_{n\times m}]$, write $f\sim g$ if $f-g \in I$. We see that
   \begin{equation}\label{eq:pf_brack}
       \{x_{ij},x_{kl}\} = \langle R^c(\nabla^L x_{ij}), \nabla^L x_{kl} \rangle - \langle R^r(\nabla^R x_{ij}), \nabla^R x_{kl} \rangle.
   \end{equation}
    Note that
    \begin{equation*}
        \nabla^L x_{ij} = \sum_{\mu=1}^m x_{i\mu} e_{j\mu}, \ \ \nabla^R x_{ij} = \sum_{\mu=1}^n x_{\mu j} e_{\mu i};
    \end{equation*}
   hence, for the first term in equation~\eqref{eq:pf_brack}, we obtain
   \begin{equation*}\begin{split}
       \langle R^c(\nabla^L x_{ij}),\nabla^L x_{kl} \rangle = &\sum_{\mu_1=1}^m\sum_{\mu_2=1}^m x_{i\mu_1}x_{k\mu_2} \langle R^c(e_{j\mu_1}),e_{l\mu_2}\rangle\sim\\ &\sim x_{ij}x_{kl}\big(\langle R^c(e_{jj}),e_{ll}\rangle + \delta_{ik}\langle R^c(e_{jl}),e_{lj}\rangle \big) = \\ &= x_{ij}x_{kl}\big( \langle R_0^c(e_{jj}),e_{ll}\rangle + \delta_{ik}[j<l] \big).
       \end{split}
   \end{equation*}
   Likewise, for the second term,
   \begin{equation*}\begin{split}
       \langle R^r(\nabla^R x_{ij}),\nabla^R x_{kl} \rangle = &\sum_{\mu_1=1}^n\sum_{\mu_2=1}^n x_{\mu_1 j} x_{\mu_2 l} \langle R^r(e_{\mu_1 i}),e_{\mu_2 k}\rangle \sim \\ &\sim x_{ij}x_{kl} \big(\langle R^r(e_{ii}),e_{kk}\rangle + \delta_{jl} \langle R^r(e_{ki}),e_{ik}\rangle \big) = \\ &=x_{ij}x_{kl} \big( \langle R_0^r(e_{ii}),e_{kk}\rangle + \delta_{jl} [k < i] \big).
   \end{split}
   \end{equation*}
   Combining these expressions, we obtain
    \begin{equation*}
        \begin{split}
            \{x_{ij},x_{kl}\} \sim x_{ij}x_{kl} \big(\langle R_0^c(e_{jj}),e_{ll}\rangle-\langle R_0^r(e_{ii}),e_{kk}\rangle  + \delta_{ik}[j<l] - \delta_{jl} [k < i]\big).
        \end{split}
    \end{equation*}
    This expression holds for any pair of BD quadruples (in particular, for the standard one). Subtracting $\{x_{ij},x_{kl}\}_{\std}$ from $\{x_{ij},x_{kl}\}$, we obtain that
    \begin{equation*}
        \{x_{ij},x_{kl}\} - \{x_{ij},x_{kl}\}_{\std} \sim x_{ij} x_{kl}\langle S_0^c(e_{jj}),e_{ll}\rangle-x_{ij} x_{kl}\langle S_0^r(e_{ii}),e_{kk}\rangle.
    \end{equation*}
    Passing from this bracket identity to the corresponding bivector identity, we obtain the statement of the proposition.
\end{proof}
\subsection{Solutions of the LM equations}\label{s:sols_lm}
In this subsection, given any BD pair $(\bg^r,\bg^c)$, we describe a subset $\mathcal{S}^\prime \subseteq \mathcal{S}(\pi_0)$ of solutions to the LM equations and prove that $\mathcal{S}^\prime$ is a basis of $\ker \bbeta$. We then prove that if $\pi_0$ is $\mathbb{T}_{(\bg^r,\bg^c)}$-log-symplectic, then $\mathcal{S}^\prime = \mathcal{S}(\pi_0)$.

For a given pair $(\bg^r,\bg^c)$ of BD triples, define the following $n \times m$ matrices:
\begin{align*}
    & \btheta_{pq} := -e_{pq} + e_{p+1,q} + e_{p,q+1} - e_{p+1,q+1}, & &p\in [1,n-1], \ q \in [1,m-1];\\
    & \btheta_{p}^r:= e_{p,1} - e_{p+1,1} - e_{q,m} + e_{q+1,m}, & & \alpha_p^r \in \Gamma_1^r, \ \alpha_q^r=\gamma_r(\alpha_p^r) \in \Gamma_2^r;\\
    & \btheta_{p}^c:= -e_{n,p} + e_{n,p+1} + e_{1,q} - e_{1,q+1},& &\alpha_p^c \in \Gamma_1^c, \ \alpha_q^c=\gamma_c(\alpha_p^c) \in \Gamma_2^c.
\end{align*}
Let $\mathcal{S}^\prime:= \mathcal{S}(\bg^r,\bg^c)$ be the set consisting of the above matrices; in other words,
\begin{equation*}
    \mathcal{S}^\prime:= \{\btheta_{pq} \mid p\in [1,n-1], \ q \in [1,m-1]\} \cup \{\btheta_{p}^r \mid \alpha_p^r \in \Gamma_1^r\} \cup \{\btheta_{p}^c \mid \alpha_p^c \in \Gamma_1^c\}.
\end{equation*}
In the proofs, we will need alternative expressions for the elements of $\mathcal{S}^\prime$. Specifically, let $e_1,\ldots,e_n$ be the standard basis of $\mathbb{R}^n$, and let $f_1,\ldots,f_m$ be the standard basis of $\mathbb{R}^m$. Then each  $e_{ij}$ can be written as $e_{ij} = e_i f_j^T$. Hence, 
\begin{align*}
    &\btheta_{pq} = -(e_p-e_{p+1})(f_q-f_{q+1})^T, & & p \in [1,n-1], \ q \in [1,m-1];\\
    &\btheta_p^r = (e_p-e_{p+1})f_1^T - (e_q-e_{q+1})f_m^T, & & \alpha_p^r \in \Gamma_1^r, \ \alpha_q^r = \gamma_r(\alpha_p^r);\\
    &\btheta_p^c = -e_n(f_p-f_{p+1})^T + e_1(f_q-f_{q+1})^T, & & \alpha_p^c \in \Gamma_1^c, \ \alpha_q^c = \gamma_c(\alpha_p^c).
\end{align*}
Setting $h_i^r := e_i - e_{i+1}$ and $h_j^c := f_j - f_{j+1}$, $i \in [1,n-1]$, $j \in [1,m-1]$, we obtain that
\begin{align*}
    &\btheta_{pq} = -h_p^r(h_q^c)^T, & & p \in [1,n-1], \ q \in [1,m-1];\\
    &\btheta_p^r = h_p^rf_1^T - h_q^rf_m^T, & & \alpha_p^r \in \Gamma_1^r, \ \alpha_q^r = \gamma_r(\alpha_p^r);\\
    &\btheta_p^c = -e_n(h_p^c)^T + e_1(h_q^c)^T, & & \alpha_p^c \in \Gamma_1^c, \ \alpha_q^c = \gamma_c(\alpha_p^c).
\end{align*}

Given the matrix of coefficients $\blambda$ of $\pi_0$, we denote by $\blambda_{-,kl}$ the column of $\blambda$ indexed by $(k,l)$.
 
\begin{proposition}\label{p:bd_lm}
    For any $(\bg^r,\bg^c)$, $\mathcal{S}^\prime\subseteq \mathcal{S}(\pi_0,\mathbb{T})$. Moreover, the elements of $\mathcal{S}^\prime$ satisfy the following equations:
    \begin{align*}
        &\sum_{k=1}^n\sum_{l=1}^m\blambda_{-,kl} (\btheta_{pq})_{kl} = e_{pq}-e_{p+1,q+1}, & & p \in [1,n-1], \ q \in [1,m-1];\\
        &\sum_{k=1}^n\sum_{l=1}^m\blambda_{-,kl}( \btheta_{p}^r)_{kl} = e_{qm} - e_{p+1,1} , & & \alpha_p^r \in \Gamma_1^r, \ \alpha_q^r= \gamma_r(\alpha_p^r) \in \Gamma_2^r\\
        &\sum_{k=1}^n\sum_{l=1}^m\blambda_{-,kl}( \btheta_{p}^c)_{kl} = e_{np}-e_{1,q+1}, & & \alpha_p^c \in \Gamma_1^c, \ \alpha_q^c = \gamma_c(\alpha_p^c) \in \Gamma_2^c.
    \end{align*}
\end{proposition}
\begin{proof}
    Let us set $\lambda_{ij,kl}^c := \langle S_0^c e_{jj}, e_{ll}\rangle$ and $\lambda_{ij,kl}^r := \langle S_0^r e_{ii},e_{kk}\rangle$. By Proposition~\ref{p:log-part},
    \begin{equation*}
        \lambda_{ij,kl} = (\lambda_{\std})_{ij,kl} + \lambda_{ij,kl}^c-\lambda_{ij,kl}^r.
    \end{equation*}
    The computations are lengthy but straightforward. Let us outline the key identities. First, observe that 
    \begin{equation}\label{eq:sgneq}
        \sgn(p+1-i) -\sgn(p-i) = \delta_{ip}+\delta_{i,p+1}
    \end{equation}
    Second, by equation~\eqref{eq:S0},
    \begin{equation}\label{eq:s0r_eq}
        S_0^r(h_p-h_q) = \frac{1}{2}(h_p+h_q)
    \end{equation}
    where $h_p = e_{pp}-e_{p+1,p+1}$ and $\alpha_q = \gamma_r(\alpha_p)$. Using identity~\eqref{eq:sgneq}, one can verify that
    \begin{equation*}
        \sum_{k=1}^n\sum_{l=1}^m (\lambda_{\std})_{ij,kl} (\btheta_{pq})_{kl} = \delta_{ip}\delta_{jq} - \delta_{i,p+1}\delta_{j,q+1};
    \end{equation*}
    \begin{equation*}
        \sum_{k=1}^n\sum_{l=1}^m \lambda_{ij,kl}^c (\btheta_{pq})_{kl} = \sum_{k=1}^n\sum_{l=1}^m \lambda_{ij,kl}^r (\btheta_{pq})_{kl} = 0.
    \end{equation*}
    Hence the first equation holds. To show the second equation, we compute
    \begin{equation*}
        \sum_{k=1}^n \sum_{l=1}^m (\lambda_{\std})_{ij,kl}(\btheta_p^r)_{kl} = \frac{1}{2}( \delta_{i,p+1}-\delta_{ip}) + \frac{1}{2}(-\delta_{iq} + \delta_{i,q+1}) - \delta_{j,1}\delta_{p+1,i} + \delta_{iq}\delta_{jm};
    \end{equation*}
    \begin{equation*}
        \sum_{k=1}^n \sum_{l=1}^m(-\lambda^r_{ij,kl})(\btheta_p^r)_{kl} = \frac{1}{2}(\delta_{ip} - \delta_{i,p+1}) + \frac{1}{2}(\delta_{qi} - \delta_{q+1,i}),
    \end{equation*}
    where we used~\eqref{eq:s0r_eq}; and lastly,
    \begin{equation*}
        \sum_{k=1}^n \sum_{l=1}^m \lambda^c_{ij,kl} (\btheta_p^r)_{kl} = 0.
    \end{equation*}
    Combining these identities, we obtain the second equation. The last equation can be proved in a similar way.

    Now let us verify that the elements in $\mathcal{S}^\prime$ belong to $\ker \bbeta$. More precisely, we verify the following equations:
    \begin{align*}
        & \sum_{k=1}^{n}\sum_{l=1}^m (\beta_{\std})_{kl}(\btheta_{pq})_{kl} = 0, & &p \in [1,n-1], \ q \in [1,m-1];\\
        & \sum_{k=1}^n \sum_{l=1}^m (\beta_{\std})_{kl}  (\btheta_p^r)_{kl} = \alpha^r_p - \alpha^r_q \in (\mathfrak{h}^r)^0, & &\alpha_p^r \in \Gamma_1^r,\ \alpha_q^r = \gamma_r(\alpha_p^r)\in \Gamma_2^r;\\
        &\sum_{k=1}^n \sum_{l=1}^m (\beta_{\std})_{kl}  (\btheta_p^c)_{kl} = \alpha^c_q - \alpha^c_p \in (\mathfrak{h}^c)^0, & &\alpha_p^c \in \Gamma_1^c,\ \alpha_q^c = \gamma_c(\alpha_p^c) \in \Gamma_2^c.
    \end{align*}
    Here, for a subspace $\mathfrak{h}^\prime$ of $\mathfrak{h}$, $(\mathfrak{h}^\prime)^0$ is the subspace of $\mathfrak{h}^*$ of linear functions vanishing on $\mathfrak{h}^\prime$. 
    Recall from Subsection~\ref{s:bd_background} that the $\mathbb{T}_{\std}:=\mathbb{T}_{(\bg_{\std},\bg_{\std})}$-weight $(\beta_{\std})_{kl}$ of a matrix entry $x_{kl}$ is given by
    \begin{equation*}
        (\beta_{\std})_{kl}= \beta_k^r + e + \beta_l^c \in X^*(\mathbb{T}_{\std}) =  X^*(\mathcal{H}_{\ses}^{(n)})\oplus \mathbb{Z}\langle e\rangle \oplus X^*(\mathcal{H}_{\ses}^{(m)})
    \end{equation*}
    where 
     $\beta_{k}^r = \varpi_k^r - \varpi_{k-1}^r \in X^*(\mathcal{H}_{\ses}^{(n)})$ and $\beta_{l}^c = \varpi_l^c - \varpi_{l-1}^c \in X^*(\mathcal{H}_{\ses}^{(m)})$ (here, $\varpi_{0}^r = \varpi_n^r = 0$ and $\varpi_0^c =\varpi_m^c = 0$). The projections
    \begin{equation*}
        \beta_{kl} = \pr_r(\beta_{k}^r) + e + \pr_c(\beta_{l}^c) \in X^*(\mathbb{T}_{(\bg^r,\bg^c)}) = X^*(\mathcal{H}_{\bg^r})\oplus \mathbb{Z}\langle e\rangle \oplus X^*(\mathcal{H}_{\bg^c})
    \end{equation*} 
    are the $\mathbb{T}_{(\bg^r,\bg^c)}$-weights of $x_{kl}$.
    
    We see that
    \begin{equation*}\begin{split}
        \sum_{k=1}^{n}\sum_{l=1}^m \beta_{k}^r (\btheta_{pq})_{kl} &= \beta_{p}^r (\btheta_{pq})_{pq} +  \beta_{p+1}^r (\btheta_{pq})_{p+1,q} + \beta_{p}^r (\btheta_{pq})_{p,q+1} + \beta_{p+1}^r (\btheta_{pq})_{p+1,q+1} = 0.
        \end{split}
    \end{equation*}
    The same holds for the corresponding sums involving $e$ and $\beta_l^c$. Therefore, $\btheta_{pq} \in \ker \bbeta$. For the elements $\btheta_p^r$, we see that
    \begin{equation*}\begin{split}
        \sum_{k=1}^n \sum_{l=1}^m \beta_{k}^r (\btheta_{p}^r)_{kl} &= \beta_{p}^r - \beta_{p+1}^r - \beta_{q}^r + \beta_{q+1}^r = \\ &= (2\varpi_p^r - \varpi_{p-1}^r - \varpi_{p+1}^r) - (2\varpi_q^r - \varpi_{q-1}^r - \varpi_{q+1}^r) = \alpha_p^r - \alpha_q^r,
        \end{split}
    \end{equation*}
    which vanishes under the projection $\pr_r$.
    Likewise, one computes that 
    \begin{equation*}
        \sum_{k=1}^n \sum_{l=1}^m \beta_{l}^c (\btheta_p^r)_{kl} = 0, \qquad \sum_{k=1}^n \sum_{l=1}^m e (\btheta_p^r)_{kl} = 0.
    \end{equation*}
    Combining these identities, we obtain that $\btheta_p^r \in \ker \bbeta$. In a similar way, we see that $\btheta_p^c \in \ker \bbeta$.
\end{proof}

\begin{proposition}\label{p:kerb_basis}
    For any $(\bg^r,\bg^c)$, the set $\mathcal{S}^\prime$ is a basis of $\ker \bbeta$.
\end{proposition}
\begin{proof}
    Notice that 
    \begin{equation*}
        \dim \ker \bbeta = nm - \dim \mathbb{T}_{(\bg^r,\bg^c)} = (n-1)(m-1) + |\Gamma_1^r| + |\Gamma_1^c| = |\mathcal{S}^\prime|;
    \end{equation*}
    hence, it suffices to show that the set $\mathcal{S}^\prime$ is linearly independent. 
    
    Let us first show that the set of $\btheta_{pq}$, for $p\in [1,n-1]$ and $q \in [1,m-1]$, is linearly independent. Suppose that
    \begin{equation*}
        \sum_{p=1}^{n-1}\sum_{q = 1}^{m-1} a_{pq} \btheta_{pq} = 0
    \end{equation*}
    for some $a_{pq} \in \mathbb{C}$. For any $k \in [1,n-1]$, multiplying on the left by $-(e_1 + \cdots + e_k)^T$ gives
    \begin{equation*}
        \sum_{q=1}^{m-1} a_{kq}(h_q^c)^T = 0,
    \end{equation*}
    hence $a_{kq} = 0$ for each $q \in [1,m-1]$. We conclude that $\btheta_{pq}$ are linearly independent.

    Now let us turn to the general case. Consider a linear combination
    \begin{equation}\label{eq:lin_comb}
        \sum_{p=1}^{n-1}\sum_{q = 1}^{m-1} a_{pq} \btheta_{pq} + \sum_{\alpha_p^r \in \Gamma_1^r}b_p\btheta_p^r + \sum_{\alpha_p^c \in \Gamma_1^c} c_p\btheta_p^c = 0.
    \end{equation}
    Multiplying the above equation on the right by $f_1 + \cdots + f_m$, we arrive at the equation
    \begin{equation*}
        \sum_{\alpha_p^r \in \Gamma_1^r}b_p(h_p^r - h_{\gamma_r(p)}^r) = 0.
    \end{equation*}
    Let $\mathcal{C}(\bg^r)$ be the set of BD chains corresponding to $\bg^r$. The above sum can be rewritten as
    \begin{equation*}
        \sum_{\alpha_p^r \in \Gamma_1^r}b_p(h_p^r - h_{\gamma_r(p)}^r) = \sum_{C \in \mathcal{C}(\bg^r)}\sum_{\alpha_p^r \in \Gamma_1^r \cap C}b_p(h_p^r - h_{\gamma_r(p)}^r) = 0,
    \end{equation*}
    and since the chains do not intersect, for every $C \in \mathcal{C}(\bg^r)$,
    \begin{equation*}
        \sum_{\alpha_p^r \in \Gamma_1^r \cap C}b_p(h_p^r - h_{\gamma_r(p)}^r) = 0.
    \end{equation*}
    Now let us write $C= \{\alpha_{p_1}^r,\alpha_{p_2}^r,\ldots,\alpha_{p_k}^r\}$, where $\gamma_r(\alpha_{p_i}^r) = \alpha_{p_{i+1}}^r$. We see that
    \begin{equation*}\begin{split}
        \sum_{i=1}^{k-1}& b_{p_i}(h_{p_i}^r - h_{\gamma_r(p_i)}^r) =\\&= b_{p_1}h_{p_1}^r + (b_{p_2}-b_{p_1})h_{p_2}^r + \cdots +(-b_{p_{k-1}})h_{p_k}^r = 0,
        \end{split}
    \end{equation*}
    hence $b_{p_1} = b_{p_2} = \cdots = b_{p_{k-1}} = 0$. We conclude that for all $p \in \Gamma_1^r$, $b_p = 0$. To show that the coefficients $c_p$ vanish, we apply the same argument, now multiplying equation~\eqref{eq:lin_comb} by $(e_1+\cdots + e_n)^T$ on the left in order to isolate the sum over $\Gamma_1^c$. Finally, since $b_p=c_p=0$ and the elements $\btheta_{pq}$ are linearly independent, we conclude that $a_{pq}=0$ for all $p,q$. Thus $\mathcal{S}^\prime$ is a basis of $\ker \bbeta$.
\end{proof}
\begin{corollary}
    For any $(\bg^r,\bg^c)$, the bivector $\pi_0$ is $\mathbb{T}_{(\bg^r,\bg^c)}$-log-symplectic if and only if $\Gamma(\mathcal{S}^\prime)$ has no cycles.
\end{corollary}
\begin{proof}
    This follows from Corollary~\ref{c:cycles_tlogs}.
\end{proof}

\begin{proposition}\label{p:sp_eqls_sp0}
    For any $(\bg^r,\bg^c)$, if $\pi_0$ is $\mathbb{T}_{(\bg^r,\bg^c)}$-log-symplectic, then $\mathcal{S}^\prime = \mathcal{S}(\pi_0)$.
\end{proposition}
\begin{proof}
    Let $\Gamma(\mathcal{S}^\prime)$ be the smoothing graph associated with $\mathcal{S}^\prime$. Let $V:=[1,n]\times [1,m]$ be its vertex set, and let $E$ be its edge set. For each edge $\varepsilon\in E$, let $\btheta^{\varepsilon} \in \mathcal{S}^\prime$ denote the corresponding element. Set $V^\circ :=[2,n-1]\times[2,m-1]$ and $\partial V:= V \setminus V^{\circ}$.
    Throughout the proof, the degree of a vertex means its degree in the smoothing graph $\Gamma(\mathcal{S}^\prime)$.
    Recall from Proposition~\ref{p:lm_bunch} that each vertex of $V$ has degree at most $2$. Clearly, each vertex of $V^{\circ}$ has degree $2$, for if $(p,q) \in V^{\circ}$, then $\btheta_{p-1,q-1},\btheta_{pq}\in \mathcal{S}^\prime$. The vertices in $\partial V$ may have degree $0$, $1$, or $2$. Denote by $(\partial V)_1$ and $(\partial V)_2$ the subsets of boundary vertices of degrees $1$ and $2$, respectively.
    
    Suppose, for contradiction, that there exists $\btheta^\prime \in \mathcal{S}(\pi_0)\setminus \mathcal{S}^\prime$, and let $\varepsilon^\prime$ be the corresponding edge in $\Gamma(\pi_0,\mathbb{T})$. Since $\pi_0$ is $\mathbb{T}$-log-symplectic and $\mathcal{S}^\prime$ is a basis of $\ker \bbeta$, the smoothing graph $\Gamma(\mathcal{S}^\prime\cup\{\btheta^\prime\})$ contains a cycle $c$ passing through $\varepsilon^\prime$. We orient the cycle so that $a^{\varepsilon} = 1$ for the edges $\varepsilon$ of $c$ belonging in $\Gamma(\mathcal{S}^\prime)$ (see Section~\ref{s:cycles} for orientation convention and Proposition~\ref{p:bd_lm} for the LM equations). Let $V(c)$ and $E(c)$ denote the vertex set and edge set of $c$, respectively. By Proposition~\ref{p:ker_cycles}, we obtain a linear relation
    \begin{equation}\label{eq:cycle2}
        \btheta(c) = (a^{\varepsilon^\prime})^{-1}\btheta^\prime + \sum_{\varepsilon \in E(c)\setminus\{\varepsilon^\prime\}} \btheta^\varepsilon = 0.
    \end{equation}
    Since $\btheta^\prime$ has exactly $2$ negative entries, and both are equal to $-1$, from~\eqref{eq:cycle2} we see that $(a^{\varepsilon^\prime})^{-1} \in \mathbb{Z}$, and moreover, by Proposition~\ref{p:lm_bunch} (or by \cite[Lemma 2.2.9]{lm-deform}), $(a^{\varepsilon^\prime})^{-1}$ is a positive integer. Define
    \begin{equation}\label{eq:sum_ecP}
        P:=\sum_{\varepsilon \in E(c)\setminus\{\varepsilon^\prime\}} \btheta^\varepsilon
    \end{equation}
    It follows that 
    \begin{equation}\label{eq:P}
        P_{ij} > 0 \ \Leftrightarrow \ (i,j) \in (\partial V)_1 \cap V(c).
    \end{equation}
    Let us observe the following statements:
    \begin{itemize}
        \item If $(1,1) \in V(c)$, then $(1,1)$ has degree $1$, $\alpha_1^r \in \Gamma_1^r$ and $\alpha_1^c \in \Gamma_2^c$, and $(1,2),(2,1) \in V(c)$. Indeed, $(1,1)$ has degree $1$ in $\Gamma(\mathcal{S}^\prime)$ for any pair $(\bg^r,\bg^c)$. Therefore, if $(1,1) \in V(c)$, equation~\eqref{eq:P} implies that $P_{11} > 0$. The edge $(1,1)\edge (2,2)$ contributes $\btheta_{11}$ to the right-hand side of~\eqref{eq:sum_ecP}, and $\btheta_{11}$ in turn contributes $-e_{11}$. Since $P_{11} > 0$, the other terms must contribute at least $2e_{11}$. This can only happen if $\alpha_1^r \in \Gamma_1^r$ and $\alpha_1^c \in \Gamma_2^c$, which also implies that $(1,2),(2,1) \in V(c)$ (in this case, $a^{\varepsilon^\prime}= 1$).
        \item For any vertex $(1,q) \in V(c) \cap (\partial V)_2$ with $q \in [2,m-1]$, we see that $(1,q-1),(1,q+1) \in V(c)$, and $\alpha_{q-1}^c \in \Gamma_2^c$. Indeed, in this case, the cycle $c$ passes through $(n,p) \edge (1,q) \edge (2,q+1)$, where $q-1 = \gamma_c(p)$ for some $p \in \Gamma_1^c$. The term $\btheta_{1,q}$ in the sum on the right-hand side of~\eqref{eq:sum_ecP} contributes $e_{1,q+1}$. If another term $\btheta^\varepsilon$, $\varepsilon \in E(c)\setminus\{\varepsilon^\prime,\varepsilon_{1,q}\}$, contributes $-e_{1,q+1}$, then $(1,q+1) \in V(c)$. Otherwise, $P_{1,q+1} > 0$, and equation~\eqref{eq:P} still implies that $(1,q+1) \in V(c)$. Likewise, $\btheta_p^c$ contributes $e_{1,q-1}$, and the same reasoning shows that $(1,q-1) \in V(c)$.
        \item For any vertex $(1,q) \in V(c) \cap (\partial V)_1$ with $q \in [2,m-1]$, we see that $\alpha_{q-1}^c \notin \Gamma_2^c$, $\alpha_q^c \in \Gamma_2^c$, and $(1,q-1),(1,q+1) \in V(c)$. Indeed, by equation~\eqref{eq:P}, $P_{1,q} > 0$,  and $c$ passes through $(1,q) \edge (2,q+1)$, which contributes $-e_{1,q}$ to~\eqref{eq:sum_ecP}. Therefore, two additional terms must each contribute $e_{1,q}$. The only possibility is that the sum in~\eqref{eq:sum_ecP} has a contribution of $\btheta_{p}^c$ with $q =\gamma_c(p)$ and a contribution of $\btheta_{1,q-1}$, but no contribution of $\btheta_{p^\prime}^c$ with $q -1 = \gamma_c(p^\prime)$. This gives the statement of this item.
        \item If $(1,m) \in V(c)$, then $(1,m)$ has degree $2$, $(1,m-1), (2,m) \in V(c)$, $\alpha_{m-1}^c \in \Gamma_2^c$, and $\alpha_1^r \in \Gamma_2^r$.
    \end{itemize}
    As a consequence, if $(1,q) \in V(c)$ for some $q \in [1,m]$, we obtain that $(1,q) \in V(c)$ for all $q \in [1,m]$. A similar argument shows that if $(p,1) \in V(c)$ for some $p \in [1,n]$, then $(p,1) \in V(c)$ for all $p \in [1,n]$. Moreover, if either $(1,q) \in V(c)$ for some $q$ or $(p,1) \in V(c)$ for some $p$, then $(1,1) \in V(c)$, and hence all vertices $(1,q)$ and $(p,1)$ belong to $V(c)$.

    Now observe that $V(c)$ contains at least $3$ vertices of degree $1$. Indeed, $(1,1)$ always has degree $1$. If every vertex $(1,q)$ for $q \in [2,m]$ had degree $2$, then $\Gamma_2^c = [1,m-1]$, contradicting the nilpotency of $\gamma_c$. Therefore, there exists $q^\prime \in [2,m]$ such that $(1,q^\prime)$ has degree $1$. Likewise, there exists $p^\prime \in [2,n]$ such that $(p^\prime,1)$ has degree $1$. Since there can be at most $2$ vertices of degree $1$, we obtain a contradiction.
\end{proof}    

\subsection{\texorpdfstring{Property~\eqref{eq:W1}}{Property W1}}
In this subsection, we prove that the set $\mathcal{S}^\prime$ defined in the previous subsection satisfies Property~\eqref{eq:W1}.

\begin{proposition}\label{p:difw}
    For any BD pair $(\bg^r,\bg^c)$, the weights $\beta_{kl}$ are pairwise distinct.
\end{proposition}
\begin{proof}
    Clearly, the weights $\beta_{k}^r$ and the weights $\beta_{l}^c$ belong to different direct summands of the character lattice. Recall that the pairing between the fundamental weights $\varpi_k$ and the fundamental coweights $\varpi_l^{\vee}$ is given by $\langle \varpi_k,\varpi_l^{\vee}\rangle = \min \{k,l\}$, and that $\beta_{k}^r = \varpi_k^r - \varpi_{k-1}^r$. Recall also that the Lie algebra $\mathfrak{h}_{\bg^r}$ is spanned by the elements $h_C$ associated with the BD chains $C$, where
    \begin{equation*}
        h_C = \sum_{i:\, \alpha_i \in C}\varpi_i^{r,\vee}.
    \end{equation*}
    Evaluating $\beta_k^r$ on $h_C$, we see that
    \begin{equation*}
        \beta_{k}^r(h_C) = \sum_{i:\, \alpha_i \in C} (\min\{k,i\}-\min\{k-1,i\}) = |\{i \in [1,n-1] \mid \alpha_i \in C, \ i \geq k\}|.
    \end{equation*}
    Hence, to verify that $\pr_r(\beta_{k}^r)$ and $\pr_r(\beta_{k^\prime}^r)$ are distinct, let $i := \min\{k,k^\prime\}$ and choose the BD chain $C$ containing~$\alpha_i$. We then have $\beta_{k}^r(h_C) \neq \beta_{k^\prime}^r(h_C)$. By the same argument, the weights associated with~$\bg^c$ are pairwise distinct. By Proposition~\ref{p:w1check}, we conclude that $\mathcal{S}^\prime$ satisfies Property~\eqref{eq:W1}.
\end{proof}

\begin{proposition}\label{p:BD_W1}
    For any BD pair $(\bbg^r,\bbg^c)$, let $\pi_0 := (\pi_{(\bbg^r,\bbg^c)})_0$. Then the set $\mathcal{S}^\prime$ satisfies Property~\eqref{eq:W1}. In particular, if $\pi_0$ is $\mathbb{T}_{(\bg^r,\bg^c)}$-log-symplectic, then $\mathcal{S}(\pi_0)$ satisfies Property~\eqref{eq:W1}.
\end{proposition}
\begin{proof}
    Recall that the BD torus $\mathbb{T}_{(\bg^r,\bg^c)}$ acting on $\Mat_{n\times m}(\mathbb{C})$ has a  $\mathbb{C}^\times$-factor acting by scalar multiplication, which is a scaling subtorus. By Proposition~\ref{p:difw}, the weights $\beta_{kl}$ are pairwise distinct. Moreover, since $\mathcal{S}^\prime$ is a basis of $\ker\bbeta$, it is linearly independent, and hence $0 \notin \mathcal{S}^\prime_{\geq 1}$. Thus all the assumptions of Proposition~\ref{p:w1check} are satisfied, and $\mathcal{S}^\prime$ satisfies Property~\eqref{eq:W1}. In particular, if $\pi_0$ is $\mathbb{T}_{(\bg^r,\bg^c)}$-log-symplectic, Proposition~\ref{p:sp_eqls_sp0} gives $\mathcal{S}^\prime = \mathcal{S}(\pi_0)$, so $\mathcal{S}(\pi_0)$ also satisfies Property~\eqref{eq:W1}.
\end{proof}

Proposition~\ref{p:difw} can also be used to give a finer description of $\mathcal{S}(\pi_0)$ for any $(\bg^r,\bg^c)$.
\begin{proposition}
    For any BD pair $(\bg^r,\bg^c)$, 
    \begin{equation*}
        \mathcal{S}(\pi_0)  = \{\btheta \in (\mathbb{Z}_{\geq -1})^{nm}\mid \,|J_{\btheta}|=2, \, \blambda \btheta \in (\mathbb{C}^{J_{\btheta}})_0, \, \bbeta 
        \btheta = 0\}.
    \end{equation*}
\end{proposition}
Here, we regard each nonzero $(\mathbb{C}^\times)^{nm}$-weight $\btheta \in \Mat_{n\times m}(\mathbb{Z}_{\geq -1})$ as a column vector.
\begin{proof}
    The statement follows from Proposition~\ref{p:difw} and Proposition~\ref{prop:SisLM}.
\end{proof}

\subsection{Belavin--Drinfeld bivectors as Poisson deformations}\label{s:pi0_piBD}

In this subsection, we prove that if $(\pi_{(\bbg^r,\bbg^c)})_0$ is $\mathbb{T}_{(\bg^r,\bg^c)}$-log-symplectic, then $\pi_{(\bbg^r,\bbg^c)}$ is the maximal normalized admissible $\mathbb{T}_{(\bg^r,\bg^c)}$-invariant algebraic Poisson deformation of $(\pi_{(\bbg^r,\bbg^c)})_0$.

We use the expressions for the elements of $\mathcal{S}^\prime$ given in Subsection~\ref{s:sols_lm}. We adopt the following notation:
\begin{itemize}
    \item If $\alpha = \alpha_p + \cdots + \alpha_{p+k-1} \in \Phi_+(A_{n-1})$, set $h_{\alpha}^r := e_p - e_{p+k}$;
    \item If $\alpha = \alpha_q + \cdots + \alpha_{q+k-1} \in \Phi_+(A_{m-1})$, set $h_{\alpha}^c := f_q - f_{q+k}$.
\end{itemize}
In the next lemma, we describe certain elements of the set $\mathcal{S}^\prime_{\geq 1}$.

\begin{lemma}\label{l:sp_w}
    For any BD pair $(\bg^r,\bg^c)$, the following statements hold:
    \begin{enumerate}[1)]
        \item For any $\alpha \in \Phi_+(A_{n-1})$ and $\beta \in \Phi_+(A_{m-1})$, $-h_{\alpha}^r(h_{\beta}^c)^T \in \mathcal{S}^\prime_{\geq 1}$;
        \item For any $\alpha,\beta \in \Phi_+(A_{n-1})$ such that $\gamma_r^d(\alpha) = \beta$ for some $d \geq 1$, and any $i,j \in [1,n]$, \begin{equation*}
            h_{\alpha}^r f_i^T - h_{\beta}^r f_{j}^T \in \mathcal{S}^\prime_{\geq 1};
        \end{equation*}
        \item For any $\alpha,\beta \in \Phi_+(A_{m-1})$ such that $\gamma_c^d(\alpha) = \beta$ for some $d \geq 1$, and any $i,j \in [1,m]$,
        \begin{equation*}
            -e_{i} (h_{\alpha}^c)^T + e_j (h_{\beta}^c)^T \in\mathcal{S}^\prime_{\geq 1}.
        \end{equation*}
    \end{enumerate}
\end{lemma}
\begin{proof}
    \begin{enumerate}[1)]
        \item If $\alpha \in \Phi_+(A_{n-1})$ and $\beta \in \Phi_+(A_{m-1})$, write 
        \begin{equation*}
            \alpha = \alpha_p + \cdots + \alpha_{p+k-1}, \quad \beta = \alpha_q + \cdots + \alpha_{q+l-1}.
        \end{equation*}
        Then
        \begin{equation*}
            -h_{\alpha}^r(h^c_{\beta})^T = -(h^r_{p} + \cdots + h^r_{p+k-1})(h_{q}^c + \cdots + h_{q+l-1}^c)^T = \sum_{i=p}^{p+k-1} \sum_{j=q}^{q+l-1} \btheta_{ij} \in \mathcal{S}^\prime_{\geq 1}.
        \end{equation*}
        \item If $\alpha, \beta \in \Phi_+(A_{n-1})$ are such that $\gamma_r^d(\alpha) = \beta$ for some $d \geq 1$, and if $i,j \in [1,n]$, observe that
        \begin{equation*}
            (h_{\alpha}^r f_i^T - h_{\beta}^r f_j^T) = (h_{\alpha}^rf_1^T-h_{\beta}^rf_m^T)+ (-h_{\alpha}^r(f_1-f_i)^T) + (-h_{\beta}^r(f_j-f_m)^T).
        \end{equation*}
        Hence, if we show that $h_{\alpha}^r f_1^T - h_{\beta}^r f_m^T \in \mathcal{S}^\prime_{\geq 1}$, then $h_{\alpha}^r f_i^T - h_{\beta}^r f_j^T \in \mathcal{S}^\prime_{\geq 1}$. For $d > 1$, observe that
        \begin{equation*}
            h_{\alpha}^rf_1^T - h_{\gamma_r^2(\alpha)}^r f_m^T = (h_{\alpha}^r f_1^T-h_{\gamma_r(\alpha)}^r f_m^T) +(-h_{\gamma_r(\alpha)}^r(f_1-f_m)^T) + (h_{\gamma_r(\alpha)}^r f_1^T - h_{\gamma_r^2(\alpha)}^rf_m^T).
        \end{equation*}
        Iterating this argument, we see that it suffices to consider the case $d=1$. If $\alpha = \alpha_p + \cdots +\alpha_{p+k-1}$, we see that
        \begin{equation*}\begin{split}
            h_{\alpha}^r f_1^T - h_{\beta}^r f_m^T = (h_{\alpha_{p}}^r f_1^T - h_{\gamma_r(\alpha_p)}^r f_m^T) + (h_{\alpha_{p+1}}^r f_1^T - h_{\gamma_r(\alpha_{p+1})}^r f_m^T) + \cdots\\ \cdots + (h_{\alpha_{p+k-1}}^r f_1^T - h_{\gamma_r(\alpha_{p+k-1})}^r f_m^T) \in \mathcal{S}^\prime_{\geq 1}.
            \end{split}
        \end{equation*}
        Thus the second statement holds.
        \item The proof of this statement is similar to that of the previous statement.
    \end{enumerate}
    Thus all three statements hold.
\end{proof}
\begin{proposition}
    For any pair $(\bbg^r,\bbg^c)$ such that $\pi_0:=(\pi_{(\bg^r,\bg^c)})_0$ is $\mathbb{T}_{(\bg^r,\bg^c)}$-log-symplectic, $\pi_{(\bg^r,\bg^c)}$ is the normalized maximal admissible $\mathbb{T}_{(\bg^r,\bg^c)}$-invariant Poisson deformation of $\pi_0$.
\end{proposition}

\begin{proof}
    Set $\pi := \pi_{(\bg^r,\bg^c)}$, and let $\mathsf{W}(\pi)$ be the set of all the nonzero $(\mathbb{C}^{\times})^{nm}$-weights appearing in the monomial terms of $\pi$. By Proposition~\ref{p:sp_eqls_sp0}, $\mathcal{S}^\prime = \mathcal{S}(\pi_0)$, and by Proposition~\ref{p:BD_W1}, $\mathcal{S}(\pi_0)$ has Property~\eqref{eq:W1}. To establish the statement, it remains to show the following:
    \begin{enumerate}[1)]
        \item $\mathcal{S}(\pi_0)\subseteq \mathsf{W}(\pi) \subseteq \mathcal{S}(\pi_0)_{\geq 1}$ (in particular, $\pi$ is $\mathcal{S}(\pi_0)$-admissible);
        \item $\pi^{\mathcal{S}(\pi_0)} = \sum_{\btheta \in \mathcal{S}(\pi_0)}V_{\btheta}^\prime$ (that is, $\pi$ is normalized);
    \end{enumerate}
    (see equation~\eqref{eq:Vprime} for the definition of $V_{\btheta}^\prime$). Once established, it follows from Theorem~\ref{thm:W1} that the Poisson deformation of $\pi_0$ along $\pi^{\mathcal{S}(\pi_0)}$ equals $\pi$. 
   
    Let us recall that the Poisson bivector $\pi$ can be written as
    \begin{equation}\label{eq:piw}
        \pi = \frac{1}{2}\sum_{i,k=1}^n \sum_{j,l=1}^m  \langle R^c(\nabla^L x_{ij}),\nabla^L x_{kl} \rangle \frac{\partial}{\partial x_{ij}}\wedge \frac{\partial}{\partial x_{kl}} - \frac{1}{2}\sum_{i,k=1}^n \sum_{j,l=1}^m  \langle R^r(\nabla^R x_{ij}),\nabla^R x_{kl} \rangle \frac{\partial}{\partial x_{ij}}\wedge \frac{\partial}{\partial x_{kl}}.
    \end{equation}
    Note that the factor $1/2$ is due to the fact that the sums run over all indices. We also keep $R^c$ and $R^r$ in the formula for convenience, rather than switching to their skew-symmetric parts. As we already mentioned in Subsection~\ref{s:bd_background}, writing $R = \frac{1}{2}\id + S$ (for both the column and the row labels), it's not difficult to see that the symmetric parts cancel each other out.
    
    Expanding the terms with the left gradients, we see that
    \begin{equation}\label{eq:rc_1}
    \begin{split}
        \langle R^c(\nabla^L x_{ij}),\nabla^L x_{kl} \rangle &= \sum_{\mu_1=1}^m\sum_{\mu_2=1}^m x_{i\mu_1}x_{k\mu_2} \langle R^c(e_{j\mu_1}),e_{l\mu_2}\rangle = \\ &=x_{ij}x_{kl}\langle R_0^c(e_{jj}),e_{ll}\rangle + \sum_{\mu_1=1}^m\sum_{\mu_2=1}^mx_{i\mu_1}x_{k\mu_2}\langle \pr_{>}e_{j\mu_1},e_{l\mu_2}\rangle + \\ &+\sum_{\mu_1=1}^m \sum_{\mu_2=1}^m x_{i\mu_1}x_{k\mu_2}\langle \frac{\gamma_c}{1-\gamma_c}\pr_>e_{j\mu_1},e_{l\mu_2}\rangle - \\&-\sum_{\mu_1=1}^m \sum_{\mu_2=1}^m x_{i\mu_1}x_{k\mu_2}\langle \frac{\gamma_c^*}{1-\gamma_c^*}\pr_{<}e_{j\mu_1},e_{l\mu_2}\rangle.
        \end{split}
    \end{equation}
    Let us analyze the last term. Recall that
    \begin{equation*}
         \frac{\gamma_c^*}{1-\gamma_c^*}\pr_{<} = \sum_{d \geq 1} (\gamma_c^*)^d \pr_{<},
    \end{equation*}
    where the sum is finite, since $\gamma_c^*\pr_{<}$ is nilpotent. Hence the last term in~\eqref{eq:rc_1} contributes bivectors of the form
    \begin{equation}\label{eq:rc_4}
        -x_{i\mu_1}x_{k\mu_2}\langle (\gamma_c^*)^d\pr_{<}e_{j\mu_1},e_{l\mu_2}\rangle\frac{\partial}{\partial x_{ij}}\wedge \frac{\partial}{\partial x_{kl}}.
    \end{equation}
    If $j < \mu_1$ or $l > \mu_2$, the above term is zero. Otherwise, let $j > \mu_1$, $l < \mu_2$, and assume there are positive roots $\alpha,\beta \in \Phi_+(A_{m-1})$ such that $e_{l\mu_2} = e_{\alpha}$, $e_{-\beta} = e_{j\mu_1}$, and $(\gamma_c^*)^d \beta = \alpha$. We see that the weight $\bw$ of the above term is given by
    \begin{equation*}\begin{split}
        \bw = e_{i\mu_1} + e_{k\mu_2} - e_{ij} - e_{kl} &= e_i(f_{\mu_1}-f_{j})^T - e_k(f_{l}-f_{\mu_2})^T = \\ &=e_i(h_{\beta}^c)^T - e_k (h_{\alpha}^c)^T.
        \end{split}
    \end{equation*}
    By Lemma \ref{l:sp_w}, we conclude that $\bw  \in \mathcal{S}(\pi_0)_{\geq 1}$.

    In a similar way, one analyzes the other terms. For example, for the third term, if $j <\mu_1$, $l > \mu_2$, $e_{\alpha}:=e_{j\mu_1}$, $e_{-\beta}:= e_{l\mu_2}$, and $\gamma_c^d\alpha = \beta$, then
    \begin{equation}\label{eq:rc_3}
        \mathrm{wt}\Big(x_{i\mu_1}x_{k\mu_2}\langle \gamma_c^d\pr_>e_{j\mu_1},e_{l\mu_2}\rangle \frac{\partial}{\partial x_{ij}}\wedge \frac{\partial}{\partial x_{kl}}\Big) =  -e_i(h_{\alpha}^c)^T + e_k(h_{\beta}^c)^T \in \mathcal{S}(\pi_0)_{\geq 1}.
    \end{equation}
    Here, $\mathrm{wt}$ stands for the $(\mathbb{C}^{\times})^{nm}$-weight of the monomial bivector. The above equation is also valid for $d=0$, which determines the weights contributed by the second sum in~\eqref{eq:rc_1}. The first term in~\eqref{eq:rc_1} contributes a log-canonical bivector.

    The terms in~\eqref{eq:piw} involving the right gradients can be expanded as
    \begin{equation}\label{eq:rr_1}
    \begin{split}
        \langle R^r(\nabla^R x_{ij}),\nabla^R x_{kl} \rangle &= \sum_{\mu_1=1}^n\sum_{\mu_2=1}^n x_{\mu_1 j} x_{\mu_2 l} \langle R^r(e_{\mu_1 i}),e_{\mu_2 k}\rangle = \\ &=x_{ij}x_{kl}\langle R_0^r(e_{ii}),e_{kk}\rangle +  \sum_{\mu_1=1}^n\sum_{\mu_2=1}^n x_{\mu_1 j} x_{\mu_2 l}\langle \pr_{>} e_{\mu_1 i}, e_{\mu_2 k}\rangle + \\ &+ \sum_{\mu_1=1}^n\sum_{\mu_2=1}^n x_{\mu_1 j} x_{\mu_2 l}\langle \frac{\gamma_r}{1-\gamma_r}\pr_{>}e_{\mu_1 i},e_{\mu_2 k}\rangle  - \\ &-\sum_{\mu_1=1}^n\sum_{\mu_2=1}^n x_{\mu_1 j} x_{\mu_2 l}\langle \frac{\gamma_r^*}{1-\gamma_r^*}\pr_{<}e_{\mu_1 i},e_{\mu_2 k}\rangle.
        \end{split}
    \end{equation}
    For the last term, we see that for any $d \geq 1$, if $\mu_1 > i$, $\mu_2 < k$, and $\alpha,\beta \in \Phi_+(A_{n-1})$ are such that $e_{\alpha} = e_{\mu_2 k}$, $e_{-\beta} = e_{\mu_1 i}$, and $(\gamma_r^*)^d\beta = \alpha$, then
    \begin{equation}\label{eq:rr_4}
        \mathrm{wt} \Big( x_{\mu_1 j} x_{\mu_2 l} \langle (\gamma_r^*)^d \pr_{<} e_{\mu_1 i},e_{\mu_2 k}\rangle \frac{\partial}{\partial x_{ij}}\wedge \frac{\partial}{\partial x_{kl}} \Big) = h_{\alpha}^r f_l^T - h_{\beta}^r f_j^T \in \mathcal{S}(\pi_0)_{\geq 1}.
    \end{equation}
    For the third sum in~\eqref{eq:rr_1}, if $\mu_1 < i$, $\mu_2 > k$, and $\alpha,\beta \in \Phi_+(A_{n-1})$ are such that $e_{\alpha} = e_{\mu_1 i}$, $e_{-\beta} = e_{\mu_2 k}$, and $\gamma_r^d\alpha = \beta$, then
    \begin{equation}\label{eq:rr_3}
        \mathrm{wt} \Big(x_{\mu_1 j} x_{\mu_2 l}\langle \gamma_r^d\pr_> e_{\mu_1i},e_{\mu_2 k}\rangle \frac{\partial}{\partial x_{ij}}\wedge \frac{\partial}{\partial x_{kl}} \Big) = h_{\alpha}^r f_j^T - h_{\beta}^r f_l^T \in \mathcal{S}(\pi_0)_{\geq 1}.
    \end{equation}
    The above equation is also valid for $d=0$, which determines the weights contributed by the second sum in~\eqref{eq:rr_1}.
    The first term in~\eqref{eq:rr_1} is log-canonical.

    The above computations show that $\mathsf{W}(\pi)\subseteq \mathcal{S}(\pi_0)_{\geq 1}$. Now we show that for each $\btheta \in \mathcal{S}(\pi_0)$, the corresponding monomial bivector in $\pi$ equals $V_{\btheta}^\prime$.
    \begin{description}
        \item[Terms of weights $\btheta_{pq}$.] Let us consider the contribution of the second terms in \eqref{eq:rc_1} and \eqref{eq:rr_1} to \eqref{eq:piw}; that is,
        \begin{equation*}
            \frac{1}{2}\sum_{i,k=1}^n\sum_{j,l=1}^m [j < l] x_{il}x_{kj} \frac{\partial}{\partial x_{ij}}\wedge \frac{\partial}{\partial x_{kl}} - \frac{1}{2}\sum_{i,k=1}^n\sum_{j,l=1}^m [k<i]x_{kj}x_{il}\frac{\partial}{\partial x_{ij}}\wedge \frac{\partial}{\partial x_{kl}}
        \end{equation*}
        where $[j<l]$ and $[k<i]$ denote predicates which are equal to $1$ or $0$ depending on whether the corresponding condition is true or false. Given the weight $\btheta_{pq}$, we can locate a monomial in the first sum with the $(\mathbb{C}^{\times})^{nm}$-weight $\btheta_{pq}$ by setting $(i,j,k,l) = (p,q,p+1,q+1)$. In the second sum, we locate the monomial with weight $\btheta_{pq}$ by setting $(i,j,k,l) = (p+1,q+1,p,q)$. Combining the two monomials, we obtain 
        \begin{equation*}\begin{split}
            \frac{1}{2}x_{p,q+1}x_{p+1,q}\frac{\partial}{\partial x_{pq}} \wedge \frac{\partial}{\partial x_{p+1,q+1}} - \frac{1}{2}x_{p,q+1}x_{p+1,q}\frac{\partial}{\partial x_{p+1,q+1}}\wedge \frac{\partial}{\partial x_{pq}} = \\= x_{p,q+1} x_{p+1,q}\frac{\partial}{\partial x_{pq}}\wedge \frac{\partial}{\partial x_{p+1,q+1}} = V_{\btheta_{pq}}^\prime.
            \end{split}
        \end{equation*}
        In particular, we see that $\pi$ is normalized at $V_{\btheta_{pq}}^\prime$, for the order in the wedge product together with the coefficient agrees with the right-hand side of the LM equation:
        \begin{equation*}
            \sum_{k=1}^n\sum_{l=1}^m\blambda_{-,kl} (\btheta_{pq})_{kl} = e_{pq}-e_{p+1,q+1}.
        \end{equation*}
        \item[Terms of weights $\btheta_p^c$.] These weights come from the contribution of the third and the fourth terms in~\eqref{eq:rc_1} to~\eqref{eq:piw}. Specifically, by setting $\mu_1 = j-1$, $\mu_2 = l + 1$, and $d=1$ in~\eqref{eq:rc_4}, we see that the fourth term in~\eqref{eq:rc_1} contributes to~\eqref{eq:piw} terms of the form
        \begin{equation}\label{eq:rc_tpc1}
            -\frac{1}{2}x_{i,j-1} x_{k,l+1} \frac{\partial}{\partial x_{ij}}\wedge \frac{\partial}{\partial x_{kl}}, \qquad \gamma_c^*(\alpha_{j-1}) = \alpha_l.
        \end{equation}
        From equation~\eqref{eq:rc_3}, by setting $\mu_1 = j+1$, $\mu_2 = l-1$, and $d = 1$, we find that the third term in~\eqref{eq:rc_1} contributes to~\eqref{eq:piw} terms of the form
        \begin{equation}\label{eq:rc_tpc2}
            \frac{1}{2}x_{i,j+1}x_{k,l-1} \frac{\partial}{\partial x_{ij}}\wedge \frac{\partial}{\partial x_{kl}}, \qquad \gamma_c(\alpha_j) = \alpha_{l-1}. 
        \end{equation}
        Now, suppose that $\gamma_c(\alpha_p) = \alpha_q$ for some $\alpha_p \in \Gamma_1^c$. Setting $(i,j,k,l) = (1,q+1,n,p)$ in~\eqref{eq:rc_tpc1} and $(i,j,k,l) = (n,p,1,q+1)$ in~\eqref{eq:rc_tpc2}, we see that the two terms combine into a monomial bivector of the form
        \begin{equation*}
            x_{n,p+1} x_{1,q} \frac{\partial}{\partial x_{n,p}}\wedge \frac{\partial}{\partial x_{1,q+1}}.
        \end{equation*}
        The coefficient and the order in the wedge product agree with the right-hand side of the LM equation
        \begin{equation*}                 \sum_{k=1}^n\sum_{l=1}^m\blambda_{-,kl}( \btheta_{p}^c)_{kl} = e_{np}-e_{1,q+1}.
        \end{equation*}
        Therefore, $\pi$ is normalized at $V_{\btheta_{p}^c}^\prime$.
        \item[Terms of weights $\btheta_p^r$.] As in the previous case, these weights come from the contributions of the third and the fourth terms in~\eqref{eq:rr_1} to~\eqref{eq:piw}. Specifically, setting $\mu_1 = i+1$, $\mu_2 = k-1$, and $d = 1$ in~\eqref{eq:rr_4}, we obtain the term
        \begin{equation}\label{eq:rr_tpc1}
            \frac{1}{2}x_{i+1,j}x_{k-1,l} \frac{\partial}{\partial x_{ij}}\wedge\frac{\partial}{\partial x_{kl}}, \qquad \gamma_r^*(\alpha_i) = \alpha_{k-1}.
        \end{equation}
        By setting $\mu_1 = i-1$ and $\mu_2 = k+1$ in~\eqref{eq:rr_3}, we obtain the term  
        \begin{equation}\label{eq:rr_tpc2}
            -\frac{1}{2}x_{i-1,j}x_{k+1,l} \frac{\partial}{\partial x_{ij}}\wedge \frac{\partial}{\partial x_{kl}}, \qquad \gamma_r(\alpha_{i-1}) = \alpha_k.
        \end{equation}
        Now, to obtain the term with the $(\mathbb{C}^{\times})^{nm}$-weight equal to $\btheta_p^r$, we set $(i,j,k,l) = (q,m,p+1,1)$ in~\eqref{eq:rr_tpc1} and $(i,j,k,l) = (p+1,1,q,m)$ in~\eqref{eq:rr_tpc2}. Their sum yields
        \begin{equation*}
            x_{q+1,m}x_{p,1} \frac{\partial}{\partial x_{qm}}\wedge \frac{\partial}{\partial x_{p+1,1}},
        \end{equation*}
        which agrees with the right-hand side of the LM equation
        \begin{equation*}
            \sum_{k=1}^n\sum_{l=1}^m\blambda_{-,kl}( \btheta_{p}^r)_{kl} = e_{qm} - e_{p+1,1} , \quad p \in \Gamma_1^r.
        \end{equation*}
    \end{description}
    Thus, $\pi$ is the normalized maximal admissible $\mathbb{T}_{(\bg^r,\bg^c)}$-invariant Poisson deformation of $\pi_0$.
\end{proof}

\subsection{Combinatorial models for verifying aperiodicity}\label{s:meander}
In this subsection, we develop a generalization of the aperiodicity condition of BD pairs $(\bg^r,\bg^c)$ from reductive algebraic groups (and in particular, from square matrices) to rectangular matrices. The key to the generalization is Corollary~\ref{c:cycles_tlogs}, which states that $\pi_0 := (\pi_{(\bbg^r,\bbg^c)})_0$ is $\mathbb{T}_{(\bg^r,\bg^c)}$-log-symplectic if and only if the smoothing graph $\Gamma(\mathcal{S}^\prime)$ is acyclic. 

The notion of aperiodic pairs of BD triples $(\bg^r,\bg^c)$ arising from the same root system was introduced in~\cite{gsv-plethora}. As we mentioned in the introduction, it is expected that a simple complex algebraic group $G$ has an open $\mathbb{T}_{(\bg^r,\bg^c)}$-leaf if and only if $(\bg^r,\bg^c)$ is aperiodic. Moreover, if $G$ is simply connected, $\mathbb{C}[G]$ is expected to carry the structure of a (non-generalized) upper cluster algebra compatible with $\pi_{(\bbg^r,\bbg^c)}$. If $(\bg^r,\bg^c)$ is not aperiodic, then $\mathbb{C}[G]$ is expected to carry the structure of a generalized upper cluster algebra compatible with $\pi_{(\bbg^r,\bbg^c)}$.

The original definition of aperiodicity is as follows. 

\begin{definition}\label{d:aperiodic}
    A pair of BD triples $(\bg^r,\bg^c)$ arising from the same root system is called \emph{aperiodic} if the map $\upsilon:=\gamma_r^*(-w_0)\gamma_c(-w_0)$ is nilpotent as a map on simple roots.
\end{definition}
Here, the map $\upsilon$ is defined on its largest possible domain $\dom \upsilon$, and the nilpotency condition means that for every $\alpha \in \dom \upsilon$, there exists $k > 0$ such that $\upsilon^k(\alpha) \notin \dom \upsilon$.

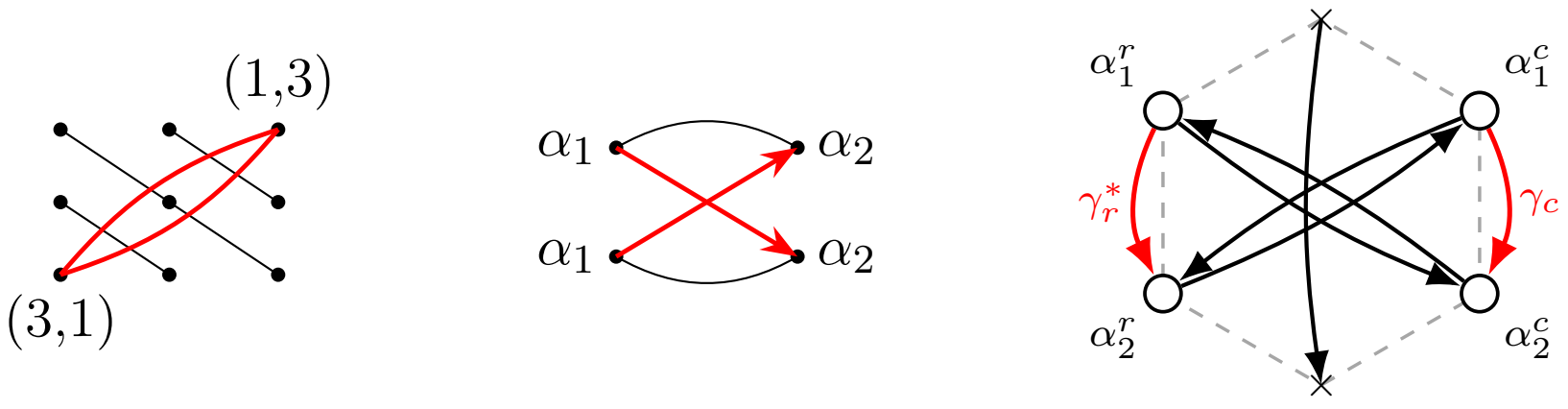
\begin{figure}[t]
\centering
\begin{tikzpicture}[baseline=(current bounding box.center)]
    \filldraw[black] (0,1) circle (\bdvsize);
    \filldraw[black] (0.75,1) circle (\bdvsize);
    \filldraw[black] (1.5,1) circle (\bdvsize)
    node[anchor=south]{(1,3)};

    \filldraw[black] (0,0.5) circle (\bdvsize);
    \filldraw[black] (0.75,0.5) circle (\bdvsize);
    \filldraw[black] (1.5,0.5) circle (\bdvsize);

    \filldraw[black] (0,0) circle (\bdvsize) node[anchor=north]{(3,1)};
    \filldraw[black] (0.75,0) circle (\bdvsize);
    \filldraw[black] (1.5,0) circle (\bdvsize);

    \draw (0,0.5) -- (0.75,0);
    \draw (0,1)   -- (1.5,0);
    \draw (0.75,1)-- (1.5,0.5);

    \draw[red, line width=0.9pt]
        (0,0) to[bend left=16] (1.5,1);

    \draw[red, line width=0.9pt]
        (0,0) to[bend right=16] (1.5,1);
        
\end{tikzpicture}
    \hspace{25pt}
    \begin{tikzpicture}[baseline=(current bounding box.center)]
    
    \filldraw[black] (0,0) circle (\bdvsize) node[anchor=east]{$\alpha_1$};
    \filldraw[black] (1.25,0) circle (\bdvsize) node[anchor=west]{$\alpha_2$};
    \filldraw[black] (0,0.75) circle (\bdvsize) node[anchor=east]{$\alpha_1$};
    \filldraw[black] (1.25,0.75) circle (\bdvsize) node[anchor=west]{$\alpha_2$};

    \draw (0,0) to[bend right=30] (1.25,0);
    \draw (0,0.75) to[bend left=30] (1.25,0.75);

    \draw[-{Stealth},red, line width = 0.9pt] (0,0.75) to (1.25,0);
    \draw[-{Stealth},red, line width = 0.9pt] (0,0) to (1.25, 0.75) ;
        
\end{tikzpicture}
\hspace{25pt}
\begin{tikzpicture}[
    baseline=(current bounding box.center),
    >=Latex,
    scale=0.72,
    every node/.style={font=\scriptsize},
    vertex/.style={
        circle,
        draw,
        fill=white,
        inner sep=0pt,
        minimum size=2.5mm,
        line width=0.7pt
    },
    bdarrow/.style={->, line width=0.8pt},
    redarrow/.style={->,red, line width=0.95pt},
    affedge/.style={dashed, gray!70, line width=0.65pt}
]
\def\R{1.75}

\foreach \i/\ang in {
  0/90,
  1/150,
  2/210,
  3/270,
  4/330,
  5/30
}{
  \coordinate (p\i) at (\ang:\R);
}

\draw[affedge] (p0)  -- (p1);
\draw[affedge] (p1)  -- (p2);
\draw[affedge] (p2)  -- (p3);
\draw[affedge] (p3)  -- (p4);
\draw[affedge] (p4)  -- (p5);
\draw[affedge] (p5)  -- (p0);

\node at (p0) {$\times$};
\node at (p3) {$\times$};
\node[vertex] (v1) at (p1) {};
\node[vertex] (v2) at (p2) {};
\node[vertex] (v4) at (p4) {};
\node[vertex] (v5) at (p5) {};
\node[above left=1pt] at (v1) {$\alpha_1^r$};
\node[below left=1pt]       at (v2) {$\alpha_2^r$};
\node[below right=1pt]       at (v4) {$\alpha_2^c$};
\node[above right=1pt]       at (v5) {$\alpha_1^c$};

\draw[bdarrow] (p0)  to[bend right=8] (p3);

\draw[bdarrow] (v1)  to[bend right=8] (v4);
\draw[bdarrow] (v2)  to[bend right=8] (v5);
\draw[bdarrow] (v5)  to[bend right=8] (v2);
\draw[bdarrow] (v4)  to[bend right=8] (v1);

\draw[redarrow] (v1) to[bend right=26] (v2);
\draw[redarrow] (v5) to[bend left=26] (v4);

\node[red] at (-2.1,0) {$\gamma_r^*$};
\node[red] at (2.1,0) {$\gamma_c$};

\end{tikzpicture}

\caption{Smoothing graph (left), BD graph (middle), affine BD meander (right) for
$n=m=3$, with $\gamma_r(\alpha_2^r) = \alpha_1^r$ and $\gamma_c(\alpha_1^c) = \alpha_2^c$.}
\label{fig:3graphs}
\end{figure}

The aperiodicity condition can be illustrated with the so-called \emph{BD graph} $G_{(\bg^r,\bg^c)}$. For a root system of rank $r$, its vertex set consists of two copies of the interval $[1,r]$. In each copy, one draws an undirected edge between $i$ and $i^*$ (where $i^* \in [1,r]$ is such that $\alpha_{i^*} = (-w_0) \alpha_i$ for a simple root~$\alpha_i$). Such edges are called \emph{horizontal}, and they depict the action of $(-w_0)$ in $\upsilon$. One then draws oriented edges between the two copies of $[1,r]$. Specifically, there is a directed edge from $i$ in the upper half to $j$ in the lower half if $\alpha_j = \gamma_c(\alpha_i)$. Likewise, there is a directed edge from $j$ in the lower half to $i$ in the upper half if $\alpha_i = \gamma_r^*(\alpha_j)$. 
The edges that are not horizontal are called \emph{inclined}. The map $\upsilon$ is nilpotent if and only if $G_{(\bg^r,\bg^c)}$ has no alternating cycles (that is, cycles in which horizontal and inclined edges alternate). For examples of BD graphs drawn using a different convention,
see~\cite{gsv-plethora,multdual}. A simple example is provided in Figure~\ref{fig:3graphs}, where we also compare the BD graph with the corresponding smoothing graph and the BD meander (defined below).

Now assume that $\bg^r$ and $\bg^c$ are given with respect to the simple roots of types $A_{n-1}$ and $A_{m-1}$, respectively. We embed $A_{n-1}$ and $A_{m-1}$ into the affine Dynkin diagram $\hat{A}_{n+m-1}$ via
\begin{align*}
    &\imath_r : A_{n-1} \rightarrow \hat{A}_{n+m-1}, \quad \imath_r(\alpha_i) = \alpha_{n-i}, \ i \in [1,n-1];\\
    &\imath_c : A_{m-1} \rightarrow \hat{A}_{n+m-1}, \quad \imath_c(\alpha_i) = \alpha_{n+i}, \ i \in [1,m-1].
\end{align*}
Define $\alpha_i^r := \imath_r(\alpha_i)$, $i \in [1,n-1]$, and $\alpha_i^c :=\imath_c(\alpha_i)$, $i \in [1,m-1]$. By means of the two embeddings, we realize a pair $(\bg^r,\bg^c)$ of BD triples as a single BD triple $\bg$ in $A_{n+m-1}$ (not affine). More precisely, we define $\bg:=(\Gamma_1,\Gamma_2,\gamma)$, where 
\begin{equation*}
    \Gamma_1 := \imath_r(\Gamma_2^r)\cup \imath_c(\Gamma_1^c), \quad \Gamma_2 := \imath_r(\Gamma_1^r)\cup \imath_c(\Gamma_2^c),
\end{equation*}
and the map $\gamma$ is uniquely determined by the requirement that the following diagram commutes:
\begin{equation}\label{eq:defg}
    \xymatrix{
    \Gamma_1^r \ar[r]^{\imath_r} & \Gamma_2 & \Gamma_2^c\ar[l]_{\imath_c}\\
    \Gamma_2^r\ar[u]^{\gamma^*_r} \ar[r]^{\imath_r} & \Gamma_1 \ar[u]^{\gamma} & \Gamma_1^c \ar[l]_{\imath_c} \ar[u]^{\gamma_c}
    }
\end{equation}
Let $\tau_n$ be the cyclic shift by $n$ on $\hat{A}_{n+m-1}$; that is, $\tau_n(\alpha_i) = \alpha_{i+n}$, where the indices are taken modulo $n+m$. Note that $\tau_n$ is an isometry.
\begin{definition}\label{d:aperiodic_mat}
The pair $(\bg^r,\bg^c)$ is called \emph{aperiodic} if the map $\gamma\tau_n$ is nilpotent. 
\end{definition}
As before, the map $\gamma\tau_n$ is viewed on its largest possible domain. 

We will need the following lemma to relate the nilpotency of $\gamma\tau_n$ to Definition~\ref{d:aperiodic}, as well as to $\mathbb{T}_{(\bg^r,\bg^c)}$-log-symplecticity.

\begin{lemma}\label{l:taun_w0}
For any $p \in [1,n-1]$ and $q \in [1,m-1]$,
\begin{equation}\label{eq:taun_w0}
    \tau_n(\alpha_p^r) = 
    \begin{cases}
        \alpha_{m+p}^r &\text{if} \  p < n-m\\
        \alpha_0 &\text{if} \ p = n-m;\\
        \alpha_{n-p}^c &\text{if} \ p > n - m.
    \end{cases} 
    \quad 
    \tau_n(\alpha_n) = 
    \begin{cases}
        \alpha_n^c &\text{if} \  n < m \\
        \alpha_0 &\text{if} \  n = m \\
        \alpha_m^r &\text{if} \  n > m
    \end{cases}
    \quad
    \tau_n(\alpha_q^c) = 
    \begin{cases}
    \alpha_{n+q}^c &\text{if} \ q <m-n;\\
    \alpha_0 &\text{if} \ q = m-n;\\
    \alpha_{m-q}^r &\text{if} \ q > m-n.
    \end{cases}
\end{equation}

\end{lemma}
\begin{proof}
Recall that $\alpha_p^r = \imath_r(\alpha_p)= \alpha_{n-p} $ and $\alpha_q^c = \imath_c(\alpha_q) = \alpha_{n+q}$. For $\tau_n(\alpha_p^r)$, we see that if $n-p \in [1,m-1]$, then
\begin{equation*}
    \tau_n(\alpha_p^r) = \alpha_{n+(n-p)} = \alpha_{n-p}^c;
\end{equation*}
if $p \in [1,n-1]$ is such that $n-p > m$, then
\begin{equation*}
    \tau_n(\alpha_p^r) = \alpha_{n+(n-p)} = \alpha_{n-p-m} = \alpha_{m+p}^r,
\end{equation*}
where we used $n+(n-p) \equiv n-(m+p) \pmod{n+m}$. If $n-p = m$, then
\begin{equation*}
    \tau_n(\alpha_p^r) = \tau_n(\alpha_{n-p}) = \alpha_{n-p+n} = \alpha_0.
\end{equation*}
The equations for $\tau_n(\alpha_n)$ and $\tau_n(\alpha_q^c)$ follow in a similar way.
\end{proof}

The next proposition shows that if $n=m$, then Definition~\ref{d:aperiodic_mat} agrees with Definition~\ref{d:aperiodic}.

\begin{proposition}
    If $n = m$, the following identities hold:
    \begin{align*}
        \imath_r^{-1}(\gamma \tau_n)^2\imath_r &=  \gamma_r^*(-w_0)\gamma_c(-w_0), \\
        \imath_c^{-1}(\gamma \tau_n)^2\imath_c &=  \gamma_c(-w_0)\gamma^*_r(-w_0).
    \end{align*}
    In particular, $\gamma\tau_n$ is nilpotent if and only if $\upsilon$ is nilpotent. 
\end{proposition}

\begin{proof}
Since $n=m$, equations~\eqref{eq:taun_w0} can be recast as 
\begin{equation*}
    \tau_n \imath_r = \imath_c(-w_0), \quad \tau_n \imath_c = \imath_r(-w_0).
\end{equation*}
By the commutative diagram \eqref{eq:defg}, on the appropriate domains, we have
\begin{equation*}\begin{split}
    \imath_r^{-1}\gamma \tau_n \gamma \tau_n \imath_r = \imath_r^{-1}\gamma \tau_n \gamma \imath_c(-w_0) = \imath_r^{-1}\gamma \tau_n \imath_c \gamma_c(-w_0)  \\ = \imath_r^{-1}\gamma \imath_r(-w_0) \gamma_c(-w_0) = \gamma_r^*(-w_0)\gamma_c(-w_0).
    \end{split}
\end{equation*}
This proves the first equation. The second equation follows in a similar way. Now note that the compositions on the right-hand sides of
the equations in the statement are either both nilpotent or both non-nilpotent.
\end{proof}

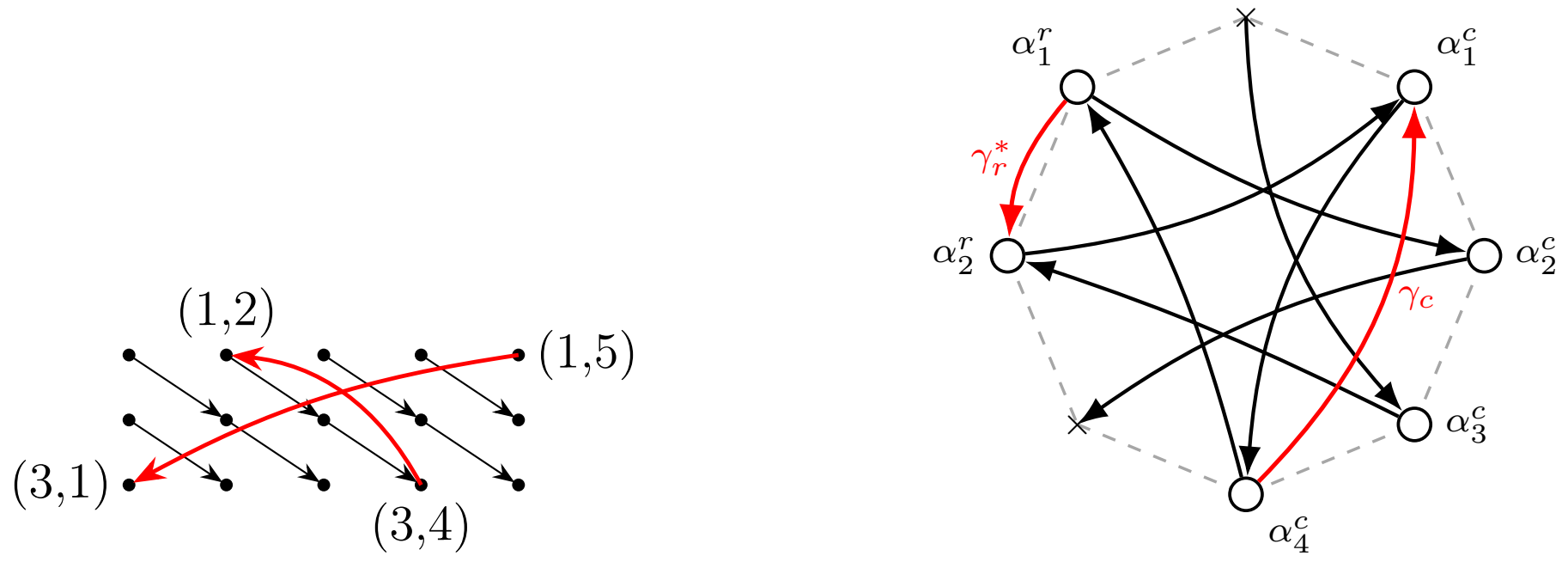
\begin{figure}[t]
\centering
\begin{tikzpicture}
    \filldraw[black] (0,1) circle (\bdvsize);
    \filldraw[black] (0.75,1) circle (\bdvsize) node[anchor=south]{(1,2)};
    \filldraw[black] (1.5,1) circle (\bdvsize);
    \filldraw[black] (2.25,1) circle (\bdvsize);
    \filldraw[black] (3,1) circle (\bdvsize) node[anchor=west]{(1,5)};

    \filldraw[black] (0,0.5) circle (\bdvsize);
    \filldraw[black] (0.75,0.5) circle (\bdvsize);
    \filldraw[black] (1.5,0.5) circle (\bdvsize);
    \filldraw[black] (2.25,0.5) circle (\bdvsize);
    \filldraw[black] (3,0.5) circle (\bdvsize);

    \filldraw[black] (0,0) circle (\bdvsize) node[anchor=east]{(3,1)};
    \filldraw[black] (0.75,0) circle (\bdvsize);
    \filldraw[black] (1.5,0) circle (\bdvsize);
    \filldraw[black] (2.25,0) circle (\bdvsize) node[anchor=north]{(3,4)};
    \filldraw[black] (3,0) circle (\bdvsize);

    \draw[-{Stealth}, shorten >=1pt] (0,0.5) -- (0.75,0);
    \draw[-{Stealth}, shorten >=1pt] (0,1) -- (0.75,0.5);
    \draw[-{Stealth}, shorten >=1pt] (0.75,0.5) -- (1.5,0);
    \draw[-{Stealth}, shorten >=1pt] (0.75,1)-- (1.5,0.5);
    \draw[-{Stealth}, shorten >=1pt] (1.5,0.5)-- (2.25,0);
    \draw[-{Stealth}, shorten >=1pt] (1.5,1) -- (2.25,0.5);
    \draw[-{Stealth}, shorten >=1pt] (2.25,0.5) -- (3,0);
    \draw[-{Stealth}, shorten >=1pt] (2.25,1)-- (3,0.5);

    \draw[-{Stealth}, shorten >=1pt, red, line width=0.9pt]
        (3,1) to[bend right=10] (0,0);

    \draw[-{Stealth}, shorten >=1pt, red, line width=0.9pt]
        (2.25,0) to[bend right=28] (0.75,1);
\end{tikzpicture}
    \hspace{50pt}
\begin{tikzpicture}[
    >=Latex,
    scale=0.72,
    every node/.style={font=\scriptsize},
    vertex/.style={
        circle,
        draw,
        fill=white,
        inner sep=0pt,
        minimum size=2.5mm,
        line width=0.7pt
    },
    bdarrow/.style={->, line width=0.8pt},
    redarrow/.style={->,  draw=red, line width=0.95pt},
    affedge/.style={dashed, gray!70, line width=0.65pt}
]
\def\R{2.55}

\foreach \i/\ang in {
  0/90,
  1/135,
  2/180,
  3/225,
  4/270,
  5/315,
  6/0,
  7/45
}{
  \coordinate (p\i) at (\ang:\R);
}

\foreach \i/\j in {0/1,1/2,2/3,3/4,4/5,5/6,6/7,7/0}{
  \draw[affedge] (p\i) -- (p\j);
}

\node at (p0) {$\times$};
\node[vertex] (v1) at (p1) {};
\node[vertex] (v2) at (p2) {};
\node at (p3) {$\times$};
\node[vertex] (v4) at (p4) {};
\node[vertex] (v5) at (p5) {};
\node[vertex] (v6) at (p6) {};
\node[vertex] (v7) at (p7) {};

\node[above left=1pt]  at (v1) {$\alpha_1^r$};
\node[left=3pt]        at (v2) {$\alpha_2^r$};
\node[below right=1pt] at (v4) {$\alpha_4^c$};
\node[right=3pt]       at (v5) {$\alpha_3^c$};
\node[right=3pt]       at (v6) {$\alpha_2^c$};
\node[above right=1pt] at (v7) {$\alpha_1^c$};

\draw[bdarrow] (v4) to[bend right=8]  (v1);
\draw[bdarrow] (v5) to[bend right=4]  (v2);
\draw[bdarrow] (v7) to[bend right=18] (v4);
\draw[bdarrow] (p0) to[bend right=20] (v5);
\draw[bdarrow] (v1) to[bend right=10] (v6);
\draw[bdarrow] (v2) to[bend right=16] (v7);
\draw[bdarrow] (v6) to[bend right=12] (p3);

\draw[redarrow]
  (v1) to[bend right=18]
  node[pos=.52,xshift=-7pt,yshift=2pt, red]
  {$\gamma_r^*$}
  (v2);

\draw[redarrow]
  (v4) to[bend right=22]
  node[pos=.52,xshift=8pt,yshift=1pt, red]
  {$\gamma_c$}
  (v7);

\%node[font=\scriptsize] at (2.15,3.0)
\end{tikzpicture}
\caption{Oriented smoothing graph $\Gamma(\mathcal{S}')$ (left) and affine BD meander (right) on $\hat{A}_{n+m-1}$ for
$n=3$, $m=5$, $\bg^r = (\{2\},\{1\},2\mapsto 1)$, $\bg^c = (\{4\},\{1\},4\mapsto 1)$.}
\label{fig:affine-bd-meander}
\end{figure}

As the nilpotency of $\upsilon$ can be tested with a BD graph, the nilpotency of $\gamma\tau_n$ can be tested by means of an affine BD meander:

\begin{definition}
    Given a pair $(\bg^r,\bg^c)$ of BD triples in types $A_{n-1}$ and $A_{m-1}$, respectively, an \emph{affine Belavin--Drinfeld meander} (an \emph{affine BD meander}) is an oriented graph $M_{(\bg^r,\bg^c)}$ obtained as follows. The vertex set of $M_{(\bg^r,\bg^c)}$ is identified with the nodes of $\hat{A}_{n+m-1}$. The edge set consists of black and red arrows, which are drawn according to the following steps:
    \begin{itemize}
        \item For each $\alpha_i \in A_{n+m-1}$, draw a black arrow from $\alpha_i$ to $\tau_n(\alpha_i)$;
        \item For each $\alpha_i \in \Gamma_1$, draw a red arrow from $\alpha_i$ to $\gamma(\alpha_i)$.
    \end{itemize}
\end{definition}
 An affine BD meander is illustrated in Figure~\ref{fig:affine-bd-meander}. We cross out the nodes $\alpha_0$ and $\alpha_n$, and we do not draw the black arrow $\alpha_0 \rightarrow \alpha_n$. These simple roots are not in the domain of $\gamma$ and do not affect the nilpotency. 
 
 The way the red arrows are drawn can be unpacked using the commutative diagram~\eqref{eq:defg}:
 \begin{itemize}
     \item There is a red arrow $\alpha_q^r \rightarrow \alpha_p^r$ if and only if $\gamma_r^*(\alpha_q^r) = \alpha_p^r$;
     \item There is a red arrow $\alpha_p^c \rightarrow \alpha_q^c$ if and only if $\gamma_c(\alpha_p^c) = \alpha_q^c$.
 \end{itemize}

\begin{definition}
    Given an affine BD meander $M_{(\bg^r,\bg^c)}$, an \emph{alternating path} (resp. \emph{cycle}) is a directed path (resp. directed cycle) in which red and black arrows alternate. An alternating path is called \emph{maximal} if it cannot be extended to a longer alternating path.
\end{definition}

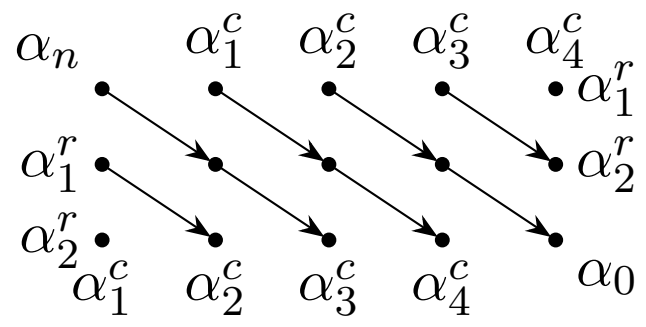
\begin{figure}[t]
\centering
\begin{tikzpicture}
    \filldraw[black] (0,1) circle (\bdvsize) node[anchor=south east]{$\alpha_n$};
    \filldraw[black] (0.75,1) circle (\bdvsize) node[anchor=south]{$\alpha_1^c$};
    \filldraw[black] (1.5,1) circle (\bdvsize) node[anchor=south]{$\alpha_2^c$};
    \filldraw[black] (2.25,1) circle (\bdvsize) node[anchor=south]{$\alpha_3^c$};
    \filldraw[black] (3,1) circle (\bdvsize) node[anchor=south]{$\alpha_4^c$} node[anchor=west]{$\alpha_1^r$};

    \filldraw[black] (0,0.5) circle (\bdvsize) node[anchor=east]{$\alpha_1^r$};
    \filldraw[black] (0.75,0.5) circle (\bdvsize);
    \filldraw[black] (1.5,0.5) circle (\bdvsize);
    \filldraw[black] (2.25,0.5) circle (\bdvsize);
    \filldraw[black] (3,0.5) circle (\bdvsize)  node[anchor=west]{$\alpha_2^r$};

    \filldraw[black] (0,0) circle (\bdvsize) node[anchor=east]{$\alpha_2^r$} node[anchor=north]{$\alpha_1^c$};
    \filldraw[black] (0.75,0) circle (\bdvsize) node[anchor=north]{$\alpha_2^c$};
    \filldraw[black] (1.5,0) circle (\bdvsize) node[anchor=north]{$\alpha_3^c$};
    \filldraw[black] (2.25,0) circle (\bdvsize) node[anchor=north]{$\alpha_4^c$};
    \filldraw[black] (3,0) circle (\bdvsize)node[anchor=north west]{$\alpha_0$};

    \draw[-{Stealth}, shorten >=1pt] (0,0.5) -- (0.75,0);
    \draw[-{Stealth}, shorten >=1pt] (0,1) -- (0.75,0.5);
    \draw[-{Stealth}, shorten >=1pt] (0.75,0.5) -- (1.5,0);
    \draw[-{Stealth}, shorten >=1pt] (0.75,1)-- (1.5,0.5);
    \draw[-{Stealth}, shorten >=1pt] (1.5,0.5)-- (2.25,0);
    \draw[-{Stealth}, shorten >=1pt] (1.5,1) -- (2.25,0.5);
    \draw[-{Stealth}, shorten >=1pt] (2.25,0.5) -- (3,0);
    \draw[-{Stealth}, shorten >=1pt] (2.25,1)-- (3,0.5);

    \end{tikzpicture}
    \caption{Assignment of roots of $\hat{A}_{n+m-1}$ to the boundary of the oriented smoothing graph $\Gamma(\mathcal{S}_{\std})$ for $n=3$, $m=5$. Illustration to the proof of Proposition~\ref{p:aperiodic_tlogsympl}.}
    \label{f:pf_illustration}
\end{figure}

\begin{proposition}
    For any $(\bg^r,\bg^c)$, the map $\gamma \tau_n$ is nilpotent if and only if the affine BD meander $M_{(\bg^r,\bg^c)}$ has no alternating cycles.
\end{proposition}
\begin{proof}
    By definition, a black (resp. red) arrow corresponds to the action of $\tau_n$ (resp. $\gamma$).
\end{proof}
\begin{proposition}\label{p:aperiodic_tlogsympl}
    For any $(\bbg^r,\bbg^c)$, the bivector $\pi_0 := (\pi_{(\bbg^r,\bbg^c)})_0$ is $\mathbb{T}_{(\bg^r,\bg^c)}$-log-symplectic if and only if $(\bg^r,\bg^c)$ is aperiodic.
\end{proposition}
\begin{proof}
    By Corollary~\ref{c:cycles_tlogs}, $\pi_0$ is $\mathbb{T}_{(\bg^r,\bg^c)}$-log-symplectic if and only if the smoothing graph $\Gamma(\mathcal{S}^\prime)$ has no cycles. The idea is to establish a correspondence between the cycles in $\Gamma(\mathcal{S}^\prime)$ and the alternating cycles in $M_{(\bg^r,\bg^c)}$. 

    By construction, the smoothing graph $\Gamma(\mathcal{S}^\prime)$ is an undirected graph, whereas the BD meander $M_{(\bg^r,\bg^c)}$ is oriented. In order to align arrows in $M_{(\bg^r,\bg^c)}$ with consecutive sequences of undirected edges of $\Gamma(\mathcal{S}^\prime)$, we orient each edge $\varepsilon$ of $\Gamma(\mathcal{S}^\prime)$ as follows. Given the corresponding element $\btheta^{\varepsilon} \in \mathcal{S}^\prime$, we let the source $s(\varepsilon)$ and target $t(\varepsilon)$ of $\varepsilon$ be such that $\blambda \btheta^{\varepsilon} = e_{s(\varepsilon)} - e_{t(\varepsilon)}$ (for example, see Figure~\ref{fig:affine-bd-meander}). The explicit values of $\blambda \btheta^{\varepsilon}$ are given in Proposition~\ref{p:bd_lm}.

    Let us consider the subgraph $\Gamma(\mathcal{S}_{\std})$ of the smoothing graph $\Gamma(\mathcal{S}^\prime)$ that consists of the same vertex set and whose edges correspond to $\btheta_{pq}$ for $p\in [1,n-1]$ and $q\in [1,m-1]$ (in other words, $\Gamma(\mathcal{S}_{\std})$ is the smoothing graph for the trivial BD triples). Each connected component of $\Gamma(\mathcal{S}_{\std})$ is a path (possibly of length $0$). Let us assign a simple root from $\hat{A}_{n+m-1}$ to each node of degree $1$ of $\Gamma(\mathcal{S}_{\std})$, and let us assign two simple roots to each node of degree $0$, in the manner illustrated in Figure~\ref{f:pf_illustration}. More precisely,
    \begin{itemize}
        \item For each $p,q \in [1,n-1]$, the nodes $(p+1,1)$ and $(q,m)$ are labeled by $\alpha_p^r$, $\alpha_q^r$, respectively;
        \item For each $p,q \in [1,m-1]$, the nodes $(1,q+1)$ and $(n,p)$ are labeled by $\alpha_q^c$, $\alpha_p^c$, respectively;
        \item The nodes $(1,1)$ and $(n,m)$ are labeled by $\alpha_n$ and $\alpha_0$, respectively.
    \end{itemize}
    Note that the nodes $(n,1)$ and $(1,m)$ receive $2$ labels. We claim there is a path in $\Gamma(\mathcal{S}_{\std})$ with the endpoints labeled by $\alpha_i$ and $\alpha_j$ if and only if there is a black arrow between $\alpha_i$ and $\alpha_j$ in $M_{(\bg^r,\bg^c)}$.
    Indeed, $i-j$ is constant along each path in $\Gamma(\mathcal{S}_{\std})$.
    For instance,
    \begin{itemize}
        \item If $n \leq m$ and $q > m - n$, there is a path $(1,q+1) \rightedge (2,q+2) \rightedge \cdots \rightedge (m-q,m)$. The vertex $(1,q+1)$ is labeled by $\alpha_q^c$, and $(m-q,m)$ is labeled by $\alpha_{m-q}^r$. By Lemma~\ref{l:taun_w0}, $\tau_n(\alpha_q^c) = \alpha_{m-q}^r$, so there is a black arrow $\alpha_q^c \rightarrow \alpha_{m-q}^r$ in $M_{(\bg^r,\bg^c)}$;
        \item If $n \geq m $ and $p < n - m$, there is a path $(p+1,1) \rightedge (p+2,2) \rightedge \cdots \rightedge (m+p,m)$. The vertex $(p+1,1)$ is labeled by $\alpha^r_p$ and $(m+p,m)$ is labeled by $\alpha_{m+p}^r$. By Lemma~\ref{l:taun_w0}, there is a black arrow $\alpha_p^r \rightarrow \alpha_{m+p}^r$ in $M_{(\bg^r,\bg^c)}$.
        \item The vertex $(n,1)$ has labels $\alpha_{n-1}^r$ and $\alpha_1^c$, and it corresponds to the black arrow $\alpha_{n-1}^r \rightarrow \alpha_1^c$.
    \end{itemize}
    The remaining cases in Lemma~\ref{l:taun_w0} can be analyzed in a similar way. Now, we observe the correspondence between the red arrows in $M_{(\bg^r,\bg^c)}$ and the edges of $\Gamma(\mathcal{S}^\prime)$ corresponding to $\btheta_p^r$ and~$\btheta_p^c$.
    \begin{itemize}
        \item For $\gamma_r^*(\alpha_q^r) = \alpha_p^r$, $\alpha_p^r \in \Gamma_1^r$, there is an edge $(q,m)\rightedge(p+1,1)$ in $\Gamma(\mathcal{S}^\prime)$. Since $(p+1,1)$ is labeled by $\alpha_p^r$ and $(q,m)$ is labeled by $\alpha_q^r$, there is a red arrow $\alpha_q^r \rightarrow \alpha_p^r$;
        \item Likewise, for $ \gamma_c(\alpha_p^c) = \alpha_q^c $, the edge $(n,p) \rightedge (1,q+1)$ corresponds to the red arrow $\alpha_p^c \rightarrow \alpha_q^c$.
    \end{itemize}
    It is now easy to see that $M_{(\bg^r,\bg^c)}$ has an alternating cycle if and only if $\Gamma(\mathcal{S}^\prime)$ has a cycle.
\end{proof}

\begin{remark}
   There is more to BD meanders than what is developed in this paper.
    In a sequel, as well as in the doctoral thesis of the first author, we will introduce and study the notion of a \emph{Belavin--Drinfeld action pair}, which is a certain twist of the $\mathbb{T}$-action pair studied in~\cite{lm-deform}. Deformations of such twists recover not only symmetric Poisson CGL extensions and aperiodic BD bivectors on $\Mat_{n\times m}(\mathbb{C})$, but a wider class of examples that includes $Q$-systems and Bottacin's Poisson structure on Hilbert schemes of points in the log-canonical plane studied in~\cite{mps-hilbert}. For instance, the latter is related to the BD triple
    $\bg = ([1,2k] \backslash\{k,k+1\}, [1,2k] \backslash \{k,2k\}, \gamma (i) = i+k \pmod{2k+1})$.
\end{remark}

\section{\texorpdfstring{Poisson rational map}{Poisson rational map}}\label{s:birat}

In this section, given a pair of BD quadruples $(\bbg^r,\bbg^c)$, we prove that there exists a $\mathbb{T}_{(\bg^r,\bg^c)}$-equivariant Poisson rational map
\begin{equation}\label{eq:map}
    \mathcal{U}_0 : \big(\Mat_{n\times m}(\mathbb{C}),(\pi_{(\bbg^r,\bbg^c)})_0\big) \dashrightarrow \big(\Mat_{n\times m}(\mathbb{C}),\pi_{(\bbg^r,\bbg^c)}\big).
\end{equation}
Moreover, we show that if $(\bg^r,\bg^c)$ is aperiodic, then the map $\mathcal{U}_0$ is birational.

Such maps are relevant to the theory of cluster algebras. Indeed, if $\mathcal{U}_0^{-1}$ exists, since all matrix entries of $\Mat_{n\times m}(\mathbb{C})$ are, by definition, log-canonical with respect to $(\pi_{(\bg^r,\bg^c)})_0$, their pullbacks by $\mathcal{U}_0^{-1}$ yield a rational log-canonical system for $(\Mat_{n\times m}(\mathbb{C}),\pi_{(\bbg^r,\bbg^c)})$. Extracting numerators, one obtains an initial cluster for a cluster structure in $\mathbb{C}[\Mat_{n\times m}]$ compatible with $\pi_{(\bbg^r,\bbg^c)}$. This approach is one of the tools used in work on the GSV conjecture (see, for instance, \cite{gsv-unified,multdual}, and the general theory in~\cite{zametka}). In the context of Poisson deformation theory, we believe that the existence of $\mathcal{U}_0$ sharpens the statement that $\pi_{(\bbg^r,\bbg^c)}$ is a Poisson deformation of its log-canonical term.

We conjecture that, more generally, given any $\mathbb{T}$-log-symplectic log-canonical Poisson bivector $\pi_0$ on $\mathbb{C}^n$ and a linearly independent subset $\mathcal{S}\subseteq \mathcal{S}(\pi_0)$ satisfying Property~\eqref{eq:W1}, there is $\mathbb{T}$-equivariant Poisson birational map between $(\mathbb{C}^n,\pi_0)$ and $(\mathbb{C}^n,\pi)$.

The section is organized as follows. In Subsection~\ref{s:gauss}, we review the notion of Gauss decomposition for rectangular matrices. In Subsection~\ref{s:biratstd}, given $\bbg_{\std}^r :=(\bg^r_{\std},r_0^r)$ and $\bbg_{\std}^c := (\bg^c_{\std},r_0^c)$, we construct the map $\mathcal{U}_0$ as the map
\begin{equation*}
    \mathcal{U}_0 : \big(\Mat_{n\times m}(\mathbb{C}),(\pi_{(\bbg^r_{\std},\bbg^c_{\std})})_0 \big)\dashrightarrow \big(\Mat_{n\times m}(\mathbb{C}),\pi_{(\bbg^r_{\std},\bbg^c_{\std})}\big)
\end{equation*}
(note that $(\pi_{(\bbg^r,\bbg^c)})_0 = (\pi_{(\bbg^r_{\std},\bbg^c_{\std})})_0$). In Subsection~\ref{s:biratgen}, we construct a rational map
\begin{equation*}
    \mathcal{U} : \big(\Mat_{n\times m}(\mathbb{C}),\pi_{(\bbg^r_{\std},\bbg^c_{\std})}\big) \dashrightarrow \big(\Mat_{n\times m}(\mathbb{C}),\pi_{(\bbg^r,\bbg^c)}\big).
\end{equation*}
For aperiodic BD pairs, we construct its inverse. In Subsection~\ref{s:biratp}, we prove that the map $\mathcal{U}$ is Poisson.
The map~\eqref{eq:map} is then obtained as the composition of the above two maps. It can then be used to show that $(\Mat_{n \times m}(\mathbb{C}), \pi_{(\bbg^r,\bbg^c)})$ has an open $\mathbb{T}_{(\bg^r,\bg^c)}$-leaf if $(\bg^r,\bg^c)$ is aperiodic.

\subsection{Gauss decomposition of matrices}\label{s:gauss}
In various constructions in this section, we use the notion of a Gauss decomposition (also known as an $LDU$ factorization) of a generic matrix $X \in \Mat_{n\times m}(\mathbb{C})$. Although this is a classical subject, we do not know of a reference for the particular construction needed in this paper. Indeed, the Gauss components of a rectangular matrix $X$ are not uniquely determined without additional normalization conditions. In this subsection, we specify our choice of the Gauss components. 

Setting $k:=\max(n,m)$, let us define an embedding
\begin{equation*}
    \imath : \Mat_{n\times m}(\mathbb{C}) \rightarrow \Mat_{k\times k}(\mathbb{C})
\end{equation*}
as follows.
\begin{itemize}
    \item For $n \leq m$, write
    \begin{equation*}
        \imath(X) := \begin{bmatrix}
            X^{[1,m-n]} & X^{[m-n+1,m]} \\
            I_{m-n} & 0
        \end{bmatrix};
    \end{equation*}
    \item For $n \geq m$, write
    \begin{equation*}
        \imath(X) := \begin{bmatrix}
            X_{[1,n-m]} & I_{n-m}\\
            X_{[n-m+1,n]} & 0
        \end{bmatrix}.
    \end{equation*}
\end{itemize}
We define a Gauss decomposition of a generic $X \in \Mat_{n\times m}(\mathbb{C})$ by using the embedding $\imath$ and the Gauss decomposition of square matrices, as follows.
\begin{itemize}
    \item For $n \leq m$, set
    \begin{equation*}
        [X]_- := [X^{[1,n]}]_-, \quad [X]_+ := [\imath(X)]_+;
    \end{equation*}
    \item For $n \geq m$, set
    \begin{equation*}
        [X]_- := [\imath(X)]_-, \quad [X]_+ := [X_{[1,m]}]_+.
    \end{equation*}
\end{itemize}
In both cases, we define $[X]_0$ as an element of $\Mat_{n\times m}(\mathbb{C})$ such that for $l:=\min(n,m)$,
\begin{equation*}
    ([X]_0)_{[1,l]}^{[1,l]} = \big[X^{[1,l]}_{[1,l]}\big]_0.
\end{equation*}
The rest of the entries of $[X]_0$ are zero. It is not difficult to verify that, as in the case of square matrices,
\begin{equation*}
    X = [X]_- [X]_0 [X]_+.
\end{equation*}
Moreover, since $\imath(X^T) = \imath(X)^T$, 
\begin{equation}\label{eq:transpose_gauss}
    [X^T]_- = [X]_+^T, \quad [X^T]_0 = [X]_0^T, \quad [X^T]_+ = [X]_-^T.
\end{equation}

Given a simple root $\alpha_p$ (in $A_{n-1}$ or $A_{m-1}$), denote by $U_p$ and $U_{-p}$ the one-parameter unipotent subgroups given by
\begin{equation*}
    U_p(t) := \exp(t e_{\alpha_p}), \quad U_{-p}(t) := \exp (t e_{-\alpha_p}).
\end{equation*}
Here, the root vector $e_{\alpha_p}$ (resp. $e_{-\alpha_p}$) is given by the matrix $e_{p,p+1}$ (resp. $e_{p+1,p}$). 

\begin{proposition}\label{p:gauss_perturb_matr}
    Assume $n \leq m$. In terms of the one-parameter subgroups, the Gauss decomposition of a generic $X \in \Mat_{n\times m}(\mathbb{C})$ has the following properties:
For any $p \in [1,n-1]$, $q \in [1,m-1]$, and generic $t \in \mathbb{C}$,
\begin{equation*}
    [U_{-p}(t) X]_- = U_{-p}(t) [X]_-, \quad [U_{-p}(t) X]_+ = [X]_+, \quad [U_{p}(t)X]_+ = U_{p}(t^{\prime}) [X]_+;
\end{equation*}
\begin{equation*}
    [XU_q(t)]_- =[X]_-, \quad 
    [XU_{-q}(t)]_- = \begin{cases}
        [X]_-U_{-q}(t^\prime) &\text{if} \ q < n;\\
        [X]_-&\text{if} \ q \geq n;
    \end{cases}
\end{equation*}
\begin{equation*}
[XU_q(t)]_+ =  \begin{cases}
   U_{n+q}(t^{\prime})[X]_+ U_q(t) & \text{if} \ q < m-n;\\
         [X]_+U_q(t) & \text{if} \ q \geq m-n.
    \end{cases}
\end{equation*} 
Here, the value of $t^\prime$ may differ from one equation to another, and it depends rationally on $X$ and $t$.
\end{proposition}
Identities for $n \geq m$ can be obtained from the identities for $n \leq m$ by using equations~\eqref{eq:transpose_gauss}.
\begin{proof}
    For a generic square matrix $Y$, we prove in Proposition~\ref{p:gauss_perturb} that $[U_p(t)Y]_+ = U_p(t^\prime)[Y]_+$ and $[YU_{-p}(t)]_- = [Y]_-U_{-p}(t^{\prime\prime})$. The precise expressions for $t^\prime$ and $t^{\prime\prime}$ are given in the appendix. Likewise, for any unipotent lower and upper triangular matrices $N_-$ and $N$, $[N_-YN]_- = N_-[Y]_-$ and $[N_-YN]_+ = [Y]_+ N$. Using these properties, it is not difficult to derive the corresponding identities for rectangular matrices.
    
    For example, let us prove the last equation. We see that
    \begin{equation*}
        \imath(XU_q(t)) = \begin{cases}
                U_{n+q}(-t)\imath(X) U_q(t) & \text{if} \ q < m-n;\\
                \imath(X) U_q(t) - te_{m,m-n+1} & \text{if} \ q = m-n;\\
                 \imath(X)U_q(t) & \text{if} \ q > m-n.
            \end{cases}
    \end{equation*}
    By definition, the matrix $\imath(XU_q(t))$ is square, and the $[\, \cdot\,]_+$ component of a square matrix does not depend on its last row. Taking $[\,\cdot\,]_+$ and using the properties of the Gauss decomposition for square matrices listed in the previous paragraph, we obtain the equation.
\end{proof}

Explicit expressions for the entries of the Gauss decomposition of a square matrix $Y$ are well-known (see, for instance, \cite{gantmacher}). These are given by ratios of minors of $Y$. Using these expressions, it is easy to obtain explicit expressions for the rectangular case. Let us only mention the corresponding divisors. For $n \leq m$, a matrix $X \in \Mat_{n\times m}(\mathbb{C})$ has the Gauss decomposition if and only if it belongs to the nonzero locus of the following minors:
\begin{equation*}
    \det X_{[1,i]}^{[1,i]},\ i \in [1,n]; \quad \det X_{[1,n]}^{[j+1,n+j]}, \ j \in [1,m-n].
\end{equation*}
Transposing, one obtains the divisors for $n \geq m$. 

\subsection{From log-canonical to standard}\label{s:biratstd}
Let us denote by $\mathbb{T}_{\std}$ the BD torus corresponding to the trivial BD triples. In this subsection, we describe a $\mathbb{T}_{\std}$-equivariant Poisson birational map
\begin{equation}\label{eq:U0std}
    \mathcal{U}_0 : \big(\Mat_{n\times m}(\mathbb{C}),(\pi_{(\bbg^r_{\std},\bbg^c_{\std})})_0 \big)\dashrightarrow \big(\Mat_{n\times m}(\mathbb{C}),\pi_{(\bbg^r_{\std},\bbg^c_{\std})}\big),
\end{equation}
where $\bbg^r_{\std} = (\bg_{\std},r_0^r)$ and $\bbg^c_{\std} = (\bg_{\std},r_0^c)$.

The map $\mathcal{U}_0$ is easier to understand in terms of the components of its Gauss decomposition. Let us set $l:=\min(n,m)$ and write $\mathcal{U}_0(X) = [\mathcal{U}_0(X)]_- [\mathcal{U}_0(X)]_0 [\mathcal{U}_0(X)]_+$. 
Then the components of $\mathcal{U}_0(X)$ further factor as follows:
\begin{align*}
    &[\mathcal{U}_0(X)]_- = [\mathcal{U}_0(X)]_-^{(1)}[\mathcal{U}_0(X)]_-^{(2)}\cdots [\mathcal{U}_0(X)]_-^{(l)};\\
    &[\mathcal{U}_0(X)]_+ = [\mathcal{U}_0(X)]_+^{(l)}[\mathcal{U}_0(X)]_+^{(l-1)}\cdots [\mathcal{U}_0(X)]_+^{(1)}.
\end{align*}
Here, each matrix $[\mathcal{U}_0(X)]_-^{(k)}$ (resp. $[\mathcal{U}_0(X)]_+^{(k)}$) is a unipotent lower (resp. upper) triangular matrix of size $n\times n$ (resp. $m \times m$). Their nonzero off-diagonal entries are given as follows: 
\begin{align*}
    &\big([\mathcal{U}_0(X)]_-^{(k)}\big)_{[k,n]}^{[k,n]} = \pr_{<} \begin{bmatrix}
        x_{kk}^{-1}\\
        \vdots \\ x_{nk}^{-1}
    \end{bmatrix} \begin{bmatrix}
        x_{kk} & \ldots& x_{nk}
    \end{bmatrix}, \quad k \in [1,l];\\
    &\big([\mathcal{U}_0(X)]_+^{(k)}\big)_{[k,m]}^{[k,m]} = \pr_{>} \begin{bmatrix}
        x_{kk}^{-1}\\
        \vdots \\ x_{km}^{-1}
    \end{bmatrix} \begin{bmatrix}
        x_{kk} & \ldots& x_{km}
    \end{bmatrix},  \quad k \in [1,l],
\end{align*}
where $\pr_{<}$ (resp. $\pr_{>}$) denotes the projection onto strictly lower (resp. upper) triangular matrices of the appropriate size. The component $[\mathcal{U}_0(X)]_0$ is given by
\begin{equation*}
    ([\mathcal{U}_0(X)]_0)_{[1,l]}^{[1,l]} = \mathrm{diag}(x_{11}, \ldots,x_{ll}),
\end{equation*}
and all its remaining entries are zero. One can verify that the Gauss components of $\mathcal{U}_0(X)$ (in the sense of the previous subsection) are indeed given by the above equations. However, we will not need this fact.

\begin{example}
    For $n = m = 2$, the map $\mathcal{U}_0$ is given by
    \begin{equation*}
        \mathcal{U}_0(X) = \begin{bmatrix}
            1 & \\
            \frac{x_{21}}{x_{11}} & 1
        \end{bmatrix} \begin{bmatrix}
            x_{11} &\\
            & x_{22}
        \end{bmatrix} \begin{bmatrix}
            1 & \frac{x_{12}}{x_{11}}\\
            & 1
        \end{bmatrix}.
    \end{equation*}
\end{example}
\begin{example}\label{ex:U0_23}
    For $n=2$, $m = 3$, the map $\mathcal{U}_0$ is given by
    \begin{equation*}
        \mathcal{U}_0(X) = \begin{bmatrix}
            1 & \\
            \frac{x_{21}}{x_{11}} & 1
        \end{bmatrix} \begin{bmatrix}
            x_{11} & &\;\; \\
            & x_{22} &\;\;
        \end{bmatrix}
        \begin{bmatrix}
            1 & &\\
            & 1 & \frac{x_{23}}{x_{22}}\\
            & & 1
        \end{bmatrix}
        \begin{bmatrix}
            1 & \frac{x_{12}}{x_{11}} & \frac{x_{13}}{x_{11}}\\
            & 1 & \frac{x_{13}}{x_{12}}\\
            & & 1
        \end{bmatrix}
    \end{equation*}
\end{example}

Let us introduce the following notation for minors of $X \in \Mat_{n\times m}(\mathbb{C})$:
\begin{equation*}
    \Delta_{ij}(X) = \begin{dcases}
        \det X_{[1,i]}^{[a,j]} &\text{if}\ i \leq j, \ a = j-i+1;\\[5pt]
        \det X_{[b,i]}^{[1,j]} &\text{if}\ i \geq j, \ b = i-j+1.
    \end{dcases}
\end{equation*}
Denoting by $\{\cdot,\cdot\}_{\std}$ the standard Poisson bracket on $\Mat_{n\times m}(\mathbb{C})$ (that is, corresponding to the trivial BD quadruples), it is well known that the minors $\Delta_{ij}$ are log-canonical. Writing $\Delta_{ij}(X) = \det X_{I}^J$ and $\Delta_{kl}(X) = \det X_{K}^L$, the Poisson brackets between them are given by
\begin{equation}\label{eq:bracket}
    \{\Delta_{ij},\Delta_{kl}\}_{\std}(X) = \frac{1}{2}\epsilon_{ij,kl} \big( |I\cap K| - |J \cap L|\big)\Delta_{ij}(X)\Delta_{kl}(X)
\end{equation}
where $\epsilon_{ij,kl} \in \{-1,1\}$ is the sign given by 
\begin{equation*}
    \epsilon_{ij,kl} = \begin{cases}
        \sgn((i-j)-(k-l)) & \text{if} \ \sgn((i-j)(k-l))> -1;\\
        1 & \text{if}\ i > j, \ k \leq l;\\
        -1 &\text{otherwise.}
    \end{cases}
\end{equation*}

Equation~\eqref{eq:bracket} can be viewed as a particular case of~\cite[Proposition 5.8]{gy-cgl} or \cite[Proposition 4.19]{gsv-book}. The Poisson bracket between the matrix entries is given by
\begin{equation*}
    \{x_{ij},x_{kl}\}_{\std} = \frac{1}{2}(\sgn(l-j)+\sgn(k-i))x_{il}x_{kj}.
\end{equation*}
By Proposition~\ref{p:log-part}, its log-canonical part is given by
\begin{equation}\label{eq:std0xij}
    \{x_{ij},x_{kl}\}_{\std,0} = \frac{1}{2}(\delta_{ik}\sgn(l-j)+\delta_{jl}\sgn(k-i))x_{ij}x_{kl}.
\end{equation}
\begin{proposition}\label{p:U0}
    For any choice of $(r_0^r,r_0^c)$, the rational map
    \begin{equation*}
        \mathcal{U}_0 : \big(\Mat_{n\times m}(\mathbb{C}),(\pi_{(\bbg^r_{\std},\bbg^c_{\std})})_0 \big)\dashrightarrow \big(\Mat_{n\times m}(\mathbb{C}),\pi_{(\bbg^r_{\std},\bbg^c_{\std})}\big)
    \end{equation*}
    is a $\mathbb{T}_{\std}$-equivariant Poisson birational map. Moreover, the pullbacks of $\Delta_{ij}$ under $\mathcal{U}_0$ are given~by
    \begin{equation}\label{eq:U0pullback}
        \mathcal{U}_0^*(\Delta_{ij}(X)) = \begin{dcases}
            \prod_{p=1}^i x_{p,a+p-1} &\text{if} \ i \leq j, \ a := j-i+1;\\
            \prod_{p=1}^j x_{b+p-1,p} &\text{if} \ i \geq j, \ b := i-j+1.
        \end{dcases}
    \end{equation}
\end{proposition}
In other words, equations~\eqref{eq:U0pullback} state that the pullback of $\Delta_{ij}(X)$ is the product of the diagonal elements of the corresponding submatrix of $X$. In the terminology of~\cite{gy-cgl}, these products are the \emph{leading terms} of $\Delta_{ij}(X)$.

\begin{figure}[t]
    \centering
    \begin{tikzpicture}
        \filldraw[black] (0,0) circle (\bdvsize) node[anchor=east]{$1$};
        \filldraw[black] (0,0.75) circle (\bdvsize) node[anchor=east]{$2$};
        \filldraw[black] (0,1.5) circle (\bdvsize) node[anchor=east]{$3$};
        \filldraw[black] (0,2.25) circle (\bdvsize) node[anchor=east]{$4$};

        \filldraw[black] (7,0) circle (\bdvsize) node[anchor=west]{$1$};
        \filldraw[black] (7,0.75) circle (\bdvsize) node[anchor=west]{$2$};
        \filldraw[black] (7,1.5) circle (\bdvsize) node[anchor=west]{$3$};
        \filldraw[black] (7,2.25) circle (\bdvsize) node[anchor=west]{$4$};

        \draw (0,0) -- (7,0);
        \draw (0,0.75) -- (7,0.75);
        \draw (0,1.5) -- (7,1.5);
        \draw (0,2.25) -- (7,2.25);

        \draw (5.5,1.5) -- (6,2.25);
        \draw (4.5,0.75) -- (5,1.5);
        \draw (3.5,0) -- (4,0.75);

        \node at (6.15,1.85) {$\frac{x_{14}}{x_{13}}$};
        \node at (5.15,1.1) {$\frac{x_{13}}{x_{12}}$};
        \node at (4.15,0.35) {$\frac{x_{12}}{x_{11}}$};       

        \filldraw[black] (6,2.25) circle (\bdvsize);
        \filldraw[black] (5.5,1.5) circle (\bdvsize);
        \filldraw[black] (5,1.5) circle (\bdvsize);
        \filldraw[black] (4.5, 0.75) circle (\bdvsize);
        \filldraw[black] (4,0.75) circle (\bdvsize);
        \filldraw[black] (3.5,0) circle (\bdvsize);

        \draw (2.5,1.5) -- (3,2.25);
        \draw (1.5,0.75) -- (2,1.5);

        \node at (3.15,1.85) {$\frac{x_{24}}{x_{23}}$};
        \node at (2.15,1.1) {$\frac{x_{23}}{x_{22}}$};
        
        \filldraw[black] (3,2.25) circle (\bdvsize);
        \filldraw[black] (2.5,1.5) circle (\bdvsize);
        \filldraw[black] (1.5,0.75) circle (\bdvsize);
        \filldraw[black] (2,1.5) circle (\bdvsize);

        \draw (0.5,1.5) -- (1,2.25);
        \filldraw[black] (0.5,1.5) circle (\bdvsize);
        \filldraw[black] (1,2.25) circle (\bdvsize);        
        \node at (1.15,1.85) {$\frac{x_{34}}{x_{33}}$};        
    \end{tikzpicture}
    \caption{Planar network for $[\mathcal{U}_0(X)]_+$ for $n\geq 3$ and $m=4$.}
    \label{f:network}
\end{figure}
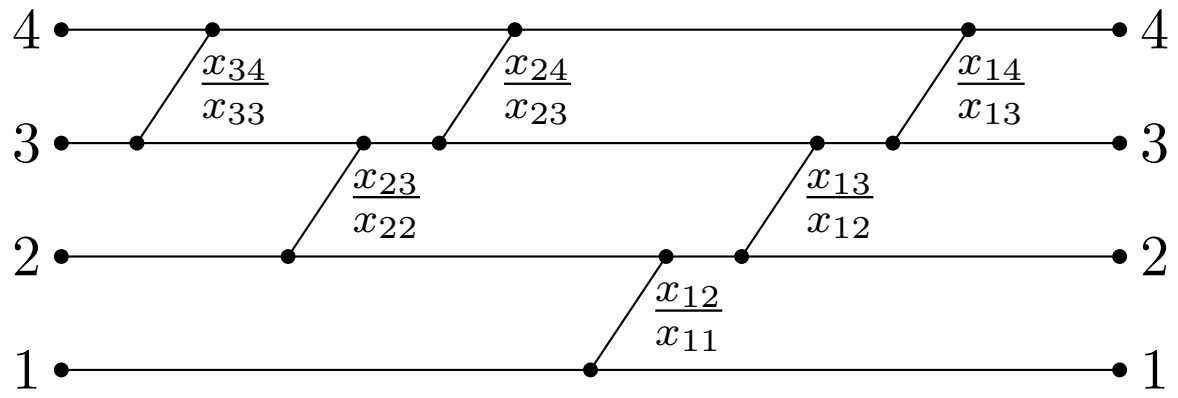

\begin{proof}
    Let us first show that equations~\eqref{eq:U0pullback} imply that $\mathcal{U}_0$ is a $\mathbb{T}_{\std}$-equivariant Poisson birational map. Recall that the minors $\Delta_{ij}(X)$, for $i\in [1,n]$ and $j \in [1,m]$, generate the field of rational functions $\mathbb{C}(\Mat_{n\times m})$. By~\eqref{eq:U0pullback}, we see that $\mathcal{U}_0^*(\Delta_{1,j}) = x_{1,j}$ and $\mathcal{U}_0^*(\Delta_{i,1}) = x_{i,1}$. For $i,j>1$, 
    \begin{equation*}
        \mathcal{U}_0^*(\Delta_{ij}(X) \Delta_{i-1,j-1}(X)^{-1}) = x_{ij}.
    \end{equation*}
    Therefore, $\mathcal{U}_0$ is a birational map. To see that $\mathcal{U}_0$ is $\mathbb{T}_{\std}$-equivariant, observe that, by equations~\eqref{eq:U0pullback}, $\Delta_{ij}(X)$ and $\mathcal{U}_0^*(\Delta_{ij}(X))$ have the same $\mathbb{T}_{\std}$-weights.

    To show that~\eqref{eq:U0pullback} implies $\mathcal{U}_0$ is Poisson, let us first consider the case of the standard Poisson brackets. It suffices to show that $\mathcal{U}_0^*\{\Delta_{ij},\Delta_{kl}\}_{\std} = \{\mathcal{U}_0^*\Delta_{ij},\mathcal{U}_0^*\Delta_{kl}\}_{\std,0}$. Since both sides are log-canonical, it suffices to compare the corresponding coefficients. Let us consider the case $i\leq j$ and $k \leq l$. In this case, by~\eqref{eq:bracket},
    \begin{equation*}
        \{\Delta_{ij},\Delta_{kl}\}_{\std}(X) = \frac{1}{2}\sgn((i-j)-(k-l))\big(\min(i,k)-|J \cap L|\big)\Delta_{ij}(X)\Delta_{kl}(X),
    \end{equation*}
    where $|J \cap L| = |[j-i+1,j]\cap [l-k+1,l]|$.
    On the other hand, by equations~\eqref{eq:std0xij}-\eqref{eq:U0pullback}, the coefficient of $\{\mathcal{U}_0^*\Delta_{ij},\mathcal{U}_0^*\Delta_{kl}\}_{\std,0}$ is given by
    \begin{equation*}\begin{split}
         \frac{1}{2}\sum_{p_1=1}^i& \sum_{p_2=1}^k (\delta_{p_1,p_2}\sgn(l-k-(j-i))+\delta_{j-i+p_1,l-k+p_2}\sgn(p_2-p_1)) = \\ &= \frac{1}{2}\sgn((i-j)-(k-l))\big( \min(i,k) - \sum_{p_1=1}^i \sum_{p_2=1}^k \delta_{j-i+p_1,l-k+p_2}\big).
         \end{split}
    \end{equation*}
    We see that the above coefficient coincides with the coefficient of $\{\Delta_{ij},\Delta_{kl}\}_{\std}$. The other cases of $i,j,k,l$ can be handled similarly. Hence $\mathcal{U}_0$ is Poisson for the standard BD quadruples. 

    To show that $\mathcal{U}_0$ is Poisson with respect to the BD quadruples $(\bbg_{\std}^r,\bbg_{\std}^c)$, observe that, since $\mathcal{U}_0$ is $\mathbb{T}_{\std}$-equivariant, $(\mathcal{U}_0)_*((s_0^c)^L) = (s_0^c)^L$ and $(\mathcal{U}_0)_*((s_0^r)^R) = (s_0^r)^R$.
    
    It remains to prove equations~\eqref{eq:U0pullback}. Let us focus on the case $i \leq j$. In this case, we see that 
    \begin{equation*}
        \Delta_{ij}(N_-X) = \Delta_{ij}(X), \quad \Delta_{ij}(TX) = \det(T)_{[1,i]}^{[1,i]} \Delta_{ij}(X)
    \end{equation*}
    for any $N_- \in \mathcal{N}_-^{(n)}$ and an invertible diagonal matrix $T$. Therefore,
    \begin{equation*}
        \mathcal{U}_0^*(\Delta_{ij}(X)) =\Big(\prod_{p=1}^i x_{pp}\Big) \Delta_{ij}([\mathcal{U}_0(X)]_+).
    \end{equation*}
    In order to compute the minor $\Delta_{ij}([\mathcal{U}_0(X)]_+)$, we draw a planar network associated with $[\mathcal{U}_0(X)]_+$, and then employ Lindstr\"om's lemma~\cite{lindstr-lemma}.
    For an overview of this technique, we refer the reader to~\cite{fz-total}. 

    Briefly, let $M$ be an $m \times m$ matrix that can be written as the product of \emph{elementary Jacobi matrices} (that is, matrices of the form $(I_m + t e_{p,p+1})$, where $I_m$ is the $m \times m$ identity matrix). To draw a planar network for $M$, one starts with two sets of vertices labeled by $[1,m]$ (called \emph{sources} and \emph{sinks}). One connects each sink $i$ to the source $i$ by a horizontal edge of weight $1$. Thus one forms $m$ horizontal \emph{levels}. For each $I_m + t e_{p,p+1}$ in the product, one places a directed edge of weight $t$ connecting level $p$ to level $p+1$. Given a path connecting a source $i$ to a sink $j$, its weight is the product of the weights of all edges along the path. More generally, the weight of a collection of vertex-disjoint paths is the product of the weights of individual paths. Lindstr\"om's lemma states that a minor $\det M_{I}^J$ equals the sum of weights of all vertex-disjoint collections of paths connecting the sources $I$ to the sinks $J$.
    
    Observe that for each $k \in [1,l]$, the matrix $[\mathcal{U}_0(X)]_+^{(k)}$ can be factored into the product of elementary Jacobi matrices:
    \begin{equation*}
        [\mathcal{U}_0(X)]_+^{(k)} = \Big(I_m + \frac{x_{k,k+1}}{x_{kk}}e_{k,k+1}\Big)\Big(I_m + \frac{x_{k,k+2}}{x_{k,k+1}}e_{k+1,k+2}\Big) \cdots \Big(I_m + \frac{x_{k,m}}{x_{k,m-1}}e_{m-1,m}\Big).
    \end{equation*}
    Since the matrix $[\mathcal{U}_0(X)]_+$ is the product of $[\mathcal{U}_0(X)]_+^{(k)}$ for $k \in [1,l]$, it can therefore be represented by a planar network (see Figure~\ref{f:network} for an illustration).  

    It is easy to see that for $I = [1,i]$ and $J = [a,j]$, $a = j-i+1$, there is a unique collection of vertex-disjoint paths connecting the sources $I$ to the sinks $J$ (see also \cite[Proof of Lemma 6]{fz-total}). In such a collection, the path connecting $p \in [1,i]$ to $q=a+p-1$ has weight 
    \begin{equation*}
        \frac{x_{p,p+1}}{x_{pp}} \cdot \frac{x_{p,p+2}}{x_{p,p+1}}\cdots \frac{x_{pq}}{x_{p,q-1}} = \frac{x_{pq}}{x_{pp}}.
    \end{equation*}
    Therefore, 
    \begin{equation*}
        \mathcal{U}_0^*(\Delta_{ij}(X)) = \Big( \prod_{p=1}^i x_{pp}\Big) \Delta_{ij}([\mathcal{U}_0(X)]_+) = \Big( \prod_{p=1}^i x_{pp}\Big)\cdot \Big( \prod_{p=1}^i \frac{x_{p,a+p-1}}{x_{pp}} \Big) = \prod_{p=1}^i x_{p,a+p-1}.
    \end{equation*}
    Thus, equations~\eqref{eq:U0pullback} hold.
\end{proof}

\begin{remark} 
    Equations~\eqref{eq:U0pullback} imply a stronger statement than birationality: the map $\mathcal{U}_0$ induces an isomorphism between the localizations
    \begin{equation*}
        \mathcal{U}_0^* : \mathbb{C}[\Mat_{n\times m}][\Delta_{ij}^{\pm 1} \mid i \in [1,n-1],\ j \in [1,m-1] ] \xrightarrow{\sim} \mathbb{C}[\Mat_{n\times m}][x_{ij}^{\pm 1} \mid i \in [1,n-1],\ j \in [1,m-1] ].
    \end{equation*}
    In particular, if we consider a cluster structure in $\mathbb{C}[\Mat_{n\times m}]$ compatible with $(\pi_{\std})_0$ in which all variables are frozen and are given by $x_{ij}$, then $\mathcal{U}_0$ is a birational quasi-isomorphism in the sense of~\cite{zametka}. The variables $\Delta_{ij}$ for $i \in [1,n-1]$ and $j \in [1,m-1]$ are the initial cluster variables in $\mathbb{C}[\Mat_{n\times m}]$, and the variables $\Delta_{nj}$ and $\Delta_{im}$ for $j \in [1,m]$ and $i \in [1,n]$ are frozen. The divisors of $\mathcal{U}_0$ are precisely the defrosted variables. 
    
    We conjecture that this is a general phenomenon. More precisely, in the setup of Subsection~\ref{s:background_deform}, let $\pi_0$ be a $\mathbb{T}$-log-symplectic log-canonical bivector, and let $\pi$ be a normalized $\mathcal{S}$-admissible $\mathbb{T}$-invariant Poisson deformation of $\pi_0$ along $\pi_1 \in \mathfrak{X}^2(\mathbb{C}^n)^{\mathcal{S}}$, with $\mathcal{S}\subseteq \mathcal{S}(\pi_0)$ satisfying Property (W0). We conjecture that $\mathbb{C}[x_1,\ldots,x_n]$ has the structure of an upper cluster algebra compatible with $\pi$ and $\mathbb{T}$ in the sense of \cite{gsv-book}. Moreover, we conjecture that there exists a $\mathbb{T}$-equivariant Poisson birational map $(\mathbb{C}^n,\pi_0)\dashrightarrow(\mathbb{C}^n,\pi)$ whose reduced pole divisor is precisely the zero divisor of defrosted variables.
\end{remark}

\subsection{From standard to nonstandard BD}\label{s:biratgen}
In this subsection, given a pair of BD quadruples $(\bbg^r,\bbg^c)$, we construct a $\mathbb{T}_{(\bg^r,\bg^c)}$-equivariant Poisson rational map
\begin{equation*}
    \mathcal{U} : (\Mat_{n \times m}(\mathbb{C}),\pi_{(\bbg_{\std}^r,\bbg_{\std}^c)}) \dashrightarrow (\Mat_{n \times m}(\mathbb{C}), \pi_{(\bbg^r,\bbg^c)}).
\end{equation*}
We prove that if $(\bg^r,\bg^c)$ is aperiodic, the map $\mathcal{U}$ is birational. We construct the inverse map by means of a finite recursive procedure. The constructions and proofs are given in Lie-theoretic terms; thus, replacing $\Mat_{n\times m}(\mathbb{C})$ with a reductive complex algebraic group $G$, we obtain a formula for the inverse of $\mathcal{U}:(G,\pi_{(\bbg^r_{\std},\bbg^c_{\std})}) \dashrightarrow (G,\pi_{(\bbg^r,\bbg^c)})$. In the next subsection, we will prove that the map $\mathcal{U}$ is Poisson.

Given a BD triple $\bg:=(\Gamma_1,\Gamma_2,\gamma)$ (in any type), let us recall that $\gamma$ can be extended to Lie algebra homomorphisms $\gamma : \mathfrak{n}_{\pm}\rightarrow \mathfrak{n}_{\pm}$, and thus can be integrated to group homomorphisms $\tilde{\gamma} : \mathcal{N}_\pm \rightarrow \mathcal{N}_\pm$ between the corresponding unipotent subgroups. The same applies to the adjoint map $\gamma^* : \Gamma_2 \rightarrow \Gamma_1$. We will not need explicit formulas in this paper;
nevertheless, in the case of matrices, the reader may find them in~\cite[Section 4.1]{multdual}.

For a given pair $(\bg^r,\bg^c)$, define the maps
\begin{align*}
    &\rho_r^* : \mathcal{N}_-^{(n)} \rightarrow \mathcal{N}_-^{(n)}, \quad \rho_r^*(N_-) = \prod_{j\geq 1}^{\leftarrow} (\tilde{\gamma}^*_r)^j(N_-) =   \cdots (\tilde{\gamma}^*_r)^2(N_-)\tilde{\gamma}^*_r(N_-);\\
    &\rho_c: \mathcal{N}^{(m)}\rightarrow \mathcal{N}^{(m)}, \quad \rho_c(N) = \prod_{j \geq 1}^{\rightarrow}(\tilde{\gamma}_c)^j(N) = \tilde{\gamma}_c(N) (\tilde{\gamma}_c)^2(N) \cdots .
\end{align*}
Note that, since the BD triples are nilpotent, the above products are finite.

Lastly, let $W_0^r$ (resp. $W_0^c$) be the $n \times n$ (resp. $m \times m$) permutation matrix corresponding to the longest permutation. We choose the nonzero entries to be equal to $1$, so that
\begin{equation*}
    (W_0^r)_{ij} = \delta_{i,n-j+1}, \quad (W_0^c)_{ij} = \delta_{i,m-j+1}.
\end{equation*}

The rational map $\mathcal{U}$ is given by
\begin{equation}\label{eq:Ubirat}
    \mathcal{U}(X) := \rho_r^*([XW_0^c]_-) X \rho_c([W_0^r X]_+).
\end{equation}

\begin{example}
    Let $n \geq 1$, $m=3$, $\bg^r = \bg_{\std}$, and $\bg^c = (\{\alpha_1\},\{\alpha_2\},\alpha_1\mapsto \alpha_2)$. Then the map $\mathcal{U}$ is given by
    \begin{equation*}
        \mathcal{U}(X) = X\big(1 + \frac{x_{n2}}{x_{n1}}e_{23}\big).
    \end{equation*}
\end{example}

\begin{example}
    Let $n=4$, $m=3$, $\bg^r = (\{\alpha_1,\alpha_2\},\{\alpha_2,\alpha_3\},\alpha_i \mapsto \alpha_{i+1})$, \mbox{$\bg^c = (\{\alpha_2\},\{\alpha_1\},\alpha_2\mapsto \alpha_1)$}. Then the map $\mathcal{U}$ is given by
    \begin{equation*}
        \mathcal{U}(X) = \begin{bmatrix}
            1 & & &\\
            \frac{\det X_{[1,2]\cup\{4\}}}{\det X_{[1,3]}} & 1 & &\\
            & & 1 &\\
            & & & 1
        \end{bmatrix}
        \begin{bmatrix}
            1 & & &\\
            \frac{\det X_{\{1,3\}}^{[1,2]}}{\det X_{[1,2]}^{[2,3]}} & 1 & &\\[10pt]
            \frac{\det X_{\{1,4\}}^{[2,3]}}{\det X_{[1,2]}^{[2,3]}} &  \frac{\det X_{[1,2]\cup\{4\}}}{\det X_{[1,3]}} & 1 &\\
            & & & 1
        \end{bmatrix}X \begin{bmatrix}
            1 & \frac{\det X_{[1,2]}^{\{1,3\}}}{\det X_{[1,2]}^{[1,2]}} &\\
            & 1 & \\
            & & 1
        \end{bmatrix}
    \end{equation*}
\end{example}

We first address the $\mathbb{T}_{(\bg^r,\bg^c)}$-equivariance of the above map. This will follow from the following properties of the maps $\rho_r^*$ and $\rho_c$:
\begin{lemma}\label{l:T_equiv_rho}
    Let $(\bg^r,\bg^c)$ be a pair of BD triples. Then for any $T_1 \in \mathcal{H}_{\bg^r}$, $T_2 \in \mathcal{H}_{\bg^c}$, $N_- \in \mathcal{N}_-^{(n)}$, and $N \in \mathcal{N}^{(m)}$,
    \begin{equation*}
        \rho_r^*(\Ad_{T_1} N_-) = \Ad_{T_1} \rho_r^*(N_-), \quad \rho_c(\Ad_{T_2} N) = \Ad_{T_2}\rho_c(N).
    \end{equation*}
\end{lemma}
\begin{proof}
For a simple root $\alpha$ in $A_{n-1}$, let $U_{-\alpha}(t):=\exp(t e_{-\alpha})$ be the corresponding one-parameter subgroup. For any $\alpha \in \Gamma_2^r$, since  $\Ad_{T_1} U_{-\alpha}(t) = U_{-\alpha}(T_1^{-\alpha}t)$ and  $T_1^{-\alpha} = T_1^{-\gamma_r^*(\alpha)}$ (by definition of $\mathcal{H}_{\bg^r}$), we see that
\begin{equation*}
    \tilde{\gamma}_r^*(\Ad_{T_1} U_{-\alpha}(t)) = \tilde{\gamma}_r^*(U_{-\alpha}(T_1^{-\alpha}t)) = U_{-\gamma_r^*(\alpha)}(T_1^{-\gamma_r^*(\alpha)}t) = \Ad_{T_1}\tilde{\gamma}_r^*(U_{-\alpha}(t)).
\end{equation*}
Thus the first equation follows. The second equation can be proved in a similar way.
\end{proof}

\begin{proposition}\label{p:Tequivar}
    The rational map $\mathcal{U}$ is $\mathbb{T}_{(\bg^r,\bg^c)}$-equivariant.
\end{proposition}
\begin{proof}
    For any $T_1 \in \mathcal{H}_{\bg^r}$, by Lemma~\ref{l:T_equiv_rho}, we see that
    \begin{equation*}
        \rho_r^*([T_1 X W_0^c]_-) = \rho_r^*(\Ad_{T_1} [XW_0^c]_-) = \Ad_{T_1}\rho_r^*([XW_0^c]_-).
    \end{equation*}
    Since $\rho_c([W_0^r T_1X]_+) = \rho_c([W_0^rX]_+)$, we see that
    \begin{equation*}
        \mathcal{U}(T_1X) = \rho_r^*([T_1XW_0^c]_-) T_1 X \rho_c([W_0^rT_1X]_+) = T_1 \mathcal{U}(X).
    \end{equation*}
    A similar proof applies to the action of $\mathcal{H}_{\bg^c}$ on the right. Thus the statement holds.
\end{proof}

Now, let us describe the idea of how to find the inverse of $\mathcal{U}$. Rewriting~\eqref{eq:Ubirat} as
\begin{equation*}
    X = \rho_r^*([XW_0^c]_-)^{-1} \mathcal{U}(X) \rho_c([W_0^rX]_+)^{-1}
\end{equation*}
and assuming the inverse exists, let us substitute $X = \mathcal{U}^{-1}(Y)$ in the above equation. We obtain the equation
\begin{equation*}
    \mathcal{U}^{-1}(Y) = \rho_r^*([\mathcal{U}^{-1}(Y)W_0^c]_-)^{-1} Y \rho_c([W_0^r \mathcal{U}^{-1}(Y)]_+)^{-1}.
\end{equation*}
Define the rational maps
\begin{equation*}
    \delta_r^*(Y) := \rho_r^*([\mathcal{U}^{-1}(Y)W_0^c]_-), \quad \delta_c(Y):=\rho_c([W_0^r \mathcal{U}^{-1}(Y)]_+).
\end{equation*}
Observe that the maps $\rho_r^*$ and $\rho_c$ satisfy the equations
\begin{equation*}
    \rho_r^*(N_-) = \tilde{\gamma}_r^*(\rho_r^*(N_-)N_-), \quad \rho_c(N) = \tilde{\gamma}_c(N \rho_c(N)), \quad N_- \in \mathcal{N}^{(n)}_-, \ N \in \mathcal{N}^{(m)}.
\end{equation*}
Substituting $N_- := [\mathcal{U}^{-1}(Y)W_0^c]_-$ and $N:= [W_0^r\mathcal{U}^{-1}(Y)]_+$, we arrive at equations for $\delta_r^*$ and $\delta_c$:
\begin{align}
    &\delta_r^*(Y) = \tilde{\gamma}_r^*(\delta_r^*(Y)\cdot[\delta_r^*(Y)^{-1}Y\delta_c(Y)^{-1}W_0^c]_-);\label{eq:deltas_r} \\
    &\delta_c(Y) = \tilde{\gamma}_c([W_0^r \delta_r^*(Y)^{-1}Y\delta_c(Y)^{-1}]_+\cdot \delta_c(Y)). \label{eq:delta_c}
\end{align}
Thus, constructing the map $\mathcal{U}^{-1}$ amounts to solving the above equations for $\delta_r^*$ and $\delta_c$. 

Note that equations~\eqref{eq:deltas_r}-\eqref{eq:delta_c} simplify in the case $n = m$. Indeed, by using properties of the Gauss decomposition, one can cancel out $\delta_r^*(Y)$ (resp. $\delta_c(Y)$) on the right-hand side of equation~\eqref{eq:deltas_r} (resp. \eqref{eq:delta_c}) and obtain the equations
\begin{equation*}
    \delta_r^*(Y) = \tilde{\gamma}_r^*([Y\delta_c(Y)^{-1}W_0^c]_-), \quad \delta_c(Y) = \tilde{\gamma}_c([W_0^rY\delta_r^*(Y)^{-1}]_+).
\end{equation*}

The above equations have a natural counterpart in the case of a reductive complex algebraic group $G$. In this case, $Y$ is replaced by a generic element $y \in G$, and the elements $W_0^c$ and $W_0^r$ are replaced by a representative $\overline{w_0}$ of the longest Weyl group element.

For the results that follow, it will be convenient for us to introduce the rational maps
\begin{equation}\label{eq:kappas}
    \kappa_r(X) := [XW_0^c]_-, \quad \kappa_c(X) := [W_0^r X]_+.
\end{equation}
The behavior of these maps under left and right multiplication of $X$ by unipotent subgroups exhibits an important connection to the combinatorics of the affine BD meander. We summarize these properties in the lemma below.

\begin{lemma}\label{l:kappa_invar}
    Assume that $n \leq m$. Let $X \in \Mat_{n\times m}(\mathbb{C})$ be a generic element, and let $\Delta^\prime \subseteq A_{n-1}$ and $\Delta^{\prime\prime}\subseteq A_{m-1}$ be any subsets of simple roots. The following statements hold:
    \begin{enumerate}[1)]
        \item For a generic $N \in \mathcal{N}^{(m)}(\Delta^{\prime\prime})$, there exists $N_-^\prime \in \mathcal{N}_-^{(n)}(\imath_r^{-1}\tau_n \imath_c(\Delta^{\prime\prime}))$ such that 
        \begin{equation*}
            \kappa_r(XN) = \kappa_r(X)N_-^\prime.
        \end{equation*}
        For any $N_- \in \mathcal{N}^{(n)}_-$,
        \begin{equation*}
            \kappa_r(N_-X) = N_- \kappa_r(X).
        \end{equation*}
        \item For a generic $N_- \in \mathcal{N}_-^{(n)}(\Delta^\prime)$, there is $N^\prime \in \mathcal{N}^{(m)}(\imath_c^{-1}\tau_n\imath_r(\Delta^\prime))$ such that
        \begin{equation*}
            \kappa_c(N_- X) = N^\prime \kappa_c(X).
        \end{equation*}
        For any a generic $N \in \mathcal{N}^{(m)}(\Delta^{\prime\prime})$, there is $N^{\prime\prime} \in \mathcal{N}^{(m)}(\imath_c^{-1}\tau_n \imath_c(\Delta^{\prime\prime}))$ such that
        \begin{equation*}
            \kappa_c(XN) = N^{\prime\prime} \kappa_c(X)N.
        \end{equation*}
    \end{enumerate}
\end{lemma}
Here, we use the notation developed for affine BD meanders in Section~\ref{s:meander}. Similar equations can be obtained for the case $n \geq m$ by transposition.
\begin{proof}
The proof follows from Proposition~\ref{p:gauss_perturb_matr} and Lemma~\ref{l:taun_w0}. For instance, let us consider the last equation. By~\eqref{eq:taun_w0}, we see that
\begin{equation*}
    \imath_c^{-1}\tau_n\imath_c(\Delta^{\prime\prime}) = \{ \alpha_{n+q} \in A_{m-1} \mid \alpha_q \in \Delta^{\prime\prime}, \ m-q  > n \}.
\end{equation*}
Now, consider a one-parameter subgroup $U_q$ for $\alpha_q\in \Delta^{\prime\prime}$. By Proposition~\ref{p:gauss_perturb_matr}, we see that
\begin{equation*}
    \kappa_c(XU_q(t)) = [W_0^r XU_q(t)]_+ = \begin{cases}
        U_{n+q}(t^\prime)\kappa_c(X) U_q(t) & q < m-n;\\
        \kappa_c(X) U_q(t) & q \geq m-n.
    \end{cases}
\end{equation*}
The condition $q < m-n$ means precisely that $\alpha_{n+q} \in \imath_c^{-1}\tau_n \imath_c(\Delta^{\prime\prime})$. The equation follows.
\end{proof}

We are now ready to address the main result of the subsection:
\begin{proposition}\label{p:birationality}
    If $(\bg^r,\bg^c)$ is aperiodic, then $\mathcal{U}$ is birational. More precisely, there exist rational maps $\delta_r^*$ and $\delta_c$ satisfying equations~\eqref{eq:deltas_r}-\eqref{eq:delta_c} such that $\mathcal{U}^{-1}$ is given by 
    \begin{equation}\label{eq:Ubirat_inv}
        \mathcal{U}^{-1}(Y) = \delta_r^*(Y)^{-1}Y\delta_c(Y)^{-1}
    \end{equation}
    for a generic $Y \in \Mat_{n\times m}(\mathbb{C})$.
\end{proposition}
It is worth abstracting the general principle of the proof, which was used before, for instance, in \cite[Lemma 4.6]{multdual}. Given a set $S$, an element $x_0 \in S$, and a map $F: S \rightarrow S$, consider the sequence $x_j := F(x_{j-1})$, $j \geq 1$ (in other words, a discrete dynamical system). Observe that if the sequence $(x_j)_{j \geq 0}$ forgets the initial point $x_0$, then the sequence stabilizes. More precisely, there exists a number $j_0$ such that for any other sequence $(\tilde{x}_j)_{j \geq 0}$ satisfying $\tilde{x}_j = F(\tilde{x}_{j-1})$, we have $x_{j_0} = \tilde{x}_{j_0}$. Then $x_{j_0} = x_{j_0+1} = x_{j_0+2} = \ldots$; in particular, $x_{j_0} = F(x_{j_0})$. Indeed, one can set $\tilde{x}_0:=x_{k}$ for any $k \geq 0$, and obtain that $x_{j_0} = \tilde{x}_{j_0} = F^{j_0}(x_k) = x_{j_0+k}$.

\begin{proof}
    Our first objective is to apply the above principle to solve the fixed-point equations~\eqref{eq:deltas_r}-\eqref{eq:delta_c} for the pair $(\delta_r^*(Y),\delta_c(Y))$, with the caveat that the map is rational. We will then use the solutions to verify that the right-hand side of~\eqref{eq:Ubirat_inv} is indeed the inverse of $\mathcal{U}$.
    
    Define a rational map 
    \begin{equation*}
        F : \mathcal{N}_-^{(n)}\times \Mat_{n\times m}(\mathbb{C})\times  \mathcal{N}^{(m)} \dashrightarrow \mathcal{N}_-^{(n)}\times  \Mat_{n\times m}(\mathbb{C}) \times \mathcal{N}^{(m)} 
    \end{equation*}
    by
       \begin{equation*}
        F(N_-,Y,N) = \big(\tilde{\gamma}_r^*(N_-\kappa_r(N_-^{-1}YN^{-1})), Y,\tilde{\gamma}_c(\kappa_c(N_-^{-1}YN^{-1})N) \big),
    \end{equation*}
    where the maps $\kappa_r$ and $\kappa_c$ are given by equation~\eqref{eq:kappas}. Consider a sequence $(\delta_{r,j}^*(Y),\delta_{c,j}(Y))_{j \geq 0}$
    such that
    \begin{equation}\label{eq:delta_dynamics}
        (\delta_{r,j}^*(Y),Y,\delta_{c,j}(Y)) = F(\delta_{r,j-1}^*(Y),Y,\delta_{c,j-1}(Y)).
    \end{equation}
    In order to prove that the sequence stabilizes after a finite number of steps, we need to show that it forgets the initial data. However, since the map $F$ is rational, we first need to check that the domain of $F^j$ is nonempty for every $j \geq 2$. Indeed, given $X \in \dom \kappa_r \cap \dom \kappa_c$, observe that $F$ fixes the point $(\rho_r^*(\kappa_r(X)),\mathcal{U}(X),\rho_c(\kappa_c(X)))$. Hence, the domain of each $F^j$ is nonempty and therefore is an open dense subset. 
    
    Now, consider a pair of sequences $({\delta}_{r,j}^*(Y),Y,{\delta}_{c,j}(Y))_{j \geq 0}$ and $(\tilde{\delta}_{r,j}^*(Y),Y,\tilde{\delta}_{c,j}(Y))_{j \geq 0}$ that both satisfy~\eqref{eq:delta_dynamics} and whose initial points belong to $\dom F^j$ for some large enough $j$ (specified below). Note that we choose $Y$ to be the same in both sequences. Define the remainders:
    \begin{equation*}
        A_{r,j}(Y) := \delta_{r,j}^*(Y)^{-1}\tilde{\delta}_{r,j}^*(Y), \quad A_{c,j}(Y) := \tilde{\delta}_{c,j}(Y)\delta_{c,j}(Y)^{-1},
    \end{equation*}
    and define the following sequence of subsets of $A_{n+m-1}$:
    \begin{equation*}
        \Delta_0 := A_{n+m-1}, \ \Delta_j := \gamma \tau_n(\Delta_{j-1}), \quad j \geq 1.
    \end{equation*}
     If we can show that
     \begin{equation}\label{eq:induct}
         A_{r,j}(Y) \in \mathcal{N}_-^{(n)}(\imath_r^{-1}\Delta_{j}), \quad A_{c,j}(Y) \in \mathcal{N}^{(m)}(\imath_c^{-1}\Delta_{j}),
     \end{equation}
     it will follow that $(\delta_{r,j}^*(Y),\delta_{c,j}(Y))_{j \geq 0}$ stabilizes. Indeed, since the BD pair is assumed to be aperiodic, $\Delta_{j} = \emptyset$ for large enough $j$. Therefore, for large enough $j$, $(\delta_{r,j}^*(Y),\delta_{c,j}(Y)) = (\tilde{\delta}_{r,j}^*(Y),\tilde{\delta}_{c,j}(Y))$. Moreover, for the same $j$, we see that $\dom F^j = \dom F^{j+1} = \cdots$, and therefore, the initial data can be taken from $\dom F^j$.

    We prove~\eqref{eq:induct} assuming that $n \leq m$. The proof for $n \geq m$ is similar. By Lemma~\ref{l:kappa_invar}, the map $F$ simplifies to
    \begin{equation*}
        F(N_-,Y,N) = \big( \tilde{\gamma}_r^*(\kappa_r(YN^{-1})),Y, \tilde{\gamma}_c(\kappa_c(N_-^{-1}YN^{-1})N) \big).
    \end{equation*}
    In other words, the contribution of $N_-$ cancels out in the first output argument.
    
    Evidently, \eqref{eq:induct} is satisfied for $j = 0$, for $\imath_r^{-1}\Delta_{0} = A_{n-1}$ and $\imath_c^{-1}\Delta_{0} = A_{m-1}$. Assume that \eqref{eq:induct} holds for $j-1$. We see that
    \begin{equation*}
        \tilde{\delta}_{r,j}^*(Y) = \tilde{\gamma}_r^*(\kappa_r(Y \tilde{\delta}_{c,j-1}(Y)^{-1})) = \tilde{\gamma}_r^*(\kappa_r(Y\delta_{c,j-1}(Y)^{-1} A_{c,j-1}(Y)^{-1})).
    \end{equation*}
    By Lemma~\ref{l:kappa_invar}, since $A_{c,j-1}(Y) \in \mathcal{N}^{(m)}(\imath_c^{-1}\Delta_{j-1})$, there exists $N^\prime_- \in \mathcal{N}^{(n)}_-(\imath_r^{-1}\tau_n\Delta_{j-1})$ such that
    \begin{equation*}
        \kappa_r(Y\delta_{c,j-1}(Y)^{-1} A_{c,j-1}(Y)^{-1}) = \kappa_r(Y\delta_{c,j-1}(Y)^{-1})N^\prime_-;
    \end{equation*}
    hence,
    \begin{equation*}
        A_{r,j}(Y) = \tilde{\gamma}_r^*(N^\prime_-) \in \mathcal{N}^{(n)}_-(\gamma_r^*\imath_r^{-1}\tau_n\Delta_{j-1}) = \mathcal{N}^{(n)}_-(\imath_r^{-1}\Delta_{j}).
    \end{equation*}

    It remains to prove~\eqref{eq:induct} for $A_{c,j}(Y)$. Let us introduce the following auxiliary sequences:
    \begin{equation*}
        X_j = \delta_{r,j}^*(Y)^{-1}Y\delta_{c,j}(Y)^{-1}, \quad \tilde{X}_j = \tilde{\delta}_{r,j}^*(Y)^{-1}Y\tilde{\delta}_{c,j}(Y)^{-1}.
    \end{equation*}
    Note that $\tilde{X}_j = A_{r,j}(Y)^{-1}X_j A_{c,j}(Y)^{-1}$. By Lemma~\ref{l:kappa_invar}, there is $N^\prime \in \mathcal{N}^{(m)}(\imath_c^{-1}\tau_n \Delta_{j-1})$ and $N^{\prime\prime} \in \mathcal{N}^{(m)}(\imath_c^{-1}\tau_n\Delta_{j-1})$ such that
    \begin{equation*}
        \kappa_c(\tilde{X}_{j-1}) = \kappa_c(A_{r,j-1}(Y)^{-1}X_{j-1}A_{c,j-1}(Y)^{-1}) = N^{\prime}N^{\prime\prime} \kappa_c(X_{j-1})A_{c,j-1}(Y)^{-1}.
    \end{equation*}
    Therefore,
    \begin{equation*}\begin{split}
        \tilde{\delta}_{c,j}(Y) &= \tilde{\gamma}_c(\kappa_c(\tilde{X}_{j-1}))\tilde{\gamma}_c(\tilde{\delta}_{c,j-1}(Y)) = \\ &= \tilde{\gamma}_c(N^\prime N^{\prime\prime}) \tilde{\gamma}_c(\kappa_c(X_{j-1}))\tilde{\gamma}_c(A_{c,j-1}(Y))^{-1}\tilde{\gamma}_c(\tilde{\delta}_{c,j-1}(Y)) = \\ &=\tilde{\gamma}_c(N^\prime N^{\prime\prime}) \delta_{c,j}(Y);
        \end{split}
    \end{equation*}
    hence,
    \begin{equation*}
        A_{c,j}(Y) = \tilde{\gamma}_c(N^\prime N^{\prime\prime}) \in \mathcal{N}^{(m)}(\gamma_c\imath_c^{-1}\tau_n(\Delta_{j-1})) = \mathcal{N}^{(m)}(\imath_c^{-1}\Delta_{j}).
    \end{equation*}
    We conclude that there exist rational maps $\delta_r^*$ and $\delta_c$ satisfying equations~\eqref{eq:deltas_r}-\eqref{eq:delta_c}.

    Lastly, let us verify that the rational maps $\delta_r^*$ and $\delta_c$ indeed give rise to the inverse of $\mathcal{U}$. Define
    \begin{equation*}
        \mathcal{V}(Y) := \delta_r^*(Y)^{-1} Y \delta_c(Y)^{-1}.
    \end{equation*}
    Then the fixed-point equations for $(\delta_r^*(Y),\delta_c(Y))$ can be expressed as
    \begin{equation*}
        \delta_r^*(Y) = \tilde{\gamma}_r^*(\delta_r^*(Y)\kappa_r(\mathcal{V}(Y))), \quad \delta_c(Y) = \tilde{\gamma}_c(\kappa_c(\mathcal{V}(Y))\delta_c(Y)).
    \end{equation*}
    Iterating, we see that $(\delta_r^*(Y),\delta_c(Y))$ satisfy the equations
    \begin{equation*}
        \delta_r^*(Y) = \rho_r^*(\kappa_r(\mathcal{V}(Y))), \quad \delta_c(Y) = \rho_c(\kappa_c(\mathcal{V}(Y)));
    \end{equation*}
    hence,
    \begin{equation*}
        \mathcal{U}(\mathcal{V}(Y)) = \rho_r^*(\kappa_r(\mathcal{V}(Y))) \delta_r^*(Y)^{-1}Y \delta_c(Y)^{-1} \rho_c(\kappa_c(\mathcal{V}(Y))) = Y.
    \end{equation*}
    Since these maps are from an irreducible algebraic variety to itself, we conclude that $\mathcal{V}$ is a two-sided inverse of $\mathcal{U}$.
\end{proof}

\begin{remark}
    The proof of Proposition~\ref{p:birationality} can be easily upgraded to the setting of a complex reductive algebraic group $G$. In fact, the proof is even easier, since it follows the same lines as the case $n=m$ for matrices. In the setting of $G$, given a generic $y \in G$, the map $F$ is given by
    \begin{equation*}
        F(u_-,y,u) = \big( \tilde{\gamma}_r^*([yu^{-1}\overline{w}_0]_-), y,\tilde{\gamma}_c([\overline{w}_0u_-^{-1}y]_+)\big), \quad (u_-,y,u) \in \mathcal{N}_-\times G\times  \mathcal{N}.
    \end{equation*}
\end{remark}

\begin{remark}
    Let us comment on the existing results on the map $\mathcal{U}$. In~\cite{gsv-unified}, given a reductive complex algebraic group $G$, the authors constructed a Poisson rational map
    \begin{equation*}
        \mathcal{U} : (G,\pi_{(\bbg_{\std}^r,\bbg_{\std}^c)}) \dashrightarrow (G, \pi_{(\bbg^r,\bbg^c)}).
    \end{equation*}
    Their construction uses intermediate factorizations in the subgroups corresponding to the BD data. We found that such factorizations can be avoided, and one can write a single concise formula. In a new version of~\cite{gsv-unified}, the authors will provide a proof that our formula coincides with theirs. Birationality was shown only for square matrices, by explicitly constructing a cluster in $\GL_n(\mathbb{C})$ associated with $(\bg^r,\bg^c)$, and then evaluating the map on that cluster. We derived a formula for the inverse via a recursive procedure, which works not only for rectangular matrices but for arbitrary $G$ as well. The converse remains an open question: we do not know if birationality of $\mathcal{U}$ implies that $(\bg^r,\bg^c)$ is aperiodic.
\end{remark}

\subsection{Poisson properties}\label{s:biratp}
In the final subsection, we prove that the rational map $\mathcal{U}$ given by~\eqref{eq:Ubirat} is Poisson. Combined with Proposition~\ref{p:U0}, this will conclude the proof of the second part of Theorem~\ref{theorem}. We then use this result to prove that $(\Mat_{n\times m}(\mathbb{C}),\pi_{(\bbg^r,\bbg^c)})$ has an open $\mathbb{T}_{(\bg^r,\bg^c)}$-leaf for aperiodic $(\bg^r,\bg^c)$.

Unless stated otherwise, we use the following setup throughout the subsection. We denote by $\bbg^r:=(\bg^r,r_0^r)$ and $\bbg^c:=(\bg^c,r_0^c)$ arbitrary BD quadruples for~$\gl_n(\mathbb{C})$ and $\gl_m(\mathbb{C})$, respectively. In addition, we set $\bbg^r_{\std}:=(\bg^r_{\std},r_0^r)$ and $\bbg^c_{\std}:=(\bg^c_{\std},r_0^c)$. We denote by $r_{\std}^r$ and $r_{\std}^c$ the $r$-matrices corresponding to $\bg_{\std}^r$ and $\bg_{\std}^c$, and by $R_{\std}^r := \pr_> + R_0^r$ and $R_{\std}^c := \pr_> + R_0^c$ the corresponding linear maps. For  $k \geq 1$, we denote by $\bar\bg_{\std}^{(k-1)} = (\bg_{\std}^{(k-1)}, \frac{1}{2}\Omega_0^{(k)})$ the standard BD quadruple for~$\gl_k(\mathbb{C})$. The corresponding $r$-matrix is denoted by $r_{\std}^{(k)}$, and the corresponding linear map is denoted by $R_{\std}^{(k)}:=\pr_> + \frac{1}{2}\pr_0$.

Our strategy for proving that $\mathcal{U}$ is Poisson is as follows. First, in Lemma~\ref{l:nmPoiss}, we observe that the case $n \leq m$ is equivalent to the case $n \geq m$. This will allow us to focus on $n \leq m$. Second, we define a pair of maps $\mathcal{U}_r$ and $\mathcal{U}_c$ as the maps~\eqref{eq:Ubirat} relative to the pairs $(\bg^r,\bg_{\std}^{(m-1)})$ and $(\bg_{\std}^{(n-1)},\bg^c)$. That is,
\begin{equation*}
    \mathcal{U}_r(X) = \rho_r^*(\kappa_r(X))X, \quad \mathcal{U}_c(X) = X \rho_c(\kappa_c(X)).
\end{equation*}
Recall from~\eqref{eq:kappas} that $\kappa_r(X) = [XW_0^c]_-$ and $\kappa_c(X) = [W_0^r X]_+$. In Proposition~\ref{p:Ufac}, we observe that $\mathcal{U}$ is Poisson if and only if the maps $\mathcal{U}_r$ and $\mathcal{U}_c$ are Poisson. In other words, instead of proving that $\mathcal{U}$ is Poisson for two nontrivial BD triples, we can let one of the triples be standard, which simplifies the argument. This idea was previously used in~\cite{gsv-unified} in the setting of reductive complex algebraic groups. We will show that both maps $\mathcal{U}_r$ and $\mathcal{U}_c$ are Poisson using the result established in~\cite{gsv-unified} for square matrices, together with a series of auxiliary lemmas.

For a BD quadruple $\bbg := (\Gamma_1,\Gamma_2,\gamma,r_0)$, let us set $\bbg^{\op}=(\Gamma_2,\Gamma_1,\gamma^*,r_0^t)$. We begin with the following lemma:
\begin{lemma}\label{l:nmPoiss}
    The transposition map is a Poisson isomorphism as the map
    \begin{equation*}
        (\Mat_{n\times m}(\mathbb{C}),\pi_{(\bbg^r,\bbg^c)}) \rightarrow (\Mat_{m\times n}(\mathbb{C}),\pi_{((\bbg^c)^{\op},(\bbg^r)^{\op})}), \quad A \mapsto A^T.
    \end{equation*}
    As a consequence, the map $\mathcal{U}$ is Poisson for $n \leq m$ if and only if it is Poisson for $n \geq m$.
\end{lemma}
\begin{proof}
    Denote by $T(A) := A^T$ the transposition map. For an $r_{\bbg}$-matrix given by $\bbg$ (see equation~\eqref{eq:r_bd}), it is easy to see that $[T_*\otimes T_*](r_{\bbg}) = (r_{\bbg^{\op}})^t$. Hence,
    \begin{align*}
        &[T_*\otimes T_*][(L_X)_*\otimes (L_X)_*](r_{\bbg^c}) = [(R_{X^T})_* \otimes (R_{X^T})_*]((r_{(\bbg^c)^{\op}})^t);\\
        &[T_*\otimes T_*][(R_X)_*\otimes (R_X)_*](r_{\bbg^r}) = [(L_{X^T})_*\otimes (L_{X^T})_*]((r_{(\bbg^r)^{\op}})^t).
    \end{align*}
    Therefore, $[T_*\otimes T_*](\pi_{(\bbg^r,\bbg^c)}) = \pi_{((\bbg^c)^{\op},(\bbg^r)^{\op})}$. For the second statement, observe that the following diagram commutes:
    \begin{equation*}
        \xymatrix{
            (\Mat_{n\times m}(\mathbb{C}),\pi_{(\bbg^r_{\std},\bbg^c_{\std})})\ar[d]_T \ar@{-->}[r] &(\Mat_{n\times m}(\mathbb{C}),\pi_{(\bbg^r,\bbg^c)})\ar[d]_{T}\\
            (\Mat_{m\times n}(\mathbb{C}),\pi_{((\bbg^c_{\std})^{\op},(\bbg^r_{\std})^{\op})}) \ar@{-->}[r] &(\Mat_{m\times n}(\mathbb{C}),\pi_{((\bbg^c)^{\op},(\bbg^r)^{\op})})
        }
    \end{equation*}
    where the horizontal arrows are the maps $\mathcal{U}$ for the indicated choices of BD quadruples. Thus the statement holds.
\end{proof}

Due to Lemma~\ref{l:nmPoiss}, for the rest of the subsection, we focus on the case $n \leq m$.

For $n \leq m$, define
\begin{equation*}
    I_{n,m} := \begin{bmatrix} 
        0 & I_n
    \end{bmatrix} \in \Mat_{n\times m}(\mathbb{Z}), \quad I_{m,n}:=(I_{n,m})^T.
\end{equation*} 
\begin{align*}
    &\Mat_{n\times m}^+(\mathbb{C}):= \mathcal{B}^{(n)}I_{n,m} = \{ \begin{bmatrix}
        0 & B
    \end{bmatrix} \in \Mat_{n\times m}(\mathbb{C}) \mid B \in \mathcal{B}^{(n)} \};\\
    &\Mat_{n\times m}^-(\mathbb{C}):= I_{n,m}\mathcal{B}^{(m)}_- = \{ \begin{bmatrix}
        A & B_-
    \end{bmatrix} \in \Mat_{n\times m}(\mathbb{C}) \mid A \in \Mat_{n\times (m-n) }(\mathbb{C}), \ B_- \in \mathcal{B}_-^{(n)}  \}.
\end{align*}

\begin{lemma}\label{lm:MpmPoi}
    For $n \leq m$, 
    the subvarieties $\Mat_{n\times m}^+(\mathbb{C})$ and $\Mat_{n\times m}^-(\mathbb{C})$ of $(\Mat_{n\times m}(\mathbb{C}),\pi_{(\bbg^r_{\std},\bbg^c_{\std})})$ are Poisson. 
\end{lemma}
\begin{proof}
    We need to show that for any $f_1,f_2 \in \mathbb{C}[\Mat_{n\times m}]$ such that $f_2|_{\Mat_{n\times m}^+(\mathbb{C})} = 0$, the identity $\{f_1,f_2\}|_{\Mat_{n\times m}^+(\mathbb{C})} = 0$ holds, where the Poisson bracket corresponds to $\pi_{(\bbg^r_{\std},\bbg^c_{\std})}$ (likewise for $\Mat_{n\times m}^-(\mathbb{C})$).

    For $X \in \Mat_{n\times m}^+(\mathbb{C})$, we see that  $(\nabla_X f_2)_{[m-n+1,m]} \in \mathfrak{n}^{(n)}$; hence,
    \begin{equation*}
        (\nabla^L_Xf_2)^{[1,m-n]} = 0, \quad (\nabla^L_Xf_2)_{[m-n+1,m]}^{[m-n+1,m]}, \ \nabla^R_X f_2 \in \mathfrak{n}^{(n)}.
    \end{equation*}
    Therefore,
    \begin{equation*}
        (R^c_{\std})^*(\nabla_X^L f_2) = 0, \quad (R^r_{\std})^*(\nabla_X^R f_2) = 0,
    \end{equation*}
    and hence,
    \begin{equation*}
        \langle R^c_{\std}(\nabla^L_X f_1),\nabla^L_X f_2\rangle = \langle R^r_{\std}(\nabla^R_X f_1),\nabla^R_X f_2\rangle = 0.
    \end{equation*}
    We conclude that the subvariety $\Mat_{n\times m}^+(\mathbb{C})$ is Poisson.

    Similarly, if $f_2 \in \mathbb{C}[\Mat_{n\times m}]$ is such that $f_2|_{\Mat_{n\times m}^-(\mathbb{C})} = 0$, then $(\nabla_X f_2)_{[m-n+1,m]} \in \mathfrak{n}_-^{(n)}$ and $(\nabla_X f_2)_{[1,m-n]} = 0$ for $X \in \Mat_{n\times m}^-(\mathbb{C})$; hence,
    \begin{equation*}
        (\nabla^L_Xf_2)_{[1,m-n]} = 0, \qquad (\nabla^L_Xf_2)_{[m-n+1,m]}^{[m-n+1,m]}, \ \nabla^R_X f_2 \in \mathfrak{n}_-^{(n)};
    \end{equation*}
    therefore,
    \begin{equation*}
        R_{\std}^c(\nabla_X^L f_2) = 0, \quad R_{\std}^r(\nabla_X^R f_2) = 0.
    \end{equation*}
    It follows that $\Mat_{n\times m}^-(\mathbb{C})$ is a Poisson subvariety.
\end{proof}

\begin{lemma}\label{lm:InmPoi}
    Let $k$ and $l$ be any positive integers. Then the map
    \begin{equation}\label{eq:InmPoi}
        (\Mat_{n\times k}(\mathbb{C}),\pi_{(\bbg^r,\bbg_{\std}^{(k-1)})}) \times (\Mat_{l\times m}(\mathbb{C}), \pi_{(\bbg_{\std}^{(l-1)},\bbg^c)}) \rightarrow (\Mat_{n\times m}(\mathbb{C}),\pi_{(\bbg^r,\bbg^c)})
    \end{equation}
    given by
    \begin{equation*}
        \Mat_{n\times k}(\mathbb{C}) \times \Mat_{l\times m}(\mathbb{C}) \ni (X,Y) \mapsto X I_{k,l} Y \in \Mat_{n\times m}(\mathbb{C})
    \end{equation*}
    is Poisson.
\end{lemma}
\begin{proof}
    This is a standard argument. Denoting the map by $\mu$, we see that
    \begin{equation*}\begin{split}
        \mu_*(\pi_{(\bbg^r,\bbg_{\std}^{(k-1)})}&\times \pi_{(\bbg_{\std}^{(l-1)},\bbg^c)})_{(X,Y)} = \\&=  (R_{I_{k,l}Y})_*(\pi_{(\bbg^r,\bbg_{\std}^{(k-1)})})_X + (L_{XI_{k,l}})_*(\pi_{(\bbg_{\std}^{(l-1)},\bbg^c)})_Y =  \\ &= (R_{I_{k,l}Y})_*(L_X)_*r_{\std}^{(k)} - (R_{XI_{k,l}Y})_*r^r + (L_{XI_{k,l}Y})_*r^c - (L_{XI_{k,l}})_*(R_Y)_*r_{\std}^{(l)}.
        \end{split}
    \end{equation*}
    By explicit computation, we see that $(R_{I_{k,l}})_*r_{\std}^{(k)} = (L_{I_{k,l}})_*r_{\std}^{(l)}$. Thus the statement holds.
\end{proof}

\begin{proposition}\label{p:Ufac}
    Assume that $n \leq m$. Then the following diagram commutes:
   \begin{equation}\label{eq:Ufac}
        \scalebox{0.86}{$\xymatrix{
            (\Mat_{n\times m}^+,\pi_{(\bbg^r_{\std},\bbg_{\std}^{(m-1)})})\times (\Mat_{n\times m}^-,\pi_{(\bbg_{\std}^{(n-1)},\bbg^c_{\std})}) \ar[d] \ar[r]^-{(\mathcal{U}_r,\mathcal{U}_c)} & (\Mat_{n\times m},\pi_{(\bbg^r,\bbg_{\std}^{(m-1)})})\times (\Mat_{n\times m},\pi_{(\bbg_{\std}^{(n-1)},\bbg^c)}) \ar[d]\\
            (\Mat_{n\times m},\pi_{(\bbg^r_{\std},\bbg^c_{\std})}) \ar@{-->}[r]^{\mathcal{U}} & (\Mat_{n\times m},\pi_{(\bbg^r,\bbg^c)})
        }$}
    \end{equation}
    where the vertical arrows are given by $(X,Y) \mapsto XI_{m,n}Y$. In particular, $\mathcal{U}$ is Poisson if and only if both $\mathcal{U}_r|_{\Mat_{n\times m}^+}$ and $\mathcal{U}_c|_{\Mat_{n\times m}^-}$ are Poisson.
\end{proposition}
\begin{proof}
    For any $B \in \mathcal{B}^{(n)}$ and any $B_- \in \mathcal{B}_-^{(m)}$, observe that the maps $\kappa_r$ and $\kappa_c$ satisfy the following equations:
    \begin{equation*}
        \kappa_r(BI_{n,m}B_-) = \kappa_r(BI_{n,m}), \quad \kappa_c(BI_{n,m}B_-) = \kappa_c(I_{n,m}B_-).
    \end{equation*}
    Since $I_{n,m}I_{m,n}I_{n,m} = I_{n,m}$, observe that the left vertical map sends a pair $(BI_{n,m}, I_{n,m}B_-)$ to the element $BI_{n,m}B_-$. Therefore,
    \begin{equation*}\begin{split}
        \mathcal{U}(BI_{n,m}B_-) &= \rho_r^*(\kappa_r(BI_{n,m}))BI_{n,m}B_-\rho_c(\kappa_c(I_{n,m}B_-)) = \\ &= \rho_r^*(\kappa_r(BI_{n,m}))BI_{n,m}\cdot I_{m,n}\cdot I_{n,m}B_-\rho_c(\kappa_c(I_{n,m}B_-)) = \mathcal{U}_r(BI_{n,m})I_{m,n}\mathcal{U}_c(I_{n,m}B_-);
        \end{split}
    \end{equation*}
    hence, the diagram is commutative. Since the vertical maps are dominant and Poisson by Lemma~\ref{lm:MpmPoi} and Lemma~\ref{lm:InmPoi}, the map $\mathcal{U}$ is Poisson if and only if both $\mathcal{U}_r|_{\Mat_{n\times m}^+}$ and $\mathcal{U}_c|_{\Mat_{n\times m}^-}$ are Poisson.
\end{proof}

For the next lemma, given a subset $S$ of a complex affine Poisson variety $(V,\pi_V)$, we denote by $I_S$ its vanishing ideal, and by $N(I_S)$ the Poisson normalizer of $I_S$; that is,
\begin{equation*}
    N(I_S):=\{f \in \mathbb{C}[V] \mid \{f,I_{S}\}\subseteq I_{S}\}.
\end{equation*}

\begin{lemma}\label{l:coisop_slice}
    Let $(V,\pi_V)$, $(V^\prime,\pi_{V^\prime})$, $(B,\pi_B)$, and $(B^\prime,\pi_{B^\prime})$ be complex affine Poisson varieties, and let $S$ be a coisotropic subvariety of $(V,\pi_V)$. Consider the following diagrams:
    \begin{equation}\label{eq:coisop_slice}
    \xymatrix{
        V \ar@{-->}[r]^{F} \ar[d]_{q} & V^\prime \ar[d]^{q^\prime} \\
        B \ar@{-->}[r]_{\bar{F}} & B^\prime
        } \quad \quad \xymatrix{
        S \ar@{-->}[r]^{F|_S} \ar[d]_{q|_S} & V^\prime \ar[d]^{q^\prime} \\
        Z \ar@{-->}[r]_{\bar{F}|_Z} & B^\prime
        }
    \end{equation}
    Assume that the maps $q$ and $q^\prime$ are regular and Poisson, the maps $F$ and $\bar{F}$ are rational, and the map $F$ is Poisson and well-defined on an open dense subset of $S$. In addition, assume that $q|_S$ is dominant onto a Poisson subvariety $Z$ of $B$ such that $q^*(\mathbb{C}[Z])\subseteq N(I_S)$, and the right diagram commutes on a common domain (but not necessarily the left one). Then the map $\bar{F}|_{Z}$  is Poisson.
\end{lemma}
In the statement of the lemma, we do not need to assume that $S$ is smooth. If, however, it is the case, then the condition $q^*(\mathbb{C}[Z]) \subseteq N(I_S)$ is equivalent to the vanishing of $dq$ along the characteristic distribution of $S$. In the setting of the subsection, we  take $F$ to be the map~\eqref{eq:Ubirat} defined on $m \times m$ matrices with respect to $(\bg_{\std}^{(m-1)},\bg^c)$, and $\bar{F}|_Z$ to be the map $\mathcal{U}_c|_{\Mat^-_{n \times m}}$. The lemma is a generalization of \cite[Proposition 6.7, Part 1]{lu-evens} and therefore can also be verified using similar technics at the linear level, assuming everything is smooth.
\begin{proof}
    Since $S$ is coisotropic, $I_S$ is a Poisson subalgebra of $\mathbb{C}[V]$, and therefore, $I_S \subseteq N(I_S)$. Moreover, $N(I_S)/I_S$ is a Poisson algebra. Since $q^*(\mathbb{C}[Z])\subseteq N(I_S) \subseteq \mathbb{C}[V]$, the map $q^*$ induces a map $(q|_S)^* : \mathbb{C}[Z] \rightarrow N(I_S)/I_S\subseteq \mathbb{C}[S]$. Moreover, since $q|_S$ is dominant onto $Z$, the map $(q|_S)^*$ is injective.

    Let us denote by $[f] \in N(I_S)/I_S$ the class of an element $f \in N(I_S)$. Under appropriate localizations, since the right diagram commutes, for any $h \in \mathbb{C}[B^\prime]$,
    \begin{equation*}
     [q^*(\bar{F}|_Z)^*(h)]= (q|_S)^*(\bar{F}|_Z)^*(h) = (F|_S)^*(q^\prime)^*(h) = [F^*(q^\prime)^*(h)].
    \end{equation*}
    Since $q$, $F$, and $q^\prime$ are Poisson, we see that for any $h_1,h_2 \in \mathbb{C}[B^\prime]$,
    \begin{equation*}\begin{split}
        [q^*(\bar{F}|_Z)^*\{h_1,h_2\}_{B^\prime}] &= [F^*(q^\prime)^*\{h_1,h_2\}_{B^\prime}] = [\{F^*(q^\prime)^*(h_1),F^*(q^\prime)^*(h_2)\}_{V}] =\\&= \{[F^*(q^\prime)^*(h_1)],[F^*(q^\prime)^*(h_2)]\}_{N(I_S)/I_S} = \\ &=\{[q^*(\bar{F}|_Z)^*(h_1)],[q^*(\bar{F}|_Z)^*(h_2)]\}_{N(I_S)/I_S} = [q^*\{(\bar{F}|_Z)^*(h_1),(\bar{F}|_Z)^*(h_2)\}_{B}].
        \end{split}
    \end{equation*}
    Since $(q|_S)^*$ is injective, we conclude that $\bar{F}|_Z$ is Poisson.
\end{proof}

\begin{lemma}\label{l:wbm}
    Let $G$ be a reductive complex algebraic group equipped with a Poisson bivector $\pi_{(\bbg_{\std}^r,\bbg_{\std}^c)}$ for $\bbg_{\std}^r:=(\bg_{\std},r_0^r)$ and $\bbg_{\std}^c:=(\bg_{\std},r_0^c)$. Then for any $w \in W$, $w\mathcal{B}_-$ is a coisotropic submanifold of $(G,\pi_{(\bbg_{\std}^r,\bbg_{\std}^c)})$.
\end{lemma}
\begin{proof}
    Observe that if $f \in \mathbb{C}[G]$ vanishes at $w\mathcal{B}_-$, then for any $x \in w\mathcal{B}_-$, 
    \begin{equation*}
        \nabla^L_xf \in \mathfrak{n}_-, \quad \nabla^R_x f \in \Ad_w \mathfrak{n}_-.
    \end{equation*}
    Therefore, for $x \in w\mathcal{B}_-$, $R^c_{\std}(\nabla^L_x f) = 0$, and if $f_1,f_2 \in \mathbb{C}[G]$ vanish at $w\mathcal{B}_-$, then 
    \begin{equation*}
        \langle R^r_{\std}(\nabla^R_x f_1),\nabla^R_x f_2\rangle = 0,
    \end{equation*}
    for $\pr_{>}\nabla^R_xf_1 \in \mathfrak{n} \cap \Ad_w \mathfrak{n}_-$ and $\nabla^R_x f_2 \in \Ad_w \mathfrak{n}_-$, and these two are orthogonal to each other.
\end{proof}

\begin{proposition}\label{p:u_poiss}
    For any $(n,m)$ and any pair of BD quadruples $(\bbg^r,\bbg^c)$, the map $\mathcal{U}$ given by~\eqref{eq:Ubirat} is Poisson.
\end{proposition}
\begin{proof}
    By Lemma~\ref{l:nmPoiss}, it suffices to consider the case $n \leq m$, and by Proposition~\ref{p:Ufac}, it suffices to show that the maps $\mathcal{U}_r|_{\Mat_{n\times m}^+}$ and $\mathcal{U}_c|_{\Mat_{n\times m}^-}$ are Poisson. We will use the corresponding result for $n=m$ proved in~\cite{gsv-unified}.

    To show that $\mathcal{U}_r|_{\Mat_{n\times m}^+}$ is Poisson, consider the Poisson map 
    \begin{equation*}
        \mathcal{U}_r^\prime :(\Mat_{n\times n},\pi_{(\bbg^r_{\std},\bbg_{\std}^{(n-1)})}) \dashrightarrow (\Mat_{n\times n},\pi_{(\bbg^r,\bbg_{\std}^{(n-1)})})
    \end{equation*}
    given by~\eqref{eq:Ubirat} with respect to the pair $(\bg^r,\bg_{\std}^{(n-1)})$.
    \begin{align*}
        &p:(\Mat_{n\times m},\pi_{(\bbg^r_{\std},\bbg_{\std}^{(m-1)})}) \rightarrow (\Mat_{n\times n}, \pi_{(\bbg^r_{\std},\bbg_{\std}^{(n-1)})}), \quad p(X) = XI_{m,n} = X^{[m-n+1,m]};\\
        &\imath : (\Mat_{n\times n},\pi_{(\bbg^r,\bbg_{\std}^{(n-1)})}) \rightarrow (\Mat_{n\times m},\pi_{(\bbg^r,\bbg_{\std}^{(m-1)})}), \quad \imath(X) = XI_{n,m} = \begin{bmatrix}
            0 & X
        \end{bmatrix}.
    \end{align*}
    Both maps are Poisson by choosing appropriate dimensions and restricting the map~\eqref{eq:InmPoi} in Lemma~\ref{lm:InmPoi}. Now observe that the map $\mathcal{U}_r|_{\Mat_{n\times m}^+}$ can be written as the composition
    \begin{equation*}
        \mathcal{U}_r|_{\Mat_{n\times m}^+} = \imath \circ \mathcal{U}_r^\prime \circ p|_{\Mat_{n\times m}^+} : (\Mat_{n\times m}^+,\pi_{(\bbg^r_{\std},\bbg_{\std}^{(m-1)})}) \rightarrow (\Mat_{n \times m},\pi_{(\bbg^r,\bbg_{\std}^{(m-1)})}).
    \end{equation*}
    Therefore, the map $\mathcal{U}_r|_{\Mat_{n\times m}^+}$ is Poisson.

    To prove that $\mathcal{U}_c$ is Poisson, we will use Lemma~\ref{l:coisop_slice}. Let $W_0^{(m-n)}$ be the $(m-n)\times (m-n)$ matrix corresponding to the longest permutation (with all nonzero entries equal to $1$), and let $w$ be the $m \times m$ permutation matrix given by
    \begin{equation*}
        w:= \begin{bmatrix}
            W_0^{(m-n)} &\\
            & I_{n}
        \end{bmatrix}.
    \end{equation*}
    By Lemma~\ref{l:wbm}, $S:=w\mathcal{B}_-^{(m)}$ is a coisotropic submanifold of $(\Mat_{m \times m},\pi_{(\bbg^{(m-1)}_{\std},\bbg^c_{\std})})$. Denote $Z:=\Mat_{n\times m}^-$. We consider the following pair of diagrams corresponding to~\eqref{eq:coisop_slice}:
 \begin{equation}\label{eq:coisop_slice_uc}
    \xymatrix{
        (\Mat_{m\times m},\pi_{(\bbg_{\std}^{(m-1)},\bbg^c_{\std})}) \ar@{-->}[r]^{\mathcal{U}_c^\prime} \ar[d]_{P} & (\Mat_{m\times m},\pi_{(\bbg_{\std}^{(m-1)},\bbg^c)}) \ar[d]^{P} \\
        (\Mat_{n\times m},\pi_{(\bbg_{\std}^{(n-1)},\bbg^c_{\std})}) \ar@{-->}[r]_{\mathcal{U}_c} & (\Mat_{n\times m},\pi_{(\bbg_{\std}^{(n-1)},\bbg^c)}) 
        } \quad \quad 
    \xymatrix{
        S = w\mathcal{B}_-^{(m)} \ar@{-->}[r]^-{\mathcal{U}_c^\prime|_S} \ar[d]_{P|_S} & \Mat_{m\times m} \ar[d]^{P} \\
        Z = \Mat_{n\times m}^- \ar@{-->}[r]_-{\mathcal{U}_c|_{Z}} & \Mat_{n\times m}
        }
    \end{equation}
    Here, the map $P$ is given by $P(X) := I_{n,m}X$. It is easy to see that $P$ is dominant. Moreover, $P$ is Poisson by Lemma~\ref{lm:InmPoi}.

    Let us verify that $P^*(\mathbb{C}[Z]) \subseteq N(I_S)$, where the normalizer is considered relative to $\pi_{(\bbg_{\std}^{(m-1)},\bbg^c_{\std})}$. Let $f \in \mathbb{C}[Z]$, $g \in I_S$, and $B_- \in \mathcal{B}_-^{(m)}$. As in the proof of Lemma~\ref{l:wbm}, we see that
    \begin{equation*}
        \nabla^L_{wB_-}g \in \mathfrak{n}_-, \qquad \nabla^R_{wB_-}g \in \Ad_w\mathfrak{n}_-;
    \end{equation*}
    hence, 
    \begin{equation*}
        R^{c}_{\std}(\nabla^L_{wB_-}g) = 0,
    \end{equation*}
    and therefore,
    \begin{equation*}
        \langle (R^{c}_{\std})_-(\nabla^L_{wB_-}(P^*f)),\nabla^L_{wB_-}g \rangle =\langle \nabla^L_{wB_-}(P^*f), -R^{c}_{\std}(\nabla^L_{wB_-}g)\rangle =  0.
    \end{equation*}
    For the other term in the Poisson bracket, since 
    \begin{equation*}
        \nabla_{wB_-}(P^*f) = \nabla_{P(wB_-)}f \cdot I_{n,m},
    \end{equation*}
    we see that the first $m-n$ columns of $\nabla_{wB_-}(P^*f)$ (hence of $\nabla^R_{wB_-}(P^*f)$) are zero. Therefore, 
    \begin{equation*}
        \langle (R^{(m)}_{\std})_-(\nabla^R_{wB_-}(P^*f)), \nabla^R_{wB_-}g \rangle = 0.
    \end{equation*}
    We conclude that $P^*(\mathbb{C}[Z]) \subseteq N(I_S)$.

    Observe that $\mathcal{U}_c^\prime$ is defined on an open dense subset of $S$. Indeed, the domain of $\mathcal{U}_c^\prime$ includes the nonzero locus of the minors
    \begin{equation*}
        \det X^{[1,i]}_{[m-i+1,m]}, \ i \in [1,m].
    \end{equation*}
    It is easy to see that the nonzero locus of the above minors intersects $w\mathcal{B}_-^{(m)}$ in an open dense subset.

    Lastly, let us verify that the right diagram in~\eqref{eq:coisop_slice_uc} commutes. It suffices to show that
    \begin{equation*}
        [W_0^{(m)} w B_-]_+ = [W_0^{(n)}I_{n,m}wB_-]_+.
    \end{equation*}
    The equality follows from the following matrix identity:
    \begin{equation*}
        W_0^{(m)}wB_- = \begin{bmatrix}
            I_n & \\
            & (B_-)_{[1,m-n]}^{[1,m-n]}
        \end{bmatrix}
        \begin{bmatrix}
            W_0^{(n)} (B_-)_{[m-n+1,m]}^{[1,m-n]} & W_0^{(n)} (B_-)_{[m-n+1,m]}^{[m-n+1,m]}\\
            I_{m-n} & 0
        \end{bmatrix}.
    \end{equation*}
    Thus, by Lemma~\ref{l:coisop_slice}, $\mathcal{U}_c|_Z$ is Poisson.
\end{proof}

In the final part of the subsection, we address the existence of an open $\mathbb{T}_{(\bg^r,\bg^c)}$-leaf in $\Mat_{n\times m}(\mathbb{C})$ for aperiodic $(\bg^r,\bg^c)$. The result is a corollary of a more general statement.

\begin{proposition}\label{p:open_t_leaf}
    Let $(V,\pi_V)$ and $(V^\prime,\pi_{V^\prime})$ be irreducible complex Poisson varieties. Assume there is an action of a complex algebraic torus $\mathbb{T}$ on both $V$ and $V^\prime$ by Poisson automorphisms, and assume there is a $\mathbb{T}$-equivariant Poisson birational map $F : (V,\pi_V) \dashrightarrow (V^\prime,\pi_{V^\prime})$. Then $(V,\pi_V)$ has an open $\mathbb{T}$-leaf if and only if $(V^\prime,\pi_{V^\prime})$ has an open $\mathbb{T}$-leaf.
\end{proposition}

\begin{proof}
    Let us first show that the maximal domain of definition $\dom F$ of $F$ is a $\mathbb{T}$-invariant Zariski open subset of $V$. Indeed, since $\mathbb{T}$ acts by biregular isomorphisms, if $W \subseteq \dom F$ is an open subset, then $F$ is defined on $t\cdot W$ by $F(x) = t\cdot F(t^{-1}\cdot x)$, $x \in t\cdot W$. 
    
    Let us set $U:=\dom F \cap F^{-1}(\dom F^{-1})$ and $U^\prime := F(U)$. Then $F$ restricts to a $\mathbb{T}$-equivariant Poisson biregular isomorphism 
    \begin{equation*}
        F|_U:(U,\pi_U) \rightarrow (U^\prime,\pi_{U^\prime})
    \end{equation*}
    between $\mathbb{T}$-invariant open subvarieties $U$ of $V$ and $U^\prime$ of $V^\prime$.
    
    To find an open $\mathbb{T}$-leaf of $(V^\prime,\pi_{V^\prime})$, we need the notion of \emph{$\mathbb{T}$-Pfaffians}. The notion was introduced for $\mathbb{T}$-leaves in~\cite{lm-leaf}. The definition can be generalized as follows. Let $(\mathsf{Y}, \pi_{\mathsf{Y}})$ be an irreducible $\mathbb{T}$-Poisson variety, $\mathsf{Y}^{\sm}$ be its smooth locus, and  $d:=\dim \mathsf{Y}$. Let $2r$ be the rank of $\pi_{\mathsf{Y}}$ (that is, the maximal dimension of symplectic leaves of $\pi_\mathsf{Y}$ in $\mathsf{Y}^{\sm}$ is $2r$). A $\mathbb{T}$-Pfaffian is a nonzero anticanonical section on $\mathsf{Y}^{\sm}$ of the form
    \[
        \mathrm{Pf}_\mathbb{T}(\pi_{\mathsf{Y}})=(\wedge^{r} \pi_{\mathsf{Y}}) \wedge \chi({\xi_1}) \wedge \cdots \wedge \chi(\xi_{d-2r}),
    \]
    where $\chi : \mathfrak{t} \rightarrow \mathfrak{X}^1(\mathsf{Y})$ is the infinitesimal $\mathbb{T}$-action and $\xi_1,\cdots, \xi_{d-2r}$ are linearly independent elements of $\mathfrak{t}$. Assuming that a $\mathbb{T}$-Pfaffian on $\mathsf{Y}$ exists, let us define
    \begin{equation*}
        W:=\{x \in \mathsf{Y}^{\sm} \mid \Pf_{\mathbb{T}}(\pi_{\mathsf{Y}})_x \neq 0\}.
    \end{equation*}
    We show that $W$ is precisely the open $\mathbb{T}$-leaf of $(\mathsf{Y}, \pi_\mathsf{Y})$. Evidently, $W$ is Zariski open, so we only need to show that $W$ is a single $\mathbb{T}$-leaf. Since $\Pf_{\mathbb{T}}(\pi_{\mathsf{Y}})$ is invariant under the $\mathbb{T}$-action and Hamiltonian vector fields~\cite[Proposition 3.15]{py-fs}, we see that $W$ is invariant under local Hamiltonian flows and the $\mathbb{T}$-action. Therefore, $W$ is a union of its $\mathbb{T}$-leaves.
    Let $L$ be any of these $\mathbb{T}$-leaves, and let $x \in L$. Since $\Pf_{\mathbb{T}}(\pi_{\mathsf{Y}})_x\neq 0$, the inverse function theorem implies that there
    exists a complex analytic open neighborhood of $x$ that can be reached by finite compositions of local Hamiltonian flows and the $\mathbb{T}$-action. Therefore, $L$ is open in $W$ in the complex analytic topology, and the $\mathbb{T}$-leaves partition $W$ into disjoint subsets open in the complex analytic topology. On the other hand, since $\mathsf{Y}$ is irreducible, $W$ is connected in the complex analytic topology. Therefore, $W$ itself is a single $\mathbb{T}$-leaf.
    
    Let us return to our original setup. Assume $(V,\pi_V)$ has an open $\mathbb{T}$-leaf. Then $(U,\pi_U)$ has an open $\mathbb{T}$-leaf given by the intersection of $U$ with the open $\mathbb{T}$-leaf of $(V,\pi_V)$. Since $F|_U$ is Poisson and $\mathbb{T}$-equivariant, $(U^\prime,\pi_{U^\prime})$ also has an open $\mathbb{T}$-leaf, which is the image of the open $\mathbb{T}$-leaf of $(U,\pi_U)$ under $F|_U$. Choose a smooth point $p$ in the open $\mathbb{T}$-leaf of $(U^\prime,\pi_{U^\prime})$, and set $d:=\dim V$. Let $2r$ be the rank of $\pi_{(V^\prime)^{\sm}}$. Then one can choose linearly independent elements $\xi_1,\cdots, \xi_{d-2r}\in \mathfrak{t}$ such that
    \begin{equation*}
        \Imm\, \pi^\sharp_{V^\prime,p} +\sum_{j=1}^{d-2r}\mathbb{C}\chi(\xi_j)_{p} =  T_p V^\prime.
    \end{equation*}
    Since the following section is nonvanishing at $p$, it defines a $\mathbb{T}$-Pfaffian of $(V^\prime,\pi_{V^\prime})$:
    \[
        \mathrm{Pf}_\mathbb{T}(\pi_{V^\prime})=(\wedge^{r} \pi_{V^\prime}) \wedge \chi({\xi_1}) \wedge \cdots \wedge \chi(\xi_{d-2r}).
    \]
    Hence, $(V^\prime)^{\sm} \backslash \{\mathrm{Pf}_\mathbb{T}(\pi_{V^\prime}) = 0\}$ is an open $\mathbb{T}$-leaf of $(V^\prime,\pi_{V^\prime})$. Applying the same argument to the birational
    inverse $F^{-1}$, we obtain the converse implication.
\end{proof}
\begin{corollary}
    Assume that $(\bg^r,\bg^c)$ is aperiodic. Then the Poisson variety $(\Mat_{n \times m},\pi_{(\bbg^r,\bbg^c)})$ has an open $\mathbb{T}_{(\bg^r,\bg^c)}$-leaf.
\end{corollary}
\begin{proof}
    The corollary is an application of Proposition~\ref{p:open_t_leaf} to the $\mathbb{T}_{(\bg^r,\bg^c)}$-equivariant Poisson birational map
    \begin{equation*}
        \mathcal{U}_0 : \big( \Mat_{n\times m}(\mathbb{C}), (\pi_{(\bbg^r_{\std},\bbg^c_{\std})})_0 \big) \dashrightarrow \big( \Mat_{n\times m}(\mathbb{C}), \pi_{(\bbg^r,\bbg^c)} \big).
    \end{equation*}
    This map is the composition of~\eqref{eq:Ubirat} with~\eqref{eq:U0std}. 
    The properties of the composition are established in Proposition~\ref{p:U0}, Proposition~\ref{p:Tequivar}, Proposition~\ref{p:birationality}, and Proposition~\ref{p:u_poiss}.
\end{proof}

\appendix
\section{Gauss decomposition}\label{s:appendix_gauss}
In this appendix, we gather some additional background on the Gauss decomposition. 

Let $G$ be a reductive complex algebraic group with a fixed pinning as in Section~\ref{s:bd_background}, and let $r$ be its semisimple rank. Denote by $u_{i} $ (resp. $u_{-i}$) the one-parameter unipotent subgroups corresponding to a simple root $\alpha_i$, given by
\begin{equation*}
    u_{i}(t) = \exp(t e_{\alpha_i}), \quad u_{-i}(t) = \exp(te_{-\alpha_i})
\end{equation*}
where $e_{\alpha_i}$ (resp. $e_{-\alpha_i}$) are the corresponding Chevalley generators. Let us recall that the cocharacter $t \mapsto t^{\alpha_i^{\vee}}$ is uniquely determined by the requirement that for any weight $\lambda$, $(t^{\alpha_i^{\vee}})^{\lambda} = t^{\lambda(\alpha_i^{\vee})}$. The $1$-parameter subgroups satisfy the following well-known commutation relations (see, for instance, \cite[Section 2]{fz-dbl-bruhat}):
\begin{equation}\label{eq:comm_rel}
    u_i(t_1) u_{-i}(t_2) = u_{-i}\left( \frac{t_2}{1+t_1t_2}\right) (1+t_1 t_2)^{\alpha_i^\vee} u_i\left( \frac{t_1}{1+t_1 t_2}\right), \quad (1+t_1 t_2) \neq 0.
\end{equation}

Introduce the additive characters $\beta_i : \mathcal{N} \rightarrow \mathbb{C}$ and $\beta_{-i} : \mathcal{N}_- \rightarrow \mathbb{C}$ given by
\begin{equation*}
    \beta_i(u_j(t)) = \delta_{ij} t, \quad \beta_{-i}(u_{-j}(t)) = \delta_{ij} t.
\end{equation*}
Recall that for a generic element $x \in G$, we denote by $[x]_- \in \mathcal{N}_-$, $[x]_0 \in \mathcal{H}$, and $[x]_+ \in \mathcal{N}$ its Gauss components, so that $x = [x]_-[x]_0 [x]_+$. Denote by $G_0:=\mathcal{B}_-\mathcal{B}$ the set of elements $x \in G$ admitting a Gauss decomposition.

\begin{proposition}\label{p:gauss_perturb}
    Let $x \in G_0$. The following statements hold:
    \begin{enumerate}[1)]
        \item For any $i \in [1,r]$ and $t \in \mathbb{C}$ such that $1+t\beta_{-i}([x]_-) \neq 0$,
        \begin{equation*}
            [u_{i}(t)x]_+ = u_{i}\Big( \frac{t}{[x]_0^{\alpha_i}(1+t\beta_{-i}([x]_-))}\Big) [x]_+.
        \end{equation*}
        \item For any $i \in [1,r]$ and $t \in \mathbb{C}$ such that $1+\beta_i([x]_+)t \neq 0$,
        \begin{equation*}
            [xu_{-i}(t)]_- = [x]_-u_{-i}\Big( \frac{t}{[x]_0^{\alpha_i}(1+\beta_i([x]_+)t)} \Big).
        \end{equation*}
    \end{enumerate}
\end{proposition}
\begin{proof}
    The two equations can be proved in a similar way, so let us focus on the second one. Write $[x]_+$ as a product of one-parameter unipotent subgroups corresponding to simple roots
    \begin{equation*}
        [x]_+ = u_{i_1}(t_1) u_{i_2}(t_2) \cdots u_{i_s}(t_s), \quad t_i \in \mathbb{C},\; s \in \mathbb{N}.
    \end{equation*}
    Using the commutation relation~\eqref{eq:comm_rel}, we see that there is $b \in \mathcal{B}$ such that
    \begin{equation*}
        [x]_{+}u_{-i}(t) = u_{-i}\Big( \frac{t}{1+\beta_i([x]_+)t} \Big)b.
    \end{equation*}
    Commuting $u_{-i}$ past $[x]_0$, we see that
    \begin{equation*}
        [x]_0u_{-i}\Big( \frac{t}{1+\beta_i([x]_+)t} \Big) = u_{-i}\Big( \frac{t}{[x]_0^{\alpha_i}(1+\beta_i([x]_+)t)} \Big)[x]_0.
    \end{equation*}
    Thus the equation holds.
\end{proof}

        \end{document}